\documentclass[12pt]{article}

\usepackage{vmargin,epsfig}
\usepackage[francais]{babel}
\usepackage[T1]{fontenc}
\usepackage[tbtags]{amsmath}
\usepackage{amsthm,amssymb}
\usepackage{pifont}

\newtheorem{introconj}{Conjecture}[section]
\newtheorem{introtheo}[introconj]{Théorème}

\newtheorem{deftn}{Définition}[subsection]

\newtheorem{theo}[deftn]{Théorème}

\newtheorem{prop}[deftn]{Proposition}

\newtheorem{lemme}[deftn]{Lemme}

\newtheorem{cor}[deftn]{Corollaire}

\newenvironment{preuve}{
  \noindent \textbf{Démonstration.}}{\hfill \(\square\) }

\def\leq{\leqslant}
\def\geq{\geqslant}

\def\N{\mathbb{N}}
\def\Z{\mathbb{Z}}
\def\Q{\mathbb{Q}}

\def\C{\mathbb{C}}
\def\F{\mathbb{F}}
\def\O{\mathcal{O}}

\def\S{\mathcal{S}}

\def\pa#1{\left(#1\right)}
\def\cro#1{\left[ #1 \right]}
\def\ccro#1{\left[ \! \left[ #1 \right] \! \right]}
\def\acco#1{\left\{ #1 \right\}}
\def\brac#1{\left< #1 \right>}

\def\epsilon{\varepsilon}

\def\hom{{\text{Hom}}}

\def\id{{\text{id}}}
\def\ker{\text{ker}\,}

\def\im{\text{im}\,}

\def\Fil{\text{Fil}\,}

\def\tr{\text{Tr}}

\def\M{\mathcal{M}}
\def\Mr{\underline{\mathcal{M}}^r}
\def\MO{\underline{\mathcal{M}}^0}
\def\Mrtilde{\underline{\widetilde {\mathcal{M}}}^r}
\def\MOtilde{\underline{\widetilde {\mathcal{M}}}^0}

\def\MSK0{\underline{\mathcal{M}}^{S_{K_0}}}

\def\pr{\text{pr}}

\def\Tst{T_{\text{st}}}
\def\Acris{A_{\text{cris}}}

\def\Ast{\hat A_{\text{st}}}

\def\Astinf{\hat A_{\text{st}, \infty}}

\def\upi{\underline \pi}

\def\It{I_t}
\def\Is{I_s}

\def\nr{\text{nr}}
\def\mr{\text{mr}}
\def\et{\text{ét}}

\def\Phi{\phi}

\def\calB{\mathcal{B}}
\def\calM{\mathcal{M}}
\def\calN{\mathcal{N}}

\def\calC{\mathcal{C}}
\def\calF{\mathcal{F}}

\def\cris{\text{cris}}
\def\CRIS{\text{CRIS}}
\def\syncris{\text{syn-cris}}
\def\SYNCRIS{\text{SYN-CRIS}}
\def\logcris{\text{log-cris}}
\def\et{\text{ét}}
\def\ET{\text{\'ET}}
\def\syn{\text{syn}}
\def\SYN{\text{SYN}}
\def\ss{\text{ss}}
\def\DP{\text{DP}}
\def\gp{\text{gp}}
\def\spec{\text{Spec\,}}
\def\st{\text{st}}
\def\car{\text{car}}
\def\log{\text{log}}

\def\J{\mathcal{J}}

\def\arinj{\ar@{^(->}}
\def\toinj{\hookrightarrow}

\def\sym{\text{Sym}}

\def\tors{\text{tors}}
\def\free{\text{free}}
\def\triv{\text{triv}}

\def\xyhat#1{\smash{\hat{#1}}\vphantom{#1}}
\def\xytilde#1{\smash{\tilde{#1}}\vphantom{#1}}
\def\xybar#1{\smash{\bar{#1}}\vphantom{#1}}
\def\xyAst{\xyhat A_{\text{st}}}
\def\xyvarinjlim#1{\smash{\varinjlim_{#1}}\vphantom{\varinjlim}}
\def\xybigoplus#1#2{\smash{\bigoplus_{#1}^{#2}}\vphantom{\bigoplus}}

\setmarginsrb{2cm}{1.5cm}{2cm}{1.5cm}{0.5cm}{0.5cm}{0.5cm}{1.2cm}
\include{xy}
\xyoption{all}
\makeatletter
\renewcommand{\everyentry@}{\vphantom{A_{[]}^{[]}}}
\makeatother

\title{Conjecture de l'inertie modérée de Serre}
\author{Xavier Caruso}
\date{Septembre 2005}

\begin{document}

\maketitle

\renewcommand{\abstractname}{Abstract}
\begin{abstract}
Let $K$ be a local field of mixte characteristics. We assume that the 
residue field is perfect. Let $X_K$ be a proper smooth scheme over $K$ 
admitting an integer model $X$ which is proper and semi-stable. In this
article, we prove a period isomorphism linking the étale cohomology of
$X_{\bar K}$ with coefficients in $\Z/p^n\Z$ and the log-crystalline 
cohomology of the special fiber of $X$. Nevertheless, we have a 
restriction on the absolute ramification of $K$ and the degree of the 
cohomologies.

We apply the theory to deduce a complete proof of the Serre conjecture
on the tame inertia.
\end{abstract}

\medskip

\renewcommand{\abstractname}{Résumé}
\begin{abstract}
On considère $K$ un corps complet pour une valuation discrète, de 
caractéristique nulle et dont le corps résiduel est supposé parfait de 
caractéristique $p$. On appelle $\O_K$ l'anneau des entiers de $K$, et 
$\bar K$ une clôture algébrique. Soit $X_K$ un schéma propre et 
lisse sur $K$ admettant un modèle propre et semi-stable $X$ sur $\O_K$. 
Dans cet 
article, on démontre un isomorphisme de périodes reliant le $r$-ième 
groupe de cohomologie étale de $X_{\bar K}$ 
à coefficients dans $\Z/p^n\Z$ et un $r$-ième groupe 
de cohomologie log-cristalline de la fibre spéciale de $X$. Nous 
avons toutefois la restriction $er < p-1$ (et une restriction 
légèrement plus forte si $n>1$) où $e$ désigne l'indice de 
ramification absolu de $K$.

On en déduit une preuve complète de la conjecture de Serre sur l'inertie 
modérée (voir \cite{serre}).
\end{abstract}

\tableofcontents

\section{Introduction}

Tout au long de cet article, on considère $p$ un nombre premier et $k$ 
un corps parfait de caractéristique $p$. On note $W = W(k)$ l'anneau des 
vecteurs de Witt à coefficients dans $k$ et $K_0$ son corps des 
fractions. On note $\sigma$ le Frobenius sur $k$, sur $W$ et sur 
$K_0$. On considère $K$ une extension totalement ramifiée de $K_0$ 
de degré $e$. On fixe $\pi$ une uniformisante de $K$ et on note $E\pa u$ 
son polynôme minimal sur $K_0$. Il s'agit d'un polynôme d'Eisenstein.
On note de plus $\O_K$ l'anneau des entiers de $K$. Le corps résiduel 
$\O_K / \pi$ s'identifie à $k$.

\medskip

On fixe $\bar K$ (resp. $\bar k$) une clôture algébrique de $K$ (resp. 
de $k$) et 
on définit $G_K$ (resp. $G_k$) comme le groupe de Galois absolu de $K$ 
(resp. de $k$). On désigne par $K^\nr$ (resp. $K^\mr$) l'extension 
maximale non ramifiée (resp. modérément ramifiée) de $K$ et $I$ comme
(resp. $\Is$) le groupe d'inertie (resp. d'inertie sauvage), 
c'est-à-dire le groupe de Galois de $\bar K$ sur $K^\nr$ (resp. sur 
$K^\mr$). Le quotient $\It = I/\Is$ est le groupe d'inertie modérée.

\medskip

Le but de cet article est de comparer, lorsque $X$ est un schéma 
propre et semi-stable sur $\O_K$, la cohomologie log-cristalline 
(définie par 
Kato --- voir \cite{kato}) de la fibre spéciale de $X$ et la cohomologie 
étale de $X_{\bar K} = X \times_{\O_K} \bar K$. Ces théorèmes de 
comparaison s'inscrivent dans une grande lignée amorcée par 
Grothendieck, Tate et Raynaud et poursuivie par Fontaine, Messing, 
Faltings, Kato, Tsuji, Breuil... Nous nous intéressons particulièrement 
aux cas des coefficients de torsion. Précisément nous obtenons le 
théorème :

\begin{introtheo}
\label{4:th:comp}
On garde les notations précédentes et on fixe $r$ un entier vérifiant 
$er < p-1$. Si $X_n = X \times_{\O_K} \O_K/p^n$, on a un isomorphisme 
canonique de modules galoisiens :
$$H^i_\et(X_{\bar K}, \Z/p^n\Z)(r) = {\Tst}_\star 
(H^i_\logcris(X_n/(S/p^n S)))$$
pour tout $i < r$ (et aussi $i = r$ si $n=1$).
\end{introtheo}

\noindent
Dans ce théorème $S$ et ${\Tst}_\star$ désignent respectivement une 
certaine 
$W$-algèbre, et un certain foncteur d'une catégorie de $S$-modules de 
torsion $\Mr$ vers la catégorie des $\Z_p$-représentations galoisiennes, 
tous deux introduits par Breuil dans \cite{breuil-ens} (pour le cas 
$e=1$) et \cite{breuil-invent} (pour le cas général), et étudiés dans
\cite{caruso}.

\medskip

Comme conséquence du théorème \ref{4:th:comp} et des résultats de
\cite{caruso}, nous donnons une réponse affirmative à une 
question formulée par Serre dans le paragraphe 1.13 de \cite{serre} :

\begin{introtheo}
\label{4:th:serre}
On garde les notations précédentes et on fixe $r$ un entier quelconque. 
Soient $V$ la restriction au groupe d'inertie $I$ de la 
$\F_p$-représentation $H^r_\et(X_{\bar K}, \Z/p\Z)^\vee$ (où \og 
$^\vee$ \fg\ signifie que l'on prend le $\F_p$-dual) et $V^\ss$ la 
semi-simplifiée de $V$. Alors les poids de l'inertie modérée sur $V^\ss$ 
sont tous compris entre $0$ et $er$.
\end{introtheo}

\bigskip

Cet article s'articule comme suit. Les deux chapitres qui suivent cette 
introduction se bornent à présenter les principaux objets : le chapitre 
\ref{4:sec:alglin} est consacré aux objets 
d'algèbre linéaire (l'anneau $S$, les catégories $\Mr$ et les foncteurs 
$\Tst^\star$ et ${\Tst}_\star$), alors que le chapitre \ref{4:sec:synto} 
plus géométrique introduit le site log-syntomique et les faisceaux 
$\O^\st_n$ et $\J_n^{[s]}$, ainsi que certaines variantes, qui 
s'avèreront cruciaux pour la preuve du théorème \ref{4:th:comp}.

Dans le chapitre \ref{4:sec:cris}, on prouve que le groupe de 
cohomologie $H^r_\logcris(X_n/(S/p^n S))$ peut être muni de structures 
supplémentaires qui en font un objet de la catégorie $\Mr$ (auquel 
on peut alors appliquer le foncteur ${\Tst}_\star$). La démonstration se 
découpe en 
deux parties : en premier lieu, on montre le résultat lorsque $n=1$, 
puis on l'étend à tout $n$ par un dévissage.

Finalement, dans le chapitre \ref{4:sec:etale}, on étudie la cohomologie 
étale et on prouve le théorème \ref{4:th:comp}, d'où nous déduisons 
directement le théorème \ref{4:th:serre}.

\bigskip

Ce travail a été accompli dans le cadre de ma thèse de doctorat
en mathématique que je prépare sous la direction de Christophe
Breuil. Je tiens à le remercier vivement ici pour les conseils, les
explications et les réponses qu'il a toujours su me fournir, ainsi que
pour la relecture patiente des versions préliminaires de ce texte.

\section{Les objets d'algèbre linéaire}
\label{4:sec:alglin}

\subsection{La catégorie $\M^r$}
On reprend les notations du début de l'introduction et on fixe dans tout 
ce chapitre un entier $r$ positif ou nul vérifiant l'inégalité $er < 
p-1$. On rappelle que $e$ désigne le degré de l'extension $K/K_0$, 
c'est-à-dire l'indice de ramification absolue de $K$.

\subsubsection{L'anneau $S$}
\label{4:subsec:s}

Soit $W\cro u$ l'anneau des polynômes en une indéterminée $u$ à 
coefficients dans $W$. Par définition, $S$ est le complété $p$-adique
de l'enveloppe aux puissances divisées de $W\cro u$ par rapport à 
l'idéal principal engendré par $E\pa u$ compatibles aux puissances 
divisées canoniques sur $p W\cro u$. Concrètement, $S$ est la 
sous-$W$-algèbre de $K_0\ccro u$ suivante :
$$S = \acco{\sum_{i=0}^{\infty} w_i \frac{\pa{E\pa u}^i}{i!}, \, w_i 
\in W\cro u, \, \lim_{i\to \infty} w_i = 0}$$
ou encore :
$$S = \acco{\sum_{i=0}^{\infty} w_i \frac{u^i}{q\pa i!}, \, w_i 
\in W, \, \lim_{i\to \infty} w_i = 0}$$
où $q\pa i$ désigne le quotient de la division euclidienne de $i$ par 
$e$.

\medskip

On munit $S$ d'un Frobenius $\phi$ défini comme l'unique application
$\sigma$-semi-linéaire vérifiant $\phi(u^i/q(i)!) = u^{pi}/q(i)!$ et 
d'un opérateur de monodromie $N$ défini comme l'unique application 
$W$-linéaire vérifiant $N(u^i/q(i)!) = - i u^i/q(i)!$. On munit 
également $S$ d'une filtration : pour tout entier positif ou nul $n$, on 
définit $\Fil^n S$ comme le complété $p$-adique de l'idéal engendré par 
les éléments $\frac {\pa{E \pa u}^i}{i!}$ pour $i \geq n$. On a
$\Fil^0 S = S$, $\Fil^n S \subset \Fil^{n-1} S$, $\bigcap_{n \in \N} 
\Fil^n S = 0$, et certaines compatibilités vis-à-vis des opérateurs à 
savoir $N\pa{\Fil^n S} \subset \Fil^{n-1} S$ et, pour $0 \leq n \leq 
p-1$, $\phi\pa{\Fil^n S} \subset p^n S$. Cela permet de définir, pour 
$0 \leq n \leq p-1$, l'application $\phi_n = \frac \phi {p^n} : \Fil^n S 
\to S$. L'élément $\phi_1\pa{E\pa u}$ est une unité de $S$, on le 
notera $c$ par la suite.

\medskip

On note $S_n = S/p^n S$. Le Frobenius, l'opérateur de monodromie et la 
filtration passent au quotient et définissent des structures analogues 
sur $S_n$.

\subsubsection{Définition des catégories}

Par définition, un objet de la catégorie $\Mr$ est la donnée :
\begin{enumerate}
\item d'un $S$-module $\calM$ isomorphe à une somme directe (finie) de 
$S_n$ pour des entiers $n$ convenables ;
\item d'un sous-module $\Fil^r \calM$ de $\calM$ contenant $\Fil^r S 
\cdot \calM$ ;
\item d'une flèche $\phi$-semi-linéaire $\phi_r : \Fil^r \calM \to \calM$ 
vérifiant la condition :
$$\phi_r\pa {sx} = \frac 1 {c^r} \phi_r \pa{s} \phi_r\pa{\pa{E\pa u}^r 
x}$$
pour tout élément $s \in \Fil^r S$ et tout élément $x \in \calM$ et 
telle que $\im \phi_r$ engendre $\calM$ en tant que $S$-module ;
\item d'une application $W$-linéaire $N : \M \to \M$ telle que :
\begin{itemize}
\item pour tout $s \in S$ et tout $x \in \calM$, $N\pa{sx} = N\pa s x + s 
N\pa x$
\item $E\pa u N\pa{\Fil^r \calM} \subset \Fil^r \calM$
\item le diagramme suivant commute :
\[\xymatrix @C=50pt {
\Fil^r \calM \ar[r]^{\phi_r} \ar[d]_{E\pa u N} & \calM \ar[d]^{cN} \\
\Fil^r \calM \ar[r]^{\phi_r} & \calM}\]
\end{itemize}
\end{enumerate}

\medskip

Une flèche entre deux objets $\calM$ et $\calM'$ de cette catégorie est 
un morphisme $S$-linéaire de $\calM$ dans $\calM'$ respectant la filtration et 
commutant aux applications $\phi_r$ et $N$.

\bigskip

Nous renvoyons à \cite{caruso} pour l'étude de la catégorie $\Mr$. Il y 
est prouvé en particulier que $\Mr$ est une catégorie abélienne et 
artinienne.

\subsection{Les objets tués par $p$}
\label{4:sec:mrtilde}
Dans ce paragraphe, nous nous intéressons à la sous-catégorie pleine 
de $\Mr$ formée des objets tués par $p$. Cette dernière catégorie est
équivalente à une catégorie d'objets sur $k\cro u / u^{ep}$ plus simple 
à manipuler que les objets de $\Mr$ (voir \cite{caruso}). Cependant, 
nous aurons besoin d'une description encore différente utilisant des 
objets sur $k\cro u / u^p$, et c'est celle-ci que nous allons détailler 
dans ce paragraphe.

\medskip

On commence par rappeler le résultat suivant :

\begin{lemme}
Soit $\calM$ un objet de $\Mr$ tué par $p$. Alors l'application :
$$\xymatrix @C=50pt {
S_1 \otimes_{k\cro u / u^e} \Fil^r \calM / E(u) \Fil^r \calM 
\ar[r]^-{\id \otimes \phi_r} & \cal M }$$
est un isomorphisme (où $S$ est vu comme un $k\cro u / u^e$-module 
\emph{via} le Frobenius $\phi : u^i \mapsto u^{pi}$).
\end{lemme}

\begin{preuve}
Le cas général est similaire au cas $e=1$ traité dans 
\cite{breuil-ens} (corollaire 2.2.2.2).
\end{preuve}

\bigskip

On déduit directement de ce lemme le corollaire suivant :

\begin{cor}
\label{4:cor:phir}
Soient $\calM$ et $\calN$ deux objets de $\Mr$ tués par $p$. Soit $f : 
\calM \to \calN$ une application $S_1$-linéaire telle que $f(\Fil^r 
\calM) \subset \Fil^r \calN$. Alors il existe une unique application 
$S_1$-linéaire $g : \calM \to \calN$ faisant commuter le diagramme 
suivant :
\[\xymatrix @C=50pt {
\Fil^r \calM \ar[r]^f \ar[d]_{\phi_r} & \Fil^r \calN
\ar[d]^{\phi_r} \\
\calM \ar[r]^g & \calN }\]
\end{cor}

\bigskip

Soit $\tilde S = k\cro u / u^p$. On munit $\tilde S$ d'un Frobenius
$\phi$, unique application $\sigma$-semi-linéaire vérifiant $\phi(u^i) = 
u^{ip}$, et d'un opérateur de monodromie $N$, unique application 
$k$-linéaire vérifiant $N(u^i) = -i u^i$. On définit également une 
filtration sur $\tilde S$ en posant pour tout entier $n$, $\Fil^n 
\tilde S = u^{en} \tilde S$. On dispose d'un morphisme d'anneaux $S_1 
\to \tilde S$ qui envoie $u$ sur $u$ et toutes les puissances divisées 
$u^i/q(i)!$ sur $0$ pour $i \geq p$.

On définit la catégorie $\Mrtilde$ en adaptant la définition de 
la catégorie $\Mr$. Un objet de $\Mrtilde$ est la donnée suivante :

\begin{enumerate}
\item un $\tilde S$-module $\tilde \calM$ libre de rang fini ;
\item un sous-module $\Fil^r \tilde \calM$ de $\tilde \calM$ contenant 
$\Fil^r \tilde S \cdot \tilde \calM = u^{er} \tilde \calM$ ;
\item une flèche $\phi$-semi-linéaire $\phi_r : \Fil^r \tilde \calM \to 
\tilde \calM$ telle que l'image de $\phi_r$ engendre $\tilde \calM$ en 
tant que $\tilde S$-module ;
\item une application $k$-linéaire $N : \tilde \calM \to \tilde \calM$ 
telle que :
\begin{itemize}
\item pour tout $\lambda \in \tilde S$ et tout $x \in \tilde \calM$, 
$N\pa{\lambda x} = N\pa \lambda x + \lambda N\pa x$
\item $u^e N(\Fil^r \tilde \calM) \subset \Fil^r \tilde \calM$
\item le diagramme suivant commute :
\[\xymatrix @C=50pt {
\Fil^r \xytilde \calM \ar[r]^{\phi_r} \ar[d]_{u^e N} & \xytilde \calM 
\ar[d]^{c_\pi N} \\
\Fil^r \xytilde \calM \ar[r]^{\phi_r} & \xytilde \calM}\]
où $c_\pi$ est la réduction de $c$ dans $\tilde S$.
\end{itemize}
\end{enumerate}

\noindent
Les morphismes entre deux objets de $\Mrtilde$ sont les applications 
$\tilde S$-linéaires qui respectent le $\Fil^r$ et commutent au 
Frobenius et à l'opérateur de monodromie. 

\medskip

On dispose pour les objets de $\Mrtilde$ de la proposition suivante, 
fort utile pour les manipulations :

\begin{prop}
Soit $\tilde \calM$ un objet de $\Mrtilde$. Alors il existe $\pa{e_1,
\ldots, e_d}$ une $\tilde S$-base de $\tilde \calM$ et des entiers 
$n_1, \ldots, n_d$ compris entre $0$ et $er$ tels que :
$$\Fil^r \tilde \calM = \bigoplus_{i=1}^d u^{n_i} k\cro u / u^p \cdot 
e_i.$$
Une telle famille $\pa{e_1, \ldots, e_d}$ est appelée une \emph{base
adaptée} de $\tilde \calM$.
\end{prop}

\begin{preuve}
C'est une conséquence directe du théorème de structure des modules de 
type fini sur un anneau principal (l'anneau étant ici $k\cro u$).
\end{preuve}

\bigskip

Par ailleurs, on dispose d'un foncteur $T$ de la sous-catégorie pleine 
de $\Mr$ formée des objets tués par $p$ dans la catégorie $\Mrtilde$ 
défini de la façon suivante. Soit $\calM$ un objet de $\Mr$ tué par $p$. 
C'est en particulier un $S_1$-module libre de rang fini et le produit 
tensoriel $\tilde \calM = \calM \otimes_{S_1} \tilde S$ est un $\tilde 
S$-module libre de rang fini. On dispose d'une projection canonique
$\calM \to \tilde \calM$. On définit $\Fil^r \tilde\calM$ comme l'image
de $\Fil^r \calM$ par cette projection, et on vérifie facilement que
les opérateurs $\phi_r$ et $N$ définis sur $\calM$ passent au quotient
pour fournir respectivement des opérateurs $\Fil^r \tilde \calM \to 
\tilde \calM$ et $\tilde \calM \to \tilde \calM$ encore notés $\phi_r$ 
et $N$.

\bigskip

Notons $\kappa = \ker(S_1 \to \tilde S)$ et prouvons un lemme concernant 
ce noyau :

\begin{lemme}
\label{4:lem:kappa}
On suppose $r \neq 0$, c'est-à-dire $1 \leq er \leq p-2$. Alors $\kappa 
\subset \Fil^r S_1$, $\phi_r(\kappa) \subset \Fil^r S_1$ et $\phi_r 
\circ \phi_r (\kappa) = 0$.
\end{lemme}

\begin{preuve}
On constate que $\kappa$ est l'idéal engendré par les $u^i/q(i)!$ pour 
$i \geq p$. En particulier, on a bien $\kappa \subset \Fil^r S_1$. Par 
ailleurs, puisque $r > 0$, on a, pour $i \geq p$ :
\begin{equation}
\label{4:eq:phirkappa}
\phi_r\pa{\frac{u^i}{q(i)!}} =  \phi_r\pa{u^{er} 
\frac{u^{i-er}}{q(i)!}} = \phi_r\pa{u^{er}} 
\phi\pa{\frac{u^{i-er}}{q(i)!}} = c^r \frac{u^{p(i-er)}}{q(i)!}.
\end{equation}
Comme $er \leq p-2$, on a $p(i-er) \geq 2p$ et donc 
$\phi_r(\kappa)$ est inclus dans $\kappa'$, l'idéal engendré par les 
$u^i/q(i)!$ pour $i \geq 2p$. En particulier, $\phi_r(\kappa) \subset 
\Fil^r S_1$. 

D'autre part, si $i \geq 2p$, on a $p(i-er) > i + ep$ et 
$v_p((p(i-er))!) > v_p(q(i)!)$ d'où par la formule (\ref{4:eq:phirkappa})
$\phi_r(u^i/q(i)!) = 0$ (dans $S_1$). Ainsi $\phi_r(\kappa') = 0$, 
d'où $\phi_r \circ \phi_r (\kappa) = 0$.
\end{preuve}

\bigskip

On a alors la proposition suivante qui permet de réduire l'étude de la 
catégorie $\Mr$ à celle de la catégorie $\Mrtilde$ et à des dévissages.

\begin{prop}
\label{4:prop:tequiv}
Le foncteur $T$ est une équivalence de catégories
\end{prop}

\begin{preuve}
Il faut traiter à part le cas trivial $r=0$ qu'on laisse au lecteur. On 
suppose à partir de maintenant $1 \leq er \leq p-2$.

\medskip

Prouvons la pleine fidélité du foncteur $T$. Soient $\calM$ et $\calM'$ 
deux objets de $\Mr$ tués par $p$ et $f : \calM \to \calN$ un morphisme 
entre ces objets. On suppose que $f$ vaut $0$ dans $\Mrtilde$, 
c'est-à-dire que $f(\calM) \subset \kappa \calN$. En particulier, 
$f(\Fil^r \calM) \subset \kappa \calN$ et donc :
$$f\circ \phi_r  (\Fil^r \calM) = \phi_r \circ f (\Fil^r 
\calM) \subset \phi_r(\kappa) \calN.$$
Par hypothèse $\phi_r \pa{\Fil^r \calM}$ engendre tout $\calM$, et donc  
$f(\calM) \subset \phi_r(\kappa) \calN$. En réitérant l'argument, et 
puisque $\phi_r \circ \phi_r(\kappa) = 0$ (lemme \ref{4:lem:kappa}), on 
obtient $f(\calM) = 0$ et donc $f=0$, ce qui assure la fidélité du 
foncteur.

Considérons à présent $\calM$ et $\calN$ deux objets de $\Mr$ tués par 
$p$. Notons $\tilde \calM$ et $\tilde \calN$ leurs images respectives 
dans la catégorie $\Mrtilde$. Soit $\tilde f : \tilde \calM \to \tilde 
\calN$ un morphisme de la catégorie $\Mrtilde$. On veut montrer qu'il 
existe un morphisme (nécessairement unique) de la catégorie $\Mr$, 
$f : \calM \to \calN$ tel que $f \equiv \tilde f \pmod {\kappa \calN}$.
On construit $f$ par approximations successives. On considère dans un
premier temps $f_0 : \calM \to \calN$ un relevé $S_1$-linéaire 
quelconque de $\tilde f$. Comme $\kappa \calN \subset \Fil^r \calN$, le 
relevé $f_0$ est automatiquement compatible à $\Fil^r$ et les 
applications $f_0$ et $\phi_r$ commutent modulo $\kappa \calN$.
Par le corollaire \ref{4:cor:phir}, il existe une unique application 
$S_1$-linéaire $f_1 : \calM \to \calN$ faisant commuter le diagramme 
suivant :
\[\xymatrix @C=50pt {
\Fil^r \calM \ar[r]^{f_0} \ar[d]_{\phi_r} & \Fil^r \calN
\ar[d]^{\phi_r} \\
\calM \ar[r]^{f_1} & \calN }\]
On vérifie directement que $f_0 \equiv f_1 \pmod {\kappa \calN}$, ce qui
implique d'une part que $f_1$ respecte les $\Fil^r$ et d'autre part, par 
un argument analogue à celui utilisé pour la fidélité, que $\phi_r \circ 
f_0 \equiv \phi_r \circ f_1 \pmod{\phi_r(\kappa) \calN}$. Ainsi $f_1$ et 
$\phi_r$ commutent modulo $\phi_r(\kappa) \calN$. On construit de même 
$f_2$ à partir de $f_1$, et celui-ci convient.

Il reste à prouver que $f_2$ commute automatiquement à $N$. Il s'agit à
nouveau d'un argument analogue. Pour plus de précisions, on pourra
consulter la fin de la preuve du lemme 6.1.3 de \cite{caruso}.

\medskip

Montrons pour finir l'essentielle surjectivité. Soit $\tilde \calM$ un 
objet 
de la catégorie $\Mrtilde$. Considérons $\pa{e_1, \ldots, e_d}$ une base 
adaptée de $\tilde \calM$ pour les entiers $n_1, \ldots, n_d$. Notons 
$\calM$ 
le $S_1$-module engendré par des éléments $\hat e_1, \ldots, \hat e_d$ 
et définissons :
$$\Fil^r \calM = \Fil^r S_1 \cdot \calM + \sum_{i=1}^d u^{n_i} S_1 
\cdot \hat e_i \subset \calM.$$
Soit $\pr : \calM \to \tilde \calM$ l'application $S_1$-linéaire
définie par $\pr(\hat e_i) = e_i$ pour tout $i$. Elle est surjective et 
respecte les $\Fil^r$. Pour tout $i$, notons $\hat x_i$ un relevé 
(\emph{i.e.} un antécédent par $\pr$) de $\phi_r \pa{e_i}$ et 
définissons $\phi_r(\hat e_i) = \hat x_i$. On prolonge $\phi_r$ à tout 
$\Fil^r \calM$ (de façon à respecter les conditions de la catégorie 
$\Mr$) obtenant ainsi une application $\phi_r : \Fil^r \calM \to 
\calM$, dont il est facile de vérifier que l'image engendre tout 
$\calM$. 

Il reste à définir un opérateur de monodromie sur $\calM$. Pour cela,
on procède à nouveau par approximations successives. On commence par
définir $N_0$ en imposant la condition de Leibniz et que $N_0(\hat e_i)$
soit un relevé de $N(e_i)$. On vérifie immédiatement que $E(u)
N_0(\Fil^r \calM) \subset \Fil^r \calM$ et que le diagramme suivant :
\[\xymatrix @C=50pt {
\Fil^r \calM \ar[r]^{\phi_r} \ar[d]_{E(u) N_0} & \calM \ar[d]^{c N_0}
\\
\Fil^r \calM \ar[r]^{\phi_r} & \calM}\]
commute modulo $\kappa \calM$. Une variante du corollaire \ref{4:cor:phir} 
assure que l'on peut construire une application $N_1$ vérifiant la 
condition de Leibniz et faisant commuter le diagramme suivant :
\[\xymatrix @C=50pt {
\Fil^r \calM \ar[r]^{\phi_r} \ar[d]_{E(u) N_0} & \calM \ar[d]^{c N_1}
\\
\Fil^r \calM \ar[r]^{\phi_r} & \calM}\]
Mais alors $N_1$ est un autre relevé de $N$, ce qui implique que $E(u) 
N_1(\Fil^r \calM) \subset \Fil^r \calM$. Et par ailleurs, le diagramme :
\[\xymatrix @C=50pt {
\Fil^r \calM \ar[r]^{\phi_r} \ar[d]_{E(u) N_1} & \calM \ar[d]^{c N_1}
\\
\Fil^r \calM \ar[r]^{\phi_r} & \calM}\]
commute modulo $\phi_r(\kappa) \calM$. L'application $N_2$ obtenue à 
partir de $N_1$ de la même façon que $N_1$ a été obtenue à partir de 
$N_0$ répond finalement à la question.
\end{preuve}

\bigskip

\noindent 
{\it Remarque.} Contrairement à ce qui se passe pour les objets modulo 
$u^{ep}$, il n'est, à notre connaissance, pas possible de décrire un 
quasi-inverse du foncteur $T$ par une simple formule.

\subsection{Foncteurs vers les représentations galoisiennes}
Il existe deux versions du foncteur vers les représentations 
galoisiennes. La première, que nous notons ${\Tst}_\star$, est 
covariante et la seconde, $\Tst^\star$, est contravariante. Nous sommes 
dans l'obligation de présenter ici les deux foncteurs et d'établir 
les liens qui les relient, car pour ce que nous voulons faire, il sera 
plus commode d'utiliser la version covariante, mais nous aurons 
également besoin d'utiliser les résultats de \cite{caruso} où seulement 
la version covariante est présentée et étudiée.

\subsubsection{Un anneau de périodes}
\label{4:subsec:acris}

Avant de pouvoir définir ces foncteurs, il faut introduire l'anneau de 
périodes $\Ast$. Cet anneau a une interprétation cohomologique que 
nous passons sous silence pour l'instant.

Pour tout entier $n$, on considère l'application :
$$\begin{array}{rcl}
\hat \theta_n \quad : \quad W_n (\O_{\bar K}/p) & \to & \O_{\bar K}/p^n 
\\
(a_0, a_1, \ldots, a_{n-1}) & \mapsto & \hat a_0^{p^n} + p \hat 
a_1^{p^{n-1}} + \cdots + p^{n-1} \hat a_{n-1}^p \\
\end{array}$$
où $\hat a_i \in \O_{\bar K}/p^n$ est un relevé quelconque de $a_i$. On 
note  $W_n (\O_{\bar K}/p)^\DP$ l'enveloppe aux puissances divisées de $ 
W_n (\O_{\bar K}/p)$ par rapport à $\ker \hat \theta_n$ (et 
compatibles avec les puissances divisées canoniques sur $p W_n (\O_{\bar 
K}/p)$). Les $W_n (\O_{\bar K}/p)^\DP$ forment un système projectif pour 
les applications de transition données par le Frobenius sur les vecteurs 
de Witt. On note $\Acris$ la limite projective de ce système.
On voit facilement que le Frobenius sur les vecteurs de Witt induit une 
application $\phi : \Acris \to \Acris$. En outre, on définit sur 
$\Acris$ une filtration obtenue à partir des filtrations données par
les puissances divisées sur $W_n (\O_{\bar K}/p)^\DP$. Si $t < 0$, on 
pose par convention $\Fil^t \Acris = \Acris$. Par ailleurs, $\Acris$ est 
muni d'une action du groupe de Galois $G_K$.

\medskip

Par définition $\Ast$ est le complété $p$-adique de $\Acris \brac X$. On 
munit $\Ast$ d'une filtration en posant :
$$\Fil^t \Ast = \acco{\sum_{i=0}^\infty a_i \frac{X^i}{i!}, \, \lim_{i 
\to \infty} a_i = 0, \, a_i \in \Fil^{t-i} \Acris}$$
pour tout entier $t$. On étend le Frobenius à $\Ast$ en imposant 
$\phi(X) = (1+X)^p - 1$. On vérifie que $\phi(\Fil^t \Ast) \subset 
p^t\Ast$, ce qui permet de définir une application $\phi_t = \phi/p^t : 
\Fil^t \Ast \to \Ast$. D'autre part, on définit sur $\Ast$ une 
dérivation $\Acris$-linéaire, par la formule :
$$N\pa{\frac{X^i}{i!}} = (1+X) \frac{X^{i-1}}{(i-1)!}.$$
On étend également l'action de $G_K$ à tout $\Ast$. Pour cela, on 
commence par fixer\footnote{Ainsi $\Ast$ dépend \emph{a priori} de ce 
choix. Cependant, on peut montrer qu'il n'en dépend pas à isomorphisme 
canonique près.} $\upi = (\pi_n)$ un système compatible de racines 
$p^n$-ièmes de 
$\pi$. Soit $g \in G_K$. On définit $\epsilon_n(g)$ comme l'unique 
élément de $\O_{\bar K}$ vérifiant $g(\pi_n) = \epsilon_n(g) \pi_n$. La 
famille des $(\epsilon_n(g))$ forme un système compatible de racines 
$p^n$-ièmes de l'unité et par suite un élément $[\epsilon(g)] \in 
\Acris$ obtenu à partir des représentants de Teichmüller 
$[\epsilon_n(g)] \in 
W_n (\O_{\bar K}/p)$. L'action de $g$ sur $X$ est alors donnée par la 
formule $g(X) = [\epsilon(g)] (1+X) - 1$. On étend cette action 
à tout $\Ast$ par semi-linéarité.

\medskip

L'anneau $S$ n'est pas sans rapport avec $\Ast$. Si l'on note $[\upi]$ 
l'élément de $\Ast$ défini à partir du système $(\pi_n)$ fixé 
précédemment, $\Ast$ peut être vu comme une $S$-algèbre \emph{via}
l'unique morphisme $W$-linéaire $S \to \Ast$ qui envoie $u$ sur 
$\frac{[\upi]}{1+X}$. Ce morphisme est injectif et identifie $S$ aux
invariants de $\Ast$ sous l'action du groupe $G_K$.

\subsubsection{La version contravariante}

Nous commençons par la version contravariante qui est plus simple à 
définir. Soit $\calM$ un objet de $\Mr$. On pose :
$$\Tst^\star(\calM) = \hom(\calM, \Astinf)$$
où par définition $\Astinf = \Ast \otimes_W K_0/W$ et où le $\hom$ 
précédent signifie que l'on ne considère que les morphismes 
$S$-linéaires, compatibles à $\Fil^r$, à $\phi_r$ et à $N$. L'objet
$\Tst^\star(\calM)$ est un $\Z_p$-module de torsion (de type fini)
qui hérite d'une action de $G_K$. On a donc ainsi bien défini un 
foncteur de $\Mr$ dans la catégorie des $\Z_p$-représentations (de 
torsion) du groupe $G_K$.

\medskip

Ce foncteur est étudié en détail dans \cite{breuil-ens} (pour le cas 
$e=1$) et dans \cite{caruso}. Le théorème suivant résume ses 
propriétés : 

\begin{theo}
\label{4:th:tstcont}
Le foncteur $\Tst^\star$ est exact, pleinement fidèle, d'image 
essentielle stable par quotients et par sous-objets. De plus, si $\calM$ 
est un objet de $\Mr$ isomorphe en tant que $S$-module à $S_{n_1} \oplus 
\cdots \oplus S_{n_d}$, alors la représentation galoisienne 
$\Tst^\star(\calM)$ est isomorphe en tant que $\Z_p$-module à 
$\Z/p^{n_1}\Z \times \cdots \times \Z/p^{n_d}\Z$.
\end{theo}

\bigskip

On dispose en outre d'une description plus simple du foncteur 
$\Tst^\star$ pour les objets tués par $p$. Considérons pour cela
l'anneau $\tilde A = \Ast/p \otimes_{S_1} \tilde S$. Comme dans 
le paragraphe 2.3 de \cite{caruso}, on montre 
que $\tilde A$ s'identifie à $(\O_{\bar K}/\pi) \brac X$. Il est 
possible de décrire les structures 
supplémentaires sur $(\O_{\bar K}/\pi) \brac X$. Exactement, $\Fil^t 
(\O_{\bar K}/\pi) \brac X$ est l'idéal engendré par les $\pi_1^{e(t-i)} 
\frac{X^i}{i!}$. La monodromie est l'unique opérateur $\O_{\bar 
K}/\pi$-linéaire envoyant $\frac{X^i}{i!}$ sur $(1+X) 
\frac{X^{i-1}}{(i-1)!}$. Il faut toutefois faire attention à $\phi_t$ 
car si $x \in \O_{\bar K}$ est un multiple de $\pi_1^{et}$ et si $\bar 
x$ désigne la réduction modulo $\pi$ de $x$, alors $\phi_t(\bar x)$ 
est la réduction modulo $\pi$ de $(-1)^t \frac{x^p}{p^t}$ (et pas 
$\frac{x^p}{p^t}$). On a ensuite la proposition suivante :

\begin{prop}
\label{4:prop:tstp}
Soit $\calM$ un objet de $\Mr$ tué par $p$. Alors :
$$\Tst^\star(\calM) = \hom(T(\calM), \tilde A)$$
où $\hom$ signifie que l'on considère les morphismes $\tilde 
S$-linéaires et commutant à $\Fil^r$, à $\phi_r$ et à $N$.
\end{prop}

\begin{preuve}
La tensorisation par $\tilde S$ au-dessus de $S_1$ fournit une 
application :
$$\Tst^\star(\calM) \to \hom(T(\calM), \tilde A).$$
On vérifie directement que cette application commute à l'action de 
Galois.

Soit $\psi \in \Tst^\star(\calM)$ qui s'envoie sur $0$ par l'application 
précédemment définie. On a alors un diagramme commutatif :
$$\xymatrix @C=50pt {
\calM \ar[r]^-{\psi} \ar[d] & \xyAst/p\xyAst \ar[d] \\
T\pa \calM \ar[r]^-{0} & \xytilde A }$$
où les flèches verticales sont déduites de la projection $S_1 \to \tilde 
S$. Ainsi, en reprenant les notations du lemme \ref{4:lem:kappa}, on 
a $\im \psi \subset \kappa \Ast/p\Ast$. Or $\psi$ commute par 
définition à $\phi_r$, d'où on déduit $\psi \circ \phi_r (\Fil^r \calM) 
= \phi_r \circ \psi(\Fil^r \calM) \subset \phi_r(\kappa) \Ast/p\Ast$. 
Comme par hypothèse, $\phi_r (\Fil^r \calM)$ engendre $\calM$, il vient 
$\im \psi \subset \phi_r(\kappa) \Ast/p\Ast$. En appliquant à nouveau 
l'argument, et en utilisant $\phi_r \circ \phi_r(\kappa) = 0$, on 
obtient $\im \psi = 0$ et donc $\psi = 0$. Ceci démontre l'injectivité 
de la flèche.

\medskip

Pour la surjectivité, on procède par approximations successives. Soit 
$\tilde \psi : T(\cal M) \to \tilde A$ un morphisme $\tilde S$-linéaire 
compatible aux structures. On note $\psi : \calM \to \Ast/p\Ast$ un 
morphisme $S_1$-linéaire faisant commuter le diagramme :
$$\xymatrix @C=50pt {
\calM \ar[r]^-{\psi} \ar[d] & \xyAst/p\xyAst \ar[d] \\
T\pa \calM \ar[r]^-{\xytilde \psi} & \xytilde A }$$
Dans un premier temps, on vérifie qu'automatiquement $\psi$ respecte le 
$\Fil^r$ et commute à $\phi_r$ modulo $\kappa \Ast/p\Ast$. D'après une 
variante du corollaire \ref{4:cor:phir}, il existe une unique application 
$S_1$-linéaire $\psi_1$ faisant commuter le diagramme :
$$\xymatrix @C=50pt {
\Fil^r \calM \ar[r]^-{\psi} \ar[d]_{\phi_r} & \Fil^r \xyAst/p\xyAst 
\ar[d]^{\phi_r} \\
\calM \ar[r]^-{\psi_1} & \xyAst/p\xyAst }$$
L'application $\psi_1$ respecte encore le $\Fil^r$ et commute à $\phi_r$ 
modulo $\phi_r(\kappa) \Ast/p\Ast$. De même, à partir de $\psi_1$, on 
construit $\psi_2$, qui respecte le $\Fil^r$ et commute à $\phi_r$ sans 
restriction.

Par un argument analogue (voir la fin de la preuve du lemme 
6.1.2 de \cite{caruso}), on montre que $\psi_2$ commute également à $N$.
\end{preuve}

\subsubsection{La version covariante}
\label{4:subsec:tstcov}

On commence par une définition, déjà présente dans \cite{breuil-duke} 
(définition 3.2.1.1) :

\begin{deftn}
\label{4:def:filadm}
Soit $\calM$ un objet de $\Mr$ (resp. de $\Mrtilde$). On appelle 
\emph{filtration admissible} de $\calM$ toute filtration décroissante 
$(Fil^t \calM)_{0 \leq t \leq r}$ par des sous-$S$-modules (resp. des 
sous-$\tilde S$-modules) vérifiant :
\begin{enumerate}
\item $\Fil^0 \calM = \calM$ et $\Fil^r \calM$ est \og le \fg\ $\Fil^r 
\calM$ de $\calM$ ;
\item pour tous $0 \leq t \leq t' \leq r$, $\Fil^{t'-t} S \cdot \Fil^t 
\calM \subset \Fil^{t'} \calM$ (resp. $\Fil^{t'-t} \tilde S \cdot \Fil^t 
\calM \subset \Fil^{t'} \calM$) ;
\item pour tout $1 \leq t \leq r$, $N(\Fil^t \calM) \subset \Fil^{t-1} 
\calM.$
\end{enumerate}
Si $(Fil^t \calM)_{0 \leq t \leq r}$ est une filtration admissible de 
$\calM$, on définit des opérateurs $\phi_t : \Fil^t\calM \to \calM$ par 
$\phi_t(x) = c^{t-r} \phi_t(E(u)^{r-t} x)$.
\end{deftn}

\medskip

Soit $\calM$ un objet de $\Mr$. On considère le produit tensoriel $\Ast 
\otimes_{S} \calM$. Il s'agit d'un $\Ast$-module de torsion 
naturellement muni d'une action de $G_K$ (en regardant son action sur le
premier facteur). On le munit en outre d'un opérateur de monodromie $N : 
\Ast \otimes_{S} \calM \to \Ast \otimes_{S} \calM$ en posant $N(a 
\otimes x) = N(a) \otimes x + a \otimes N(x)$.

On considère sur $\calM$ une filtration admissible 
quelconque\footnote{Il en existe toujours : on peut par 
exemple prendre $\Fil^t \calM = \acco{x \in \calM \, / \, E(u)^{r-t} x 
\in \Fil^r \calM}$.}. On peut alors définir, pour tout $s$ compris entre 
$0$ et $r$ :
$$\Fil^s(\Ast \otimes_{S} \calM) = \sum_{t=0}^s \Fil^t \Ast \otimes 
\Fil^{s-t} \calM.$$
C'est un sous-$\Ast$-module de $\Ast \otimes_{S} \calM$ qui dépend de la 
filtration admissible choisie. On définit finalement $\phi_s : 
\Fil^s(\Ast \otimes_{S} \calM) \to \Ast \otimes_{S} \calM$ comme 
l'unique application additive vérifiant $\phi_s(a_t \otimes x_t) = 
\phi_t(a_t) \otimes \phi_{s-t}(x_t)$ pour $a_t \in \Fil^t\Ast$ et $x_t 
\in \Fil^{s-t} \calM$.

\medskip

On pose finalement :
$${\Tst}_\star (\calM) = \Fil^r (\Ast \otimes_{S} \calM)_{N=0}^{\phi_r = 
1}$$
où la notation signifie que l'on ne retient que les $x \in \Fil^r (\Ast 
\otimes_{S} \calM)$ pour lesquels $N(x) = 0$ et $\phi_r(x) = x$. On 
obtient alors un $\Z_p$-module galoisien qui dépend \emph{a priori} du 
choix d'une filtration admissible. Toutefois, nous allons prouver dans
la suite que ce n'est pas le cas (voir la remarque faisant suite au 
corollaire \ref{4:cor:filadm}).

\bigskip

Encore une fois, ${\Tst}_\star (\calM)$ a une description plus simple 
lorsque $\calM$ est tué par $p$. Pour la donner, posons $\tilde \calM = 
T(\calM)$ et notons $\pr : \calM \to \tilde \calM$ la projection 
canonique. On vérifie facilement que si $(\Fil^t \calM)$ est une 
filtration admissible de $\calM$, alors $(\pr(\Fil^t \calM))$ est une
filtration admissible de $\tilde \calM$. En recopiant les constructions 
précédentes, on définit le $\Z_p$-module galoisien :
$$\Fil^r (\tilde A \otimes_{\tilde S} \tilde\calM)_{N=0}^{\phi_r = 1}$$
et on a alors la proposition suivante :

\begin{prop}
Si $\calM$ est un objet de $\Mr$ tué par $p$ alors :
$${\Tst}_\star (\calM) = \Fil^r (\tilde A \otimes_{\tilde S} 
T(\calM))_{N=0}^{\phi_r = 1}.$$
\end{prop}

\begin{preuve}
Si $\calM$ est tué par $p$, le $\Ast$-module $\Ast \otimes_{S} \calM$ 
l'est également. Il s'identifie au $\Ast/p\Ast$-module $\Ast/p\Ast
\otimes_{S_1} \calM$, et on dispose d'une application canonique :
$$\Ast/p\Ast \otimes_{S_1} \calM \to \tilde A \otimes_{\tilde S} 
T(\calM).$$
Cette application est clairement surjective, et en reprenant les 
notations du lemme \ref{4:lem:kappa}, son noyau s'identifie à $\kappa 
(\Ast/p\Ast \otimes_{S_1} \calM)$. En outre, on vérifie directement 
qu'elle induit une flèche compatible à l'action de $G_K$ :
$$\Psi : \Fil^r(\Ast/p\Ast \otimes_{S_1} \calM)_{N=0}^{\phi_r=1} \to 
\Fil^r(\tilde A \otimes_{\tilde S} T(\calM))_{N=0}^{\phi_r=1}.$$
Reste à prouver que cette dernière application est un isomorphisme.
Soit $x \in \ker \Psi$. On a $x \in \Fil^r(\Ast/p\Ast \otimes_{S_1} 
\calM)_{N=0}^{\phi_r=1}$ et donc $\phi_r(x) = x$. Par ailleurs, on a $x 
\in \kappa (\Ast/p\Ast \otimes_{S_1} \calM)$, d'où on déduit $x = 
\phi_r(x) \in \phi_r(\kappa) (\Ast/p\Ast \otimes_{S_1} \calM)$ puis $x = 
\phi_r(x) = 0$, puisque $\phi_r \circ \phi_r(\kappa) = 0$ (lemme 
\ref{4:lem:kappa}). L'application $\Psi$ est donc injective.

Passons à la surjectivité. Considérons $x \in \Fil^r(\tilde A 
\otimes_{\tilde S} T(\calM))_{N=0}^{\phi_r=1}$ et notons $\hat x \in 
\Ast/p\Ast \otimes_{S_1} \calM$ un relevé quelconque de $x$. On vérifie 
alors que $\hat x \in \Fil^r(\Ast/p\Ast \otimes_{S_1} \calM)$ et que
$\phi_r(\hat x) = \hat x + y$ pour un certain $y \in \kappa (\Ast/p\Ast
\otimes_{S_1} \calM)$. Posons $\hat x_1 = \hat x + y \in \Fil^r 
(\Ast/p\Ast \otimes_{S_1} \calM)$. On a $\hat x_1 \in \Fil^r 
(\Ast/p\Ast \otimes_{S_1} \calM)$ et $\phi_r(\hat x_1) = \hat x_1 + y'$ 
pour un certain $y' \in \phi_r(\kappa) (\Ast/p\Ast \otimes_{S_1} 
\calM)$.
On pose alors $\hat x_2 = \hat x_1 + y' \in \Fil^r (\Ast/p\Ast 
\otimes_{S_1} \calM)^{\phi_r = 1}$. On vérifie finalement que $N$ annule 
$\hat x_2$, ce qui assure que $\hat x_2$ est un antécédent par $\Psi$ de 
$x$.
\end{preuve}

\subsubsection{Lien entre les foncteurs $\Tst^\star$ et ${\Tst}_\star$}

Fixons $\epsilon$ une suite de racines 
$p^n$-ièmes de l'unité et définissons $t = \log \pa{\cro{\epsilon}} \in 
\Acris$ où le $\log$ est donné par la série usuelle. On a $\Phi 
\pa t = pt$ et plus exactement l'ensemble des solutions dans $\Acris$ de 
$\Phi\pa t = pt$ est le $\Z_p$-module engendré par $t$. C'est un 
$\Z_p$-module libre de rang $1$ isomorphe à $\Z_p\pa 1$ en tant que 
représentation galoisienne. Autrement dit $(\Fil^r 
\Ast)^{\phi_r=1}_{N=0} = (\Fil^r \Acris)^{\phi_r=1}$ est isomorphe en 
tant que représentation galoisienne à $\Z_p(r)$. De même, la 
représentation galoisienne $(\Fil^r \Astinf)^{\phi_r=1}_{N=0}$ est 
isomorphe à $\Q_p/\Z_p(r)$.

On ne sait toujours pas que le foncteur ${\Tst}_\star$ est bien défini 
mais si $\calM$ est un objet de $\Mr$ muni d'une filtration admissible 
fixée, on peut définir une application :
$$\Psi : {\Tst}_\star(\calM) \to \Tst^\star(\calM)^\vee (r)$$
où par définition $\Tst^\star(\calM)^\vee (r)$ est la représentation 
galoisienne $\hom(\Tst^\star(\calM), \Q_p/\Z_p(r))$.
En effet, soient $x = \sum_i a_i \otimes x_i \in \Fil^r(\Ast 
\otimes_S \calM)^{\phi_r=1}_{N=0}$, et $f : \calM \to \Astinf$ 
compatible aux structures. On pose $\Psi(x)(f) = \sum_i a_i f(x_i)$. 
C'est un élément de $(\Fil^r \Astinf)^{\phi_r=1}_{N=0}$ soit, en vertu 
de l'isomorphisme décrit précédemment, un élément de $\Q_p/\Z_p(r)$.
On vérifie immédiatement que $\Psi$ est un morphisme $\Z_p$-linéaire
commutant à l'action de Galois.

\bigskip

On veut démontrer que $\Psi$ est un isomorphisme pour tout objet 
$\calM \in \Mr$ et pour cela on suit la démonstration du paragraphe 
3.2.1. de \cite{breuil-duke} (qui concerne le cas $e=1$) : la 
méthode consiste à prouver que $\Psi$ est un isomorphisme pour les 
objets tués par $p$, à démontrer que le foncteur ${\Tst}_\star$ est bien 
défini et exact, puis à conclure à l'aide d'un dévissage.

\medskip

On rappelle, dans un premier temps, que les objets simples de $\Mr$ ont 
une structure relativement simple (proposition 4.1.1 de \cite{caruso}) :

\begin{prop}
\label{4:prop:simple}
Soit $\calM$ un objet simple de $\Mr$. Alors $\calM$ est tué par $p$ et 
admet une base adaptée $(e_1, \ldots, e_d)$ pour des entiers $n_1, 
\ldots, n_d$. De plus $N(e_i) = 0$ pour tout $i$ et il existe une 
(unique) matrice $G$ inversible à coefficients dans $k$ telle que :
$$\pa { \begin{array}{c} \phi_r(u^{n_1} e_1) \\ \vdots \\
\phi_r(u^{n_d} e_d) \end{array} } = {}^t G \pa { \begin{array}{c} e_1 \\
\vdots \\ e_d \end{array} }.$$
\end{prop}

\bigskip

On introduit les modules suivants :
\begin{equation}
\label{4:eq:filx}
\Fil^t_X \Ast = \acco{\sum_{i=t}^m a_i \frac{X^i}{i!}, \, m 
\in \N, \, a_t \in \Acris}
\end{equation}
et on démontre (comme dans le lemme 3.1.2.1. de \cite{breuil-duke}) que 
$\Fil^t_X (\Ast / p^n) = (\Fil^t_X \Ast) / p^n$ est plat sur $S_n$. Si 
$\calM$ est un objet de $\Mr$ et $(Fil^t\calM)$ est une filtration 
admissible de $\calM$, on définit pour tout $s$ compris entre $0$ et $r$ :
$$\Fil^s_X(\Ast \otimes_{S} \calM) = \sum_{t=0}^s \Fil^t_X 
\Ast \otimes_S \Fil^{s-t} \calM.$$
On a les deux lemmes suivants :

\begin{lemme}
\label{4:lem:filconoyau}
Soit $\calM$ un objet de $\Mr$ et $(\Fil^t \calM)$ une filtration 
admissible de $\calM$. Alors, pour tout $s$ compris entre $0$ et $r$, on 
a $\Fil^s_X(\Ast \otimes_{S} \calM) = \Fil^s(\Ast \otimes_{S} \calM)$ et 
si $s \geq 1$, on a une suite exacte :
$${\small \xymatrix @C=20pt {
0 \ar[r] & \displaystyle \xybigoplus {t=1} s \Fil^t_X \xyAst 
\otimes_S \Fil^{s+1-t} \calM \ar[r] & \displaystyle 
\xybigoplus {i=0} s \Fil^t_X \xyAst \otimes_S \Fil^{s-t} \calM 
\ar[r] & \Fil^s(\xyAst \otimes_{S} \calM) \ar[r] & 0 }}$$
où la première flèche est la somme des applications suivantes :
$${\small \begin{array}{ccccc}
\Fil^t_X \Ast \otimes_S \Fil^{s+1-t} \calM& \to & 
\Fil^{t-1}_X \Ast \otimes_S \Fil^{s+1-t} \calM & \oplus &
\Fil^t_X \Ast \otimes_S \Fil^{s-t} \calM \\
x_t \otimes y_{s+1-t} & \mapsto & x_t \otimes y_{s+1-t} & \oplus &
- x_t \otimes y_{s+1-t}
\end{array}}$$
\end{lemme}

\begin{preuve}
La preuve est une adaptation de celle du lemme 3.2.1.2 de 
\cite{breuil-duke}. Pour la première assertion, on a déjà 
clairement $\Fil^s_X(\Ast \otimes_{S} \calM) \subset \Fil^s(\Ast 
\otimes_{S} \calM)$.

On a une description alternative de $\Acris$ (voir 
\cite{fontaine-ast1}). Si on note $R = \varprojlim_{n\in \N} \O_{\bar 
K}/p$ pour les applications de transitions données par le Frobenius, on 
peut définir un morphisme :
\[\begin{array}{rcl}
\hat\theta \quad : \quad W\pa R & \to & \O_{\C_p} \\
\pa{a_0, a_1, \ldots, a_n, \ldots} & \mapsto & \sum_{n \geq 0} p^n 
\hat x_n^{\pa n}
\end{array} \]
où $\C_p$ désigne le complété $p$-adique de $\bar K$ et où $\hat 
x_n^{\pa n}$ est la limite quand $m$ tend vers l'infini d'une suite 
$(\hat a_n^{\pa{n+m}})^{p^m}$, $\hat a_i^{\pa j} \in \O_{\bar K}$ 
désignant un relevé quelconque de $a_i^{\pa j}$. L'anneau $\Acris$ 
s'identifie alors à l'enveloppe aux puissances divisées de $W(R)$ 
relativement à $\ker \hat\theta$ (et compatibles aux puissances divisées 
canoniques sur $p W(R)$). On vérifie facilement que $[\upi] \in 
\Acris$ (défini à la fin du paragraphe \ref{4:subsec:acris}) correspond 
bien au représentant de Teichmüller de $\upi$ (défini au même endroit). 
Par ailleurs, on montre (voir \cite{fontaine-ast1}) que $\ker 
\hat\theta$ est un idéal principal, engendré par $E([\upi])$.

Ainsi, si $x \in \Fil^s(\Ast \otimes_{S} \calM)$, il s'écrit comme une 
somme de termes de la forme :
$$a \gamma_j(E([\upi])) \gamma_k(X) \otimes m$$
(où $\gamma_j(x) = \frac{x^j}{j!}$) avec $a \in \Acris$, $m \in 
\Fil^{s-t} \calM$ et $j+k \geq t$. Mais $[\upi] = u(1+X)$ et donc 
$E([\upi]) - E(u)$ est un multiple de $uX$. On 
peut donc écrire $E([\upi]) = E(u) + u X b$ pour un certain $b \in 
\Ast$ puis :
$$\gamma_j(E([\upi])) = \sum_{\ell=0}^j \gamma_{j-\ell}(E(u)) (u b)^\ell 
\gamma_\ell(X).$$
En réinjectant cette expression dans $x$, on obtient bien $x \in 
\Fil^s_X(\Ast \otimes_{S} \calM)$.

Pour la deuxième partie du lemme, la démonstration est exactement la 
même que celle de \cite{breuil-duke}.
\end{preuve}

\begin{lemme}
Soit $\calM$ un objet de $\Mr$. Pour tout $s$ compris entre $0$ et $r$, 
on a une suite exacte :
$$\xymatrix @C=20pt {
0 \ar[r] & \Fil^s(\xyAst \otimes_S \calM)_{N=0} \ar[r] &
\Fil^s(\xyAst \otimes_S \calM) \ar[r]^N & \Fil^{s-1}(\xyAst \otimes_S \calM)
 \ar[r] & 0 }$$
où par convention $\Fil^{-1}(\Ast \otimes_S \calM) = \Ast \otimes_S 
\calM$.
\end{lemme}

\begin{preuve}
La démonstration est la même que celle du lemme 3.2.1.3 de 
\cite{breuil-duke}. Toutefois, on se ramène à la fin, non pas à un
objet de $\underline{MF}_k^{f,r}$ (en reprenant les notations de 
l'article), mais à un objet simple de la catégorie $\Mr$ dont la
structure est connue par la proposition \ref{4:prop:simple}. Le même 
argument s'applique alors.
\end{preuve}

\medskip

On déduit des deux lemmes précédents le corollaire suivant :

\begin{cor}
\label{4:cor:filadm}
Soit $\calM$ un objet de $\Mr$. Alors $\Fil^r(\Ast \otimes_S 
\calM)_{N=0}$ ne dépend pas de la filtration admissible choisie et si :
$$\xymatrix @C=20pt {
0 \ar[r] & \calM' \ar[r] & \calM \ar[r] & \calM'' \ar[r] & 0 }$$
est une suite exacte dans $\Mr$, alors la suite induite :
$$\xymatrix @C=20pt {
0 \ar[r] & \Fil^r(\xyAst \otimes_S \calM')_{N=0} \ar[r] & \Fil^r(\xyAst 
\otimes_S \calM)_{N=0} \ar[r] & \Fil^r(\xyAst \otimes_S \calM'')_{N=0} 
\ar[r] & 0 }$$
est également exacte.
\end{cor}

\begin{preuve}
C'est la même que celle de la proposition 3.2.1.4 de \cite{breuil-duke}, 
en remplaçant à nouveau $\underline{MF}_k^{f,r}$ par la sous-catégorie
de $\Mr$ formée des objets tués par $p$, et en utilisant la proposition
\ref{4:prop:simple} qui donne la structure de tels objets.
\end{preuve}

\bigskip

\noindent
{\it Remarque.} On prouve de même que si $\xymatrix @C=20pt {
0 \ar[r] & \calM' \ar[r] & \calM \ar[r] & \calM'' \ar[r] & 0 }$ est une 
suite exacte dans $\Mr$, alors la suite :
$$\xymatrix @C=20pt {
0 \ar[r] & (\xyAst \otimes_S \calM')_{N=0} \ar[r] & (\xyAst 
\otimes_S \calM)_{N=0} \ar[r] & (\xyAst \otimes_S \calM'')_{N=0} 
\ar[r] & 0 }$$
l'est aussi.

D'autre part, le corollaire précédent prouve en particulier que 
${\Tst}_\star (\calM)$ ne dépend pas de la filtration choisie. Ainsi
le foncteur ${\Tst}_\star$ est bien défini.

\bigskip

\begin{lemme}
\label{4:lem:phirsurj}
Si $\calM$ est un objet de $\Mr$, on a une suite exacte :
$$\xymatrix @C=25pt {
0 \ar[r] & \Fil^r(\xyAst \otimes_S \calM)_{N=0}^{\phi_r=1} \ar[r] &
\Fil^r(\xyAst \otimes_S \calM)_{N=0} \ar[r]^-{\phi_r-\id} & (\xyAst 
\otimes_S \calM)_{N=0} \ar[r] & 0}.$$
\end{lemme}

\begin{preuve}
Il suffit de montrer que $\phi_r - \id$ est surjective. De même que dans 
le lemme 3.2.1.6 de \cite{breuil-duke}, on se ramène au cas d'un objet 
simple de $\Mr$.

Soit $\calM$ un objet simple de $\Mr$. D'après la proposition 
\ref{4:prop:simple}, on peut écrire $\calM = S_1 e_1 \oplus \cdots \oplus
S_1 e_d$ où $e_1, \ldots, e_d$ sont tels que $N(e_i) = 0$ et $\Fil^r 
\calM$ est le sous-module de $\calM$ engendré par $\Fil^p S \cdot 
\calM$ et les $u^{n_1} e_1, \ldots u^{n_1} e_d$ pour certains entiers 
$n_i$ compris entre $0$ et $er$. On a alors directement :
\begin{eqnarray*}
\Ast \otimes_S \calM & = & \Ast / p\Ast \cdot e_1 \oplus \cdots \Ast / 
p\Ast \cdot e_d \\
(\Ast \otimes_S \calM)_{N=0} & = & \Acris / p\Acris \cdot e_1 \oplus 
\cdots \Acris / p\Acris \cdot e_d \\
\Fil^r (\Ast \otimes_S \calM)_{N=0} & = & \sum_{i=1}^p (\Acris / p\Acris 
\cdot u^{n_i} e_1 + \Fil^p \Acris / p\Acris \cdot e_i).
\end{eqnarray*}
On note $G$ l'unique matrice inversible à coefficients dans $k$ telle 
que :
$$\pa { \begin{array}{c} \phi_r(u^{n_1} e_1) \\ \vdots \\ 
\phi_r(u^{n_d} e_d) \end{array} } = {}^t G \pa { \begin{array}{c} e_1 \\ 
\vdots \\ e_d \end{array} }$$
et on conclut de même que dans le lemme 3.2.1.6 de \cite{breuil-duke}.
\end{preuve}

\begin{cor}
Le foncteur ${\Tst}_\star$ est exact.
\end{cor}

On déduit finalement de cette étude le théorème suivant :

\begin{theo}
L'application $\Psi$ définie précédemment induit une transformation 
naturelle inversible entre les foncteurs ${\Tst}_\star$ et 
$(\Tst^\star)^\vee (r)$.
\end{theo}

\begin{preuve}
Comme la catégorie $\Mr$ est artinienne, et que les foncteurs 
$\Tst^\star$ et ${\Tst}_\star$ sont exacts, il suffit de montrer le 
résultat lorsque $\calM$ est un objet simple de $\Mr$.

Si $\calM$ est un objet simple de $\Mr$, la proposition \ref{4:prop:simple}
nous assure dans un premier temps que $\calM$ est tué par $p$. On a donc :
$$\Tst^\star(\calM) = \hom(T(\calM), \tilde A)
\qquad \text{et} \qquad
{\Tst}_\star(\calM) = \Fil^r (\tilde A \otimes_{\tilde S} 
T(\calM))^{\phi_r = 1}_{N=0}.$$
Par ailleurs la même proposition fournit une description explicite de 
$T(\cal M)$ et de ses structures supplémentaires. Précisément, il existe
des entiers $n_i$ tels que :
\begin{eqnarray*}
T(\calM) & = & \tilde S e_1 \oplus \cdots \oplus \tilde S e_d \\
\Fil^r T(\calM) & = & \tilde S u^{n_1} e_1 \oplus \cdots \oplus \tilde S 
u^{n_d} e_d
\end{eqnarray*}
avec de surcroît $N(e_i) = 0$ pour tout $i$. Par ailleurs, quitte à 
passer à une extension non ramifiée de $K$, on peut 
supposer (voir théorème 4.3.2 de \cite{caruso}) que 
$\phi_r$ est donné par $\phi_r(u^{n_i} e_i) = e_{i+1}$, les indices 
étant considérés modulo $d$.

Des descriptions précédentes, on déduit facilement :
$$\Fil^r (\tilde A \otimes_{\tilde S} T(\calM))_{N=0}
= \bar \pi_1^{n_1} \O_{\bar K}/\pi \cdot e_1 \oplus \cdots \oplus
\bar \pi_1^{n_1} \O_{\bar K}/\pi \cdot e_d$$
où $\bar \pi_1$ est la réduction modulo $\pi$ de $\pi_1$ (qui, on le
rappelle, est une racine $p$-ième de $\pi$ fixée). L'opérateur $\phi_r$ 
agit sur ce module par $\phi_r(\bar \pi_1^{n_i} e_i) = e_{i+1}$.
Tout élément de $\Fil^r (\tilde A \otimes_{\tilde S} T(\calM))_{N=0}$ 
s'écrit de façon unique $x = \sum_{i=1}^d a_i \otimes e_i$ avec $a_i = 
\pi_1^{n_i} x_i$ et un tel élément appartient à ${\Tst}_\star(\calM)$ 
si, et seulement si :
$$x_i^p = \bar\pi_1^{n_{i+1}} x_{i+1}$$
pour tout indice $i \in \Z/d\Z$.

Par ailleurs, se donner un élément de $\Tst^\star(\calM)$ revient à se 
donner l'image $b_i$ de chacun des $e_i$, ces images devant vérifier 
$N(b_i) = 0$, $u^{n_i} b_i \in \Fil^r \tilde A$ et commuter à $\phi_r$.
La première condition impose $b_i \in \O_{\bar K} / \pi$. La deuxième 
condition assure que $b_i = \pi_1^{m_i} y_i$ pour $m_i = er - n_i$ et 
$y_i \in \O_{\bar K}/\pi$. Finalement, la commutation à $\phi_r$ impose 
les relations :
$$(-1)^r y_i^p = \bar \pi_1^{m_{i+1}} y_{i+1}.$$

On est finalement ramené à prouver que l'accouplement :
$$(a_1, \ldots, a_d) \times (b_1, \ldots, b_d) \mapsto \sum_{i=1}^d 
a_i b_i$$
défini sur les couples de $d$-uplets solutions des systèmes précédents 
et à valeurs dans $(\O_{\bar K} / \pi)^{\phi_r = 1}$. Ce dernier espace 
est encore $\bar t \F_p$ (ou $\bar t$ est la réduction modulo $\pi$ 
d'une racine $(p-1)$-ième de $p^r$) est non dégénéré. Or par le lemme 
5.1.2 de \cite{caruso}\footnote{Dans cette référence, on travaillait 
non pas modulo $\pi$ mais modulo $p$. Cependant, on vérifie sans mal que 
la méthode marche dans les deux cas.}, si on choisit $\eta$ une racine 
$(p^h-1)$-ième de $\pi$, si on note $\bar \eta$ sa réduction modulo 
$\pi$, et si on pose :
\begin{eqnarray*}
s_i & = & n_i p^{d-1} + n_{i+1} p^{d-2} + \cdots + n_{i+d-1} \\
t_i & = & m_i p^{d-1} + m_{i+1} p^{d-2} + \cdots + m_{i+d-1}
\end{eqnarray*}
les solutions de ces systèmes s'écrivent :
$$a_i = a^{p^i} \bar \eta^{s_i} 
\quad \text{et} \quad
b_i = (-1)^{ri} b^{p^i} \bar \eta^{t_i}$$
où $a$ décrit $\F_q$ ($q = p^d$), l'ensemble des racines dans 
$\bar k$ de l'équation $x^q = x$, et où $b$ décrit l'ensemble des 
racines dans $\bar k$ de l'équation $x^q = (-1)^{rd} x$.

Si $rd$ est pair, $\sum_{i=1}^d a_i b_i = \tr_{\F_q/\F_p} (ab) \cdot 
\bar \eta^v$ où $v = s_i + t_i = er \cdot \frac{q-1}{p-1}$ est 
indépendant de $i$, et on conclut en remarquant que la trace de $\F_q$ 
à $\F_p$ est une forme bilinéaire non dégénérée.

Si $rd$ est impair, on note $\varepsilon \in \bar k$ une racine 
$(p-1)$-ième de $-1$, on vérifie que $\sum_{i=1}^d a_i b_i = 
\varepsilon \tr_{\F_q/\F_p} (ab/\varepsilon) \cdot \bar \eta^v$ et 
on conclut comme dans le cas précédent.
\end{preuve}

\bigskip

En vertu de ce théorème, tous les résultats démontrés sur le foncteur 
$\Tst^\star$ se transposent au foncteur ${\Tst}_\star$. On obtient ainsi 
un équivalent du théorème \ref{4:th:tstcont} :

\begin{theo}
\label{4:th:tstcov}
Le foncteur ${\Tst}_\star$ est exact, pleinement fidèle, d'image 
essentielle stable par quotients et par sous-objets. De plus, si $\calM$ 
est un objet de $\Mr$ isomorphe en tant que $S$-module à $S_{n_1} \oplus 
\cdots \oplus S_{n_d}$, alors la représentation galoisienne 
${\Tst}_\star(\calM)$ est isomorphe en tant que $\Z_p$-module à 
$\Z/p^{n_1}\Z \times \cdots \times \Z/p^{n_d}\Z$.
\end{theo}

\section{Les faisceaux sur le site log-syntomique}
\label{4:sec:synto}

\subsection{Rappels et préliminaires}
\subsubsection{Log-schémas et sites usuels}

On renvoie à \cite{kato} pour la définition et les propriétés des 
log-schémas et des morphismes de log-schémas (en particulier des 
morphismes log-lisses ou log-étales). Tous les log-schémas considérés 
dans ce papier sont intègres. Si $S$ 
est un log-schéma, on note $\dot S$ le schéma sous-jacent. Si $M$ est un 
monoïde, on note $M^\gp$  le groupe associé, c'est-à-dire l'ensemble des 
éléments de la forme $a b^{-1}$ pour $a$ et $b$ dans $M$ où deux 
éléments $a b^{-1}$ et $c d^{-1}$ sont identifiés s'il existe $s \in M$ 
tel que $sad = sbc$ (\emph{i.e.} simplement $ad=bc$ si $M$ est intègre).
Si $M$ est un monoïde et $G$ un sous-groupe de $M^\gp$, on définit le 
quotient $M/G$ comme le quotient de $M$ par la relation d'équivalence
$x \sim y \leftrightarrow x y^{-1} \in G$.

Si $S$ est un log-schéma fin dont le schéma sous-jacent est tué par un 
entier non nul et muni d'un idéal à puissances divisées et si $X$ est un 
log-schéma fin sur $T$ auquel les puissances divisées s'étendent, on 
note $(X/S)_\cris$ le petit site (log-)cristallin fin (défini dans 
\cite{kato}, chapitre 5) et $(X/S)_\CRIS$ le gros site (log-)cristallin 
fin (défini dans \cite{breuil-bull}, chapitre 3). On note $\O_{X/S}$ le 
faisceau structural sur ces sites, $\J_{X/S}$ son idéal à puissances 
divisées et $\J_{X/S}^{[n]}$ les puissances divisées successives de 
$\J_{X/S}$.

De même, si $S$ est un log-schéma fin, on note $S_\et = \dot S_\et$ le 
petit site (log-)étale de $S$ : c'est la catégorie des log-schémas $X$ 
pour lesquels $\dot X$ est étale sur $\dot S$ et la log-structure sur 
$X$ est induite par celle de $S$, les recouvrements étant les 
recouvrements étales usuels (sur les schémas sous-jacents). On note 
également $S_\ET$ le gros site (log-)étale de $S$ défini comme la 
catégorie des log-schémas fins localement de type fini sur $S$ et munie 
de la topologie étale. On note $\O_X$ le faisceau structural sur chacun 
de ces deux sites.

\subsubsection{Topologie log-syntomique}

On rappelle la définition d'un morphisme de log-schémas log-syntomiques, 
due à Kato (voir \cite{kato}) :

\begin{deftn}
Un morphisme de log-schémas fins $f : Y \to X$ est dit 
\emph{log-syntomique} s'il est intègre, si $\dot f : \dot Y \to \dot 
X$ est localement de présentation finie, et si $f$ peut s'écrire étale 
localement comme la composée d'un morphisme log-lisse avec une immersion 
fermée exacte dont l'idéal est engendré en chaque point par une suite 
transversalement régulière relativement à $X$.
\end{deftn}

On montre (voir \cite{breuil-bull}) que les morphismes log-syntomiques 
sont stables par changement de base et par composition. En outre, on
dispose de la propriété remarquable suivante :

\begin{prop}
\label{4:prop:relev}
Si $Y \to X$ est une immersion fermée exacte définie par un nil-idéal, 
on peut étale-localement relever les morphismes log-syntomiques sur $Y$ 
en des morphismes log-syntomiques sur $X$.
\end{prop}

Il est de plus possible de donner une description locale très explicite 
des morphismes log-syntomiques. Précisément, si $f : Y \to X$ est un 
morphisme log-syntomique, alors il existe une carte (locale pour la 
topologie étale) sur laquelle $f$ prend la forme suivante :
$$\xymatrix @C=50pt {
\displaystyle \frac{P \oplus \N^r} G \ar[r] & \displaystyle \frac{(A 
\otimes_{\Z[P]} \Z[(P \oplus \N^r) G] [X_1, \ldots, X_s]} I \\
P \ar[r] \ar[u]_-{f} & A \ar[u]^-{f} }$$
où $A$ est un anneau, $P$ est un monoïde intègre, $r$ et $s$ des 
entiers, $G$ un sous-groupe de $P^\gp \oplus \Z^r$ et $I$ un idéal 
contenant les $[g]-1$ (pour $g \in G$) et engendré par une suite 
transversalement régulière relativement à $A$.

\bigskip

Si $S$ est un log-schéma fin, on note $S_\SYN$ le gros site 
(log-)syntomique sur $X$, c'est-à-dire la catégorie des log-schémas fins 
localement de type fini sur $X$ munie de la topologie log-syntomique : 
une famille $(f_i : X_i \to X)$ recouvre $X$ si tous les morphismes 
$f_i$ sont log-syntomiques et si topologiquement $\dot X = \bigcup_i 
f_i(\dot X_i)$. De même on définit le petit site (log-)syntomique 
$S_\syn$ en se restreignant à la catégorie des log-schémas 
log-syntomiques sur $S$.

\subsubsection{Plusieurs bases}
\label{4:subsec:bases}

Pour la suite, on sera amené à considérer plusieurs bases qui sont :
$$
\pa{ \begin{array}{rcl}
\N & \to & \O_K/p^n \\
1 & \mapsto & \pi \end{array} }
\quad ; \quad
\pa{ \begin{array}{rcl}
\N & \to & S_n \\
1 & \mapsto & u \end{array} }
\quad ; \quad
\pa{ \begin{array}{rcl}
\N & \to & \tilde S \\
1 & \mapsto & u \end{array} }
\quad ; \quad
\pa{ \begin{array}{rcl}
\N & \to & k \\
1 & \mapsto & 0 \end{array} }
$$
On note simplement $T_n$ la première, $E_n$ la deuxième, $\tilde E$ la 
troisième et $\bar E = \bar T$ la quatrième. Ces quatre bases sont 
munies de puissances divisées : sur $T_n$, elles sont prises par rapport 
à l'idéal engendré par $p$, sur $E_n$ et $\tilde E$ par rapport à 
l'idéal engendré par les $\frac{E(u)^i}{i!}$ pour $i \geq 1$ et sur la 
dernière par rapport à l'idéal nul. On a un diagramme :
$$\xymatrix @C=7pt {
& T_1 \arinj[rrrr] \arinj[ld] \arinj[rd] & & & & T_2 \arinj[rrr] 
\arinj[d] & & & \cdots \arinj[rrr] & & & T_n \arinj[d] \arinj[rrr] & & & 
\cdots \\
\xytilde E \arinj[rr] & & E_1 \arinj[rrr] & & & E_2 \arinj[rrr] & & & 
\cdots \arinj[rrr] & & & E_n \arinj[rrr] & & & \cdots }$$
où tous les morphismes sont des épaississements, les flèches
verticales étant obtenues en envoyant $u$ sur $\pi$.

\bigskip

Les bases $E_n$, $\tilde E$ et $\bar E$ sont munies d'un relèvement du 
Frobenius : c'est la multiplication par $p$ sur les monoïdes, 
l'élévation à la puissance $p$ sur $\tilde S$ et $k$ et le 
Frobenius défini au paragraphe \ref{4:subsec:s} sur $S_n$.

\subsection{Les faisceaux $\O^\st_n$ et $\J_n^{[s]}$}
Dans ce paragraphe, on définit des faisceaux $\O^\st_n$ et $\J_n^{[s]}$ 
sur le site syntomique qui permettent de calculer la cohomologie 
cristalline. On donne ensuite une description locale explicite de ces
faisceaux, technique mais cruciale pour mener à bien les calculs.

\subsubsection{Définition et description locale}
\label{4:subsec:ostlocal}

Pour tout entier $n$ et tout entier (relatif) $s$, on définit sur 
$(T_n)_\SYN$ les préfaisceaux $\J_n^{[s]}$ par la formule :
$$\J_n^{[s]}(U) = H^0 ((U/E_n)_\cris, \J^{[s]}_{U/E_n}) = H^0 
((U/E_n)_\CRIS, \J^{[s]}_{U/E_n})$$
où $U$ sur $T_n$ est vu sur $E_n$ \emph{via} l'épaississement du
paragraphe \ref{4:subsec:bases}. On pose $\O^\st_n = \J_n^{[0]}$. On montre 
(voir \cite{breuil-bull}) que les $\J_n^{[s]}$ sont des faisceaux et 
qu'ils calculent la cohomologie log-cristalline dans le sens où :
$$H^i(X_\syn, \J_n^{[s]}) = H^i(X_\SYN, \J_n^{[s]}) = H^i 
((X/E_n)_\cris, \J^{[s]}_{X/E_n}) = H^i((X/E_n)_\CRIS, 
\J^{[s]}_{X/E_n})$$
pour tout entier $i$ et tout log-schéma $X$ fin localement de type fini 
sur $E_n$.

\bigskip

Soit $U$ un log-schéma log-syntomique sur la base $T_n$. On a vu 
qu'étale-localement, on peut trouver une carte du morphisme $U \to T_n$ 
qui prend la forme :
$$\xymatrix @C=50pt {
\displaystyle \frac{\N x_0 \oplus \N x_1 \oplus \cdots \oplus \N x_r} G 
\ar[r]^-{\alpha} & \displaystyle \frac{\O_K/p^n [(\N x_0 \oplus \N x_1 
\oplus \cdots \oplus \N x_r) G] [X_1, \ldots, X_s]} {([x_0] - \pi, f_1, 
\ldots, f_t)} \\
\N x_0 \ar[r] \ar[u] & \O_K/p^n \ar[u] }$$
où $G$ est un sous-groupe de $\Z^{r+1}$ et $([x_0] - \pi, f_1, \ldots, 
f_t)$ est une suite transversalement régulière relativement à $\O_K/p^n$ 
et telle que l'idéal engendré contienne tous les $[g]-1$, pour $g \in 
G$.

\medskip

Notons $Q = \N x_0 \oplus \N x_1 \oplus \cdots \oplus \N x_r$, $P = 
Q/G$ et $A = \frac{\O_K/p^n [Q G] [X_1, \ldots, X_s]} {([x_0] - \pi, 
f_1, \ldots, f_t)}$, de sorte que (étale-)localement $U = (\spec A, 
P)$. Décrire localement (pour la topologie syntomique) les faisceaux 
$\O^\st_n$ et $\J_n^{[s]}$ serait par exemple donner des formules 
explicites pour les modules $\O^\st_n (\spec A, P)$ et $\J_n^{[s]} (\spec 
A, P)$. Cependant, on ne sait donner de telles formules que si le 
Frobenius est surjectif sur $A$ et sur $P$, ce qui n'est \emph{a priori} 
pas le cas ici. Nous allons donc devoir considérer des ouverts encore 
plus petits (toujours pour la topologie syntomique) pour forcer cette 
condition de surjectivité.

\medskip

Notons $Q^i = \N x_0^{p^{-i}} \oplus \N x_1^{p^{-i}} \oplus \cdots 
\oplus \N x_r^{p^{-i}}$, $P^i = Q^i / G$ et $A^i = \frac{\O_K/p^n [Q^i 
G] [X_1^{p^{-i}}, \ldots, X_s^{p^{-i}}]} {([x_0] - \pi, f_1, \ldots, 
f_t)}$. On a des morphimes de log-schémas $(\spec A^{i+1}, P^{i+1}) \to 
(\spec A^i, P^i)$ qui sont des recouvrements log-syntomiques, et le 
Frobenius devient surjectif sur la limite de ces recouvrements. On est 
amené à décrire explicitement les objets suivants :
$$\O^\st_n(A^\infty, P^\infty) = \varinjlim_i \O^\st_n(\spec A^i, P^i)
\quad \text{et} \quad
\J_n^{[s]} (A^\infty, P^\infty) = \varinjlim_i \J_n^{[s]}(\spec A^i, 
P^i).$$
Notons $Q^\infty = \N x_0^{1/p^\infty} \oplus \N x_1^{1/p^\infty} 
\oplus \cdots \oplus \N x_r^{1/p^\infty}$, $P^\infty = Q^\infty/G$ 
et $A^\infty = \frac{\O_K/p^n [Q^\infty G] [X_1^{1/p^\infty}, \ldots, 
X_s^{1/p^\infty}]} {([x_0] - \pi, f_1, \ldots, f_t)}$. Par ailleurs,
si $M$ est un monoïde et $n$ un entier, posons :
$$M^{(n)} = \acco{x \in M^\gp\, / \, x^{p^n} \in M}.$$
Si de plus $M$ est muni d'un morphisme de monoïdes $\N \to M$, on
définit $\N \oplus_{(\phi^n), \N} M$ comme la limite inductive du 
diagramme $\xymatrix @C=15pt {\N & \N \ar[l]_{p^n} \ar[r] & M}$. 
Posons :
$$W_n^\st (A^\infty, P^\infty) = W_n (A^\infty/p A^\infty) \otimes_{\Z 
[P^\infty]} \Z[(\N \oplus_{(\phi^n), \N} P^\infty)^{(n)}].$$
On dispose d'un morphisme canonique surjectif :
$$\begin{array}{rcl}
W_n^\st (A^\infty, P^\infty) & \to & A^\infty \\
(a_0, \ldots, a_{n-1}) \otimes [h] & \mapsto & (\hat a_0^{p^n} + \cdots 
+ p^{n-1} \hat a^{n-1}) \cdot \alpha(h^{p^n})
\end{array}$$
où $\hat a_i \in A^\infty$ désigne un relevé quelconque de $a_i$.
On note finalement $W_n^{\st,\DP} (A^\infty, P^\infty)$ l'enveloppe aux 
puissances divisées relativement à l'ideal noyau (et compatible aux 
puissances divisées canoniques sur l'idéal $(p)$). On munit $W_n^{\st, 
\DP} (A^\infty, P^\infty)$ d'une structure de $S_n$-module en envoyant 
$u$ sur l'élément $1 \otimes [1 \oplus (0, \ldots, 0)]$. Comme dans 
l'appendice D de \cite{breuil-duke}, on montre la proposition suivante :

\begin{prop}
Avec les notations précédentes, il existe un isomorphisme $S_n$-linéaire 
canonique :
$$\xymatrix{\O^\st_n (A^\infty, P^\infty) \ar[r]^-{\sim} & W_n^{\st, \DP} 
(A^\infty, P^\infty)}.$$
\end{prop}

Il existe une autre description locale qui a l'avantage d'être 
légèrement plus simple, mais l'inconvénient d'être non canonique. Notons 
pour cela :
$$W_n^\cris (A^\infty, P^\infty) = W_n (A^\infty/p A^\infty) \otimes_{\Z 
[P^\infty]} \Z[(P^\infty)^{(n)}]$$
et $W_n^{\cris, \DP}$ son enveloppe à puissances divisées par rapport à 
l'idéal noyau de l'application $(a_0, \ldots, a_{n-1}) \otimes [h]
\mapsto (\hat a_0^{p^n} + \cdots + p^{n-1} \hat a_{n-1}) 
\alpha(h^{p^n})$ où $\hat a_i \in A^\infty$ désigne un relevé 
quelconque de $a_i$. On a alors le lemme suivant qui établit un lien 
entre les anneaux $W_n^{\cris, \DP} (A^\infty, P^\infty)$ et $W_n^{\st, 
\DP} (A^\infty, P^\infty)$ :

\begin{lemme}
On garde les notations précédentes et on note $g \in P^\infty$ l'image 
de $x_0$. Soit $h \in P^\infty$ une racine $p^n$-ième de $g$. Alors 
l'application $W_n^{\cris, \DP} (A^\infty, P^\infty) \brac X \to 
W_n^{\st, \DP} (A^\infty, P^\infty)$ qui envoie $X$ sur $(1 \otimes (1 
\oplus h)) - 1$ est un isomorphisme.
\end{lemme}

\noindent
On en déduit directement la proposition :

\begin{prop}
\label{4:prop:descost}
Avec les notations précédentes, il existe un isomorphisme $S_n$-linéaire :
$$\xymatrix{\O^\st_n (A^\infty, P^\infty) \ar[r]^-{\sim} & W_n^{\cris, 
\DP} (A^\infty, P^\infty)} \brac X$$
où la structure de $S_n$-module est donnée sur $W_n^{\cris, \DP}
(A^\infty, P^\infty) \brac X$ par $\frac{u^i}{i!} \mapsto 
\frac{T_0^i}{i!} \frac 1 {1+X}$ où $T_0$ désigne une racine $p^n$-ième
de $[x_0]$. En outre, cet isomorphisme est compatible à la filtration
donnée à gauche par les $\J^{[s]}_n (A^\infty, P^\infty)$ et à droite
par les puissances divisées.
\end{prop}

\noindent
{\it Remarque.} Attention l'isomorphisme précédent n'est pas canonique : 
il dépend du choix d'une racine $p^n$-ième de $g$, image de $x_0$ dans 
$P^\infty$.

\subsubsection{Les opérateurs $\phi_s$ et $N$}
\label{4:subsec:ostphis}

La description précédente permet de prouver la proposition suivante 
importante :

\begin{prop}
Le faisceau $\O^\st_n$ est plat sur $S_n$ et les faisceaux $\J_n^{[s]}$ 
sont plats sur $W_n$.
\end{prop}

\begin{preuve}
L'argument est le même que celui de la proposition 2.1.2.1 de 
\cite{breuil-duke}.
\end{preuve}

\bigskip

Si $n$ et $m$ sont deux entiers avec $n \leq m$, on a un épaississement 
$i : T_n \toinj T_m$. Ainsi pour tout faisceau $\calF$ sur $(T_n)_\syn$, 
on peut former le faisceau $i_\star \calF$ sur $(T_m)_\syn$. Le foncteur 
$i_\star$ est exact (c'est une conséquence de la propriété 
\ref{4:prop:relev}) et, par abus, on note encore $\calF$ le faisceau 
$i_\star \calF$. L'exactitude assure qu'il revient au même de 
calculer les cohomologies de $\calF$ sur les sites $(T_n)_\syn$ et 
$(T_m)_\syn$.

Les descriptions locales données précédemment permettent facilement de 
prouver l'exactitude des deux suites de faisceaux (sur $(T_m)_\syn$) 
suivantes :
$$\xymatrix @C=25pt @R=5pt {
0 \ar[r] & \O^\st_i \ar[r]^-{p^n} & \O^\st_{n+i} \ar[r] & \O^\st_n \ar[r] & 0 
\\
0 \ar[r] & \J^{[s]}_i \ar[r]^-{p^n} & \J^{[s]}_{n+i} \ar[r] & \J^{[s]}_n 
\ar[r] & 0 }$$
pour tout entier $i$ tel que $n+i \leq m$. D'autre part, toujours pour 
$n+i \leq m$, la multiplication par $p^i$ identifie sur le site 
$(T_m)_\syn$ les faisceaux $\O^\st_n$ et $p^i \O^\st_{n+i}$. Comme la 
base $E_n$ est munie d'un relèvement du Frobenius (voir paragraphe 
\ref{4:subsec:bases}), les groupes $\O^\st_n(U) = H^0 ((U/E_n)_\cris, 
\O_{U/E_n})$ héritent d'un opérateur de Frobenius $\phi$ qui s'étend 
immédiatement en un morphisme de faisceaux $\phi : \O^\st_n \to 
\O^\st_n$.

Par ailleurs, on vérifie directement en utilisant la platitude et la 
description locale que pour tout entier $s$, on a l'inclusion $\phi 
(\J_n^{[s]}) \subset p^s \O^\st_n$. Les suites exactes précédentes 
permettent alors de définir un morphisme de faisceaux $\phi_s$ sur 
le site $(T_m)_\syn$ (pour $m \geq n+s$) comme la composée :
$$\xymatrix @C=20pt {
\J_n^{[s]} & \J_{n+s}^{[s]} / p^n \ar[l]_-{\sim} \ar[r]^-{\phi} & p^s 
\O^\st_{n+s} & \ar[l]_-{\sim}^-{p^s} \O^\st_n }$$
après avoir vérifié que $\phi$ passe au quotient.

\bigskip

\noindent
{\it Remarque.} Les faisceaux $\J_n^{[s]}$ et $\O^\st_n$ sont définis sur 
le site $(T_n)_\syn$ mais $\phi_s$ n'est, lui, défini que sur 
$(T_{n+s})_\syn$.

\bigskip

Finalement on peut munir $\O^\st_n$ d'un opérateur $N$ qui sur la 
description locale (voir proposition \ref{4:prop:descost}) est simplement 
défini comme l'unique application $W_n^{\cris, \DP} (A^\infty, 
P^\infty)$-linéaire qui envoie $\frac{X^i}{i!}$ sur $(1+X) 
\frac{X^{i-1}}{(i-1)!}$ (cette application ne dépend pas d'un choix
d'une racine $p^n$-ième de $g$, elle est donc canoniquement définie et
peut se recoller).

\subsection{Le cas de la caractéristique $p$}
\label{4:sec:osttilde}
Dans ce paragraphe, on se concentre sur le cas $n=1$ et on donne des 
descriptions plus faciles à manipuler des faisceaux précédemment 
introduits. En effet, par la suite, nous procéderons systématiquement 
par dévissages et donc le cas $n=1$ aura toujours un statut particulier.

\subsubsection{Une nouvelle description des faisceaux $\O^\st_1$ et 
$\J_1^{[s]}$}
\label{4:subsec:ost1local}

On reprend la description donnée par la proposition \ref{4:prop:descost} 
dans le cas $n=1$ :
$$\O^\st_1 (A^\infty, P^\infty) = W_1^{\cris, \DP} (A^\infty, P^\infty) 
\brac X.$$
Un calcul facile prouve que le morphisme qui envoie $\pi$ 
sur $[x_0]$ identifie $W_1^{\cris, \DP} (A^\infty, P^\infty)$ à :
$$\pa{\frac{k [Q^\infty G^{1/p}] [X_1^{1/p^\infty}, 
\ldots, X_s^{1/p^\infty}]}{([x_0]^e, f_1, \ldots, f_t)}}^\DP$$
où l'on nomme encore de façon abusive $f_i \in k [Q^\infty G^{1/p}] 
[X_1^{1/p^\infty}, \ldots, X_s^{1/p^\infty}]$ l'image de $f_i \in
\O_K/p [Q^\infty G] [X_1, \ldots, X_s]$ et où par définition $G^{1/p}$ 
désigne le sous-groupe de $(Q^\infty)^\gp$ formé des $x$ tels que $x^p 
\in G$. On remarque en outre que la suite $([x_0]^e, f_1, \ldots, f_t)$
est encore régulière dans $k [Q^\infty G^{1/p}] [X_1^{1/p^\infty},
\ldots, X_s^{1/p^\infty}]$.

Soit $\psi_i \in k [Q^\infty G] [X_1^{1/p^\infty}, \ldots, 
X_s^{1/p^\infty}]$ vérifiant $\psi_i^p = f_i$.
En explicitant les puissances divisées, on voit que $\O^\st_1 (A^\infty, 
P^\infty)$ s'identifie à :
\begin{equation}
\label{4:eq:o}
\bigoplus_{m_0, \ldots, m_{t+1} \in \N} \calB \cdot \gamma_{p m_0} 
(u^e) \gamma_{p m_1} (\psi_1) \cdots \gamma_{p m_t} (\psi_t) 
\gamma_{p m_{t+1}} (X)
\end{equation}
où $u = \frac{[x_0^{1/p}]} {1+X}$ et où on a posé :
\begin{equation}
\label{4:eq:b}
\calB = \frac{k [Q^\infty G^{1/p}] [X_1^{1/p^\infty}, \ldots, 
X_s^{1/p^\infty}][X]}{([x_0]^e, f_1, \ldots, f_t, X^p)}.
\end{equation}

\medskip

La description précédente fournit également une description locale 
des faisceaux $\J_1^{[s]}$ : $\J_1^{[s]} (A^\infty, P^\infty)$
s'identifie au sous-$\calB$-module de $\O^\st_1 (A^\infty, P^\infty)$
engendré par les :
$$\gamma_{m_0}(u^e) \gamma_{m_1}(\psi_1) \cdots \gamma_{m_t}(\psi_t) 
\gamma_{m_{t+1}} (X) \quad \text{avec} \quad m_0 + \cdots + m_t \geq
s.$$
Les quotients $\J_1^{[s]} / \J_1^{[s+1]}$ ont également une écriture 
sympathique :
\begin{equation}
\frac{\J_1^{[s]} (A^\infty, P^\infty)} {\J_1^{[s+1]} (A^\infty, 
P^\infty)} 
= \bigoplus_{\sum m_i = s} (\calB/(u^e, \psi_1, \ldots, 
\psi_t, X)) \cdot \gamma_{m_0} (u^e) \gamma_{m_1} (\psi_1) \cdots 
\gamma_{m_t} (\psi_t) \gamma_{m_{t+1}} (X).
\end{equation}

Si $s < p$, donc en particulier pour $s \leq r$, la dernière description 
se simplifie légèrement et donne :
\begin{equation}
\label{4:eq:js}
\frac{\J_1^{[s]} (A^\infty, P^\infty)} {\J_1^{[s+1]} (A^\infty, 
P^\infty)} =
\bigoplus_{\sum m_i = s} (\calB/(u^e, \psi_1, \ldots,
\psi_t, X)) \cdot u^{e m_0} \psi_1^{m_1} \cdots \psi_t^{m_t}
X^{m_{t+1}}
\end{equation}

Comme ces quotients ne sont pas tués par $u$, on introduit de nouveaux 
faisceaux intermédiaires :

\begin{deftn}
Soit $q$ un nombre rationnel de la forme $\frac t e$ où $t$ est un 
entier positif ou nul. Si $t = s e + \delta$ est la division euclidienne 
de $t$ par $e$, on pose :
$$\J_1^{[q]} = u^\delta \J_1^{[s]} + \J_1^{[s+1]}.$$
\end{deftn}

On peut à nouveau évaluer les quotients successifs comme le résume la 
proposition suivante :

\begin{prop}
Soit $0 \leq q < p$ un nombre rationnel de la forme $\frac t e$ pour un 
certain entier $t$. Notons $t = s e + \delta$ la division 
euclidienne de $t$ par $e$. Alors, en reprenant les notations 
précédentes :
\begin{equation}
\label{4:eq:j}
\frac{\J_1^{[q]} (A^\infty, P^\infty)} {\J_1^{[q+1/e]} (A^\infty,
P^\infty)} = \bigoplus_{\sum m_i = s} \calC \cdot u^{e m_0 + \delta} 
\psi_1^{m_1} \cdots \psi_t^{m_t} X^{m_{t+1}}
\end{equation}
où :
\begin{equation}
\label{4:eq:c}
\calC = \frac \calB {(u, \psi_1, \ldots, \psi_t, X)} = \frac{k 
[Q^\infty G^{1/p}] [X_1^{1/p^\infty}, \ldots, 
X_s^{1/p^\infty}]}{([x_0^{1/p}], \psi_1, \ldots, \psi_t)}.
\end{equation}
\end{prop}

\begin{preuve}
Il est clair que $\frac{\J_1^{[q]} (A^\infty, P^\infty)} {\J_1^{[q+1/e]}
(A^\infty, P^\infty)}$ s'identifie à 
$$\frac{u^\delta (\J_1^{[s]} (A^\infty, P^\infty) / \J_1^{[s+1]}
(A^\infty, P^\infty))}{u^{\delta+1} (\J_1^{[s]} (A^\infty, P^\infty) /
\J_1^{[s+1]} (A^\infty, P^\infty))}.$$
D'après la formule (\ref{4:eq:js}), il suffit de montrer que si l'on pose
$B = \calB/(u^e, \psi_1, \ldots, \psi_t, X)$, l'application de
multiplication par $u^\delta$ induit un isomorphisme de $\calC = B/uB$ 
dans $u^\delta B / u^{\delta+1} B$. La surjectivité est claire.

Montrons l'injectivité. Soit $B' = k [Q^\infty G^{1/p}]
[X_1^{1/p^\infty}, \ldots, X_s^{1/p^\infty}]$. Il suffit de montrer
que s'il existe un entier $\delta \leq e$ et des éléments $x, y, y_1,
\ldots, y_t$ de $B'$ tels que :
\begin{equation}
\label{4:eq:hyp}
u^\delta x = u^e y + y_1 \psi_1 + \cdots + y_t \psi_t
\end{equation}
alors, il existe $y'_1, \ldots, y'_t \in B'$ tels que :
$$x = u^{e-\delta} y + y'_1 \psi_1 + \cdots + y'_t \psi_t .$$
Nous allons prouver ce résultat par récurrence sur $t$. On rappelle
que la famille $(\psi_0, \psi_1, \ldots, \psi_t)$ est une suite
régulière de $B$ et que, dans $B$, on a $u^e = \psi_0$ (puisque $X$
est nul). 

L'initialisation de la récurrence provient simplement du fait que $u$
n'est pas diviseur de $0$ dans $B$ (puisque $u^e$ ne l'est pas). Pour
l'hérédité, on suppose que l'équation (\ref{4:eq:hyp}) est vérifiée.
Cela entraîne $u^\delta x \equiv y_t \psi_t \pmod{u^e, \psi_1, \ldots,
\psi_{t-1}}$ et donc $u^{e-n\delta} y_t \psi_t \equiv 0 \pmod{u^e,
\psi_1, \ldots, \psi_{t-1}}$. Puisque la suite $(u^e, \psi_1, \ldots, 
\psi_t)$ est régulière, il vient $u^{e-\delta} y_t \equiv 0
\pmod{u^e, \psi_1, \ldots, \psi_{t-1}}$, et donc il existe $z,
z_1, \ldots, z_{t-1} \in B$ tels que :
$$u^{e-\delta} y_t = z u^e + z_1 \psi_1 + \cdots + z_{t-1} \psi_{t-1}
.$$
En appliquant l'hypothèse de récurrence, on obtient l'existence
d'éléments $z', z'_1, \ldots, z'_{t-1}$ dans $B$ tels que :
$$y_t = z' u^\delta + z'_1 \psi_1 + \cdots + z'_{t-1} \psi_{t-1} .$$
En réinjectant dans (\ref{4:eq:hyp}), il vient :
$$u^\delta (x - z' y_t) = u^e y + (y_1 + z'_1 y_t) \psi_1 + \cdots +
(y_{t-1} + z'_1 y_{t-1}) \psi_{t-1}$$
et une nouvelle application de l'hypothèse de récurrence permet de
conclure.
\end{preuve}

\bigskip

\noindent
{\it Remarque.} Sur les quotients de la proposition qui précède, la 
structure de $A^\infty$-module est simple à décrire. En effet, ces 
quotients sont tués par $u$ et par tous les $\psi_i$, et donc la 
structure de $A^\infty$-module se factorise en une structure de 
$\calC$-module, qui est décrite de façon transparente sur la somme 
directe précédente.

\bigskip

Finalement, remarquons que si $s \leq p-1$ et si $p \geq 3$, il est 
possible de donner une définition alternative de $\phi_s : \J^{[s]}_1 
\to \O^\st_1$ qui est déjà valable sur le site $(T_2)_\syn$. On remarque 
pour cela que $\phi_1$ s'annule sur $\J^{[2]}_1$ et qu'il définit ainsi 
par passage au quotient un morphisme de faisceaux (sur le site 
$(T_2)_\syn$) $\J^{[1]}_1 / \J^{[2]}_1 \to \O^\st_1$. En outre, on a
le lemme suivant :

\begin{lemme}
Le morphisme de faisceaux canoniques $\sym_{\O_1} \J^{[1]}_1 /
\J^{[2]}_1 \to \J^{[s]}_1 / \J^{[s+1]}_1$ est un isomorphisme.
\end{lemme}

\begin{preuve}
Avec les descriptions précédentes, c'est une conséquence directe de
l'alinéa I.3.4.4 de \cite{berthelot}.
\end{preuve}

\bigskip

On vérifie que l'application $\phi_s$ s'obtient comme la composée :
$$\xymatrix @C=40pt {
\J_1^{[s]} \ar[r] & \J_1^{[s]} / \J_1^{[s+1]} \ar[r]^-{\sym^s \phi_1} &
\O^\st_1 }. $$
Cette dernière formule assure, au moins pour $s \leq p-1$ et $p \geq 3$,
que $\phi_s$ ne dépend que de $\phi_1$ et peut être défini sur 
$(T_2)_\syn$.

\subsubsection{Les faisceaux $\protect\tilde{\O}^\st$, $\protect\tilde 
{\J}^{[q]}$}

Lorsque $n=1$, on définit des versions simplifiées des faisceaux 
$\O^\st_1$ et $\J_1^{[s]}$ en remplaçant la base $E_1$ par la base $\tilde 
E$ (voir le paragraphe \ref{4:subsec:bases} pour la définition).

\bigskip

Si $U$ est un log-schéma sur $T_1$, on définit pour tout entier $s$ :
$$\tilde \J^{[s]} (U) = H^0 ((U/\tilde E)_\cris, \J^{[s]}_{U/\tilde E}) 
= H^0 ((U/\tilde E)_\CRIS, \J^{[s]}_{U/\tilde E})$$
et $\tilde \O^\st = \tilde \J^{[0]}$. Comme précédemment, les 
préfaisceaux $\tilde \J^{[s]}$ sont des faisceaux sur le gros site 
syntomique $(T_1)_\SYN$ et calculent la cohomologie log-cristalline :
$$H^i(X_\syn, \tilde \J^{[s]}) = H^i(X_\SYN, \tilde \J^{[s]}) = H^i
((X/\tilde E)_\cris, \tilde \J^{[s]}_{X/\tilde E}) = 
H^i((X/\tilde E)_\CRIS, \J^{[s]}_{X/\tilde E}).$$
De même que précédemment, on pose pour $q = s + \frac \delta e$ où 
$s \geq 0$ et $0 \leq \delta < e$ sont des entiers :
$$\tilde \J^{[q]} = u^\delta \tilde \J^{[s]} + \tilde \J^{[s+1]}.$$

Il est encore possible de donner une description locale très 
explicite de ces faisceaux. En reprenant les arguments utilisés pour
les faisceaux $\O^\st_1$ et $\J_1^{[s]}$, et en reprenant les notations 
des paragraphes \ref{4:subsec:ostlocal} et \ref{4:subsec:ost1local}, on 
aboutit à :
\begin{equation}
\label{4:eq:otilde}
\tilde \O^\st (A^\infty, P^\infty) = \bigoplus_{m_1, \ldots, m_{t+1} \in 
\N} \tilde \calB \cdot \gamma_{p m_1} (\psi_1) \cdots \gamma_{p m_t} 
(\psi_t) \gamma_{p m_{t+1}} (X)
\end{equation}
avec :
\begin{equation}
\label{4:eq:btilde}
\tilde \calB = \frac{k [Q^\infty G^{1/p}] [X_1^{1/p^\infty}, \ldots, 
X_s^{1/p^\infty}][X]}{([x_0], f_1, \ldots, f_t, X^p)} = \frac{k 
[Q^\infty G^{1/p}] [X_1^{1/p^\infty}, \ldots, 
X_s^{1/p^\infty}][X]}{(u^p, f_1, \ldots, f_t, X^p)}.
\end{equation}
et, pour $q = s + \frac \delta e < \frac p e$ (par exemple si $q \leq 
r$) :
\begin{equation}
\label{4:eq:jtilde}
\frac{\tilde \J^{[q]} (A^\infty, P^\infty)} {\tilde \J^{[q+1/e]} 
(A^\infty, P^\infty)} = \bigoplus_{\sum m_i = s} \tilde \calC \cdot 
u^{e m_0 + \delta} \psi_1^{m_1} \cdots \psi_t^{m_t} X^{m_{t+1}}
\end{equation}
où $\tilde \calC = \calC$. La structure de $A^\infty$-module sur ce 
dernier quotient se factorise en une structure de $\tilde \calC$-module, 
qui est celle qui correspond à l'écriture sous forme de somme directe.

\medskip

En comparant les formules (\ref{4:eq:j}) et (\ref{4:eq:jtilde}), on en
déduit la proposition suivante :

\begin{prop}
\label{4:prop:quotj}
Pour $q \in \frac 1 e \N$, $q < \frac p e$, la projection canonique 
induit un isomorphisme de faisceaux :
$$\xymatrix @C=25pt {
\J_1^{[q]} / \J_1^{[q+1/e]} \ar[r]^-{\sim} & \xytilde \J^{[q]} / 
\xytilde \J^{[q+1/e]} }.$$
\end{prop}

Par ailleurs, pour tout entier $s < \frac p e$, l'opérateur $\phi_s : 
\J_1^{[s]} \to \O^\st_1$ induit par passage au quotient un morphisme 
$\phi_s : \J_1^{[s]} / \J_1^{[s+1/e]} \to \O^\st_1 / u^p \O^\st_1$. 
D'autre part, la projection canonique $\O^\st_1 \to \tilde \O^\st$ 
s'annule sur $u^p \O^\st_1$ et donc fournit un morphisme de faisceaux 
$\O^\st_1 / u^p \O^\st_1 \to \tilde \O^\st$ (attention, ce n'est pas un 
isomorphisme). Tout cela, avec la proposition précédente, permet de 
définir un opérateur $\tilde \phi_s$ comme la composée :
$$\xymatrix @C=25pt {
\xytilde \J^{[s]} \ar[r] & \xytilde \J^{[s]} / \xytilde \J^{[s+1/e]} 
& \J_1^{[s]} / \J_1^{[s+1/e]} \ar[l]_-{\sim} \ar[r]^-{\phi_s} & 
\O^\st_1 / u^p \O^\st_1 \ar[r] & \xytilde \O^\st }. $$
Encore une fois, le morphisme $\tilde \phi_s$ n'est pas 
défini sur $(T_1)_\syn$, mais \emph{a priori} simplement sur 
$(T_{s+1})_\syn$.

\bigskip

Finalement, on munit le faisceau $\tilde \O^\st$ d'un opérateur $\tilde 
N$, après avoir remarqué que l'opérateur $N : \O^\st_1 \to \O^\st_1$ 
passe au quotient puisque $N(u^p \O^\st_1 + \sum_{i \geq p} \gamma_i (u) 
\O^\st_1) \subset u^p \O^\st_1 + \sum_{i \geq p} \gamma_i (u) \O^\st_1$.

\subsubsection{Les faisceaux $\protect\bar{\O}^\st$ et $\protect\bar 
{\J}^{[s]}$}

On recopie les définitions précédentes en se plaçant 
désormais sur la base $\bar E$ (voir paragraphe \ref{4:subsec:bases}) :
pour tout entier $s$ et tout $U \in (T_1)_\SYN$, on pose
$$\bar \J^{[s]} (U) = H^0 ((\bar U/\bar E)_\cris, \J^{[s]}_{\bar U/\bar 
E}) = H^0 ((\bar U/\bar E)_\CRIS, \J^{[s]}_{\bar U/\bar E})$$
où $\bar U = U \times_{T_1} \bar E$. On obtient ainsi des faisceaux sur 
le site $(T_1)_\SYN$ qui, comme précédemment, calculent la cohomologie 
log-cristalline. On pose en outre $\bar \O^\st = \bar \J^{[0]}$.

On dispose encore d'une description locale du faisceau 
$\bar \O^\st$. En gardant les mêmes notations, on a :
\begin{equation}
\label{4:eq:obar}
\bar \O^\st (A^\infty, P^\infty) = \bigoplus_{m_1, \ldots, m_{t+1} \in 
\N} \bar \calB \cdot \gamma_{p m_1} (\psi_1) \cdots \gamma_{p m_t} 
(\psi_t) \gamma_{p m_{t+1}} (X)
\end{equation}
avec :
\begin{equation}
\label{4:eq:bbar}
\bar \calB = \frac{k [Q^\infty G^{1/p}] [X_1^{1/p^\infty}, \ldots, 
X_s^{1/p^\infty}][X]}{(u, f_1, \ldots, f_t, X^p)}.
\end{equation}

L'élément $u$ est nul par hypothèse dans $\bar E$ et donc également dans 
$\bar \O^\st$. On ne définit donc pas les faisceaux intermédiaires $\bar 
\J^{[q]}$ pour $q$ rationnel, mais, si $s$ est un entier, on a toujours 
une description explicite du quotient :
\begin{equation}
\label{4:eq:jbar}
\frac{\bar \J^{[s]} (A^\infty, P^\infty)} {\bar \J^{[s+1]} (A^\infty, 
P^\infty)} = 
\bigoplus_{\sum m_i = s} \bar \calC \cdot \gamma_{m_1} (\psi_1) \cdots 
\gamma_{m_t} (\psi_t) \gamma_{m_{t+1}} (X)
\end{equation}
où $\bar \calC = \tilde \calC = \calC$. Il est encore une fois possible 
de décrire la structure de $A^\infty$-module sur ce quotient : elle se
factorise par la structure de $\bar \calC$-module naturelle sur la somme 
directe. Finalement, en comparant les formules (\ref{4:eq:j}) et 
(\ref{4:eq:jbar}), on voit que le morphisme naturel $\J_1^{[s]} 
/ \J_1^{[s+1/e]} \to \bar \J^{[s]} / \bar \J^{[s+1]}$ est un 
isomorphisme pour $s < p$. Cela permet de définir un opérateur $\bar 
\phi_s$ comme la composée :
$$\xymatrix @C=25pt {
\xybar \J^{[s]} \ar[r] & \xybar \J^{[s]} / \xybar \J^{[s+1]} 
& \J_1^{[s]} / \J_1^{[s+1/e]} \ar[l]_-{\sim} \ar[r]^-{\phi_s} & 
\O^\st_1 / u^p \O^\st_1 \ar[r] & \xybar \O^\st }. $$

\subsubsection{Les faisceaux $\O^\car_1$, $\protect\tilde{\O}^\car$ et 
$\protect\bar{\O}^\car$}
\label{4:subsec:defocar}

Sur un log-schéma de caractéristique $p$, on définit le Frobenius absolu 
de la façon suivante : c'est le Frobenius absolu classique sur le schéma
sous-jacent et la multiplication par $p$ sur le monoïde (à supposer 
qu'il soit noté additivement).

Soit $U$ un log-schéma fin (localement de présentation finie) sur 
$T_1$. On voit $U$ sur $E_1$ grâce à l'épaississement $T_1 \toinj E_1$. 
On note $U' = U \times_{E_1} E_1$ où $E_1$ est vu sur lui-même par le 
morphisme de Frobenius absolu. Le Frobenius absolu sur $U$ se factorise 
par $U'$ et fournit donc un morphisme $U \to U'$ appelé Frobenius 
relatif $f_U$.

Dans la terminologie de Kato (voir \cite{kato}, paragraphe 4.9), la 
flèche $f_U$ est faiblement purement inséparable et d'après le théorème 
4.10 de \emph{loc. cit.}, elle se factorise de façon \emph{unique} sous 
la forme :
$$f_U : U \to U'' \to U'$$
où le premier morphisme est purement inséparable et le second est 
log-étale. On définit alors, comme dans le paragraphe 2.2.1 de 
\cite{breuil-duke} :
$$\O^\car_1 (U) = \Gamma(U'', \O_{U''}).$$

On prouve (voir appendice B de \cite{breuil-duke}) que l'on définit 
ainsi un faisceau $\O^\car_1$ sur 
le gros site syntomique $(T_1)_\SYN$. Par ailleurs, d'après les 
résultats de \cite{kato}, si $U = (\spec A, P)$, on a simplement :
$$\O^\car_1 (U) = (S_1 \otimes_{(\phi), k} A/\pi A) \otimes_{\Z[\N 
\oplus_{(\phi), \N} P]} \Z[(\N \oplus_{(\phi), \N}, P)^{(1)}].$$

Avec les notations introduites précédemment, ceci se réécrit :
\begin{equation}
\label{4:eq:car}
\O^\car_1 (A^\infty, P^\infty) = \bigoplus_{m_0 \in \N} \calB \cdot 
\gamma_{p m_0} (u^e)
\end{equation}
où on rappelle que $\calB$ est défini par la formule (\ref{4:eq:b}).

\bigskip

On définit de même sur le site $(T_1)_\SYN$ les faisceaux $\tilde 
\O^\car$ et $\bar \O^\car$ en remplaçant la base $E_1$ respectivement
par les bases $\tilde E$ et $\bar E$. Comme précédemment, on peut
donner une description explicite des faisceaux obtenus.
On obtient :
\begin{equation}
\label{4:eq:cartb}
\tilde \O^\car (A^\infty, P^\infty) = \tilde \calB 
\qquad \text{et} \qquad
\bar \O^\car (A^\infty, P^\infty) = \bar \calB
\end{equation}
où les anneaux $\tilde \calB$ et $\bar \calB$ sont définis 
respectivement par les formules (\ref{4:eq:btilde}) et (\ref{4:eq:bbar}).

\medskip

Pour finir, mentionnons que la structure de $A^\infty$-module sur 
l'objet $\tilde \O^\car (A^\infty, P^\infty)$ est donnée par 
l'application $A^\infty \to \tilde \calB$ déduite du Frobenius.
En particulier, on constate que cette application se factorise par 
$\calC = \tilde \calC = \bar \calC$.
Il en va de même pour $\bar \O^\car (A^\infty, P^\infty)$.

\section{Calcul de la cohomologie cristalline}
\label{4:sec:cris}

On montre dans cette partie comment associer à un log-schéma $X$
propre et log-lisse sur la base $T = (\N \to \O_K, 1 \mapsto \pi)$ dont 
la réduction modulo $p$ est du \emph{type de Cartier} (voir définition 
4.8 de \cite{kato}) sur $X \times_T T_1$, et à un entier $n$ des objets 
des catégories $\Mr$ tués par $p^n$. Précisément, il s'agit, pour tout 
$i < r$, des quadruplets :
$$(H^i ((X_n)_\syn, \O^\st_n), H^i ((X_n)_\syn, \J^{[r]}_n), \phi_r, 
N)$$
où par définition $X_n = X \times_T T_n$. On fait remarquer une fois de 
plus que le morphisme $\phi_r$ n'est défini que sur $(X_{n+r})_\syn$. 
Toutefois, l'écriture précédente du quadruplet est légitime car on 
dispose d'isomorphismes canoniques entre $(H^i ((X_n)_\syn, \O^\st_n)$
et $(H^i ((X_{n+r})_\syn, \O^\st_n)$ d'une part et $H^i ((X_n)_\syn, 
\J^{[r]}_n)$ et $H^i ((X_{n+r})_\syn, \J^{[r]}_n)$ d'autre part.

\medskip

On commence par traiter le cas $n=1$ (sous-partie \ref{4:sec:crisp}) : 
alors, l'élément de $\Mr$ est tué par $p$ et peut-être donc vu dans 
$\Mrtilde$. On procède ensuite par dévissage avant de passer un passage 
à la limite projective pour proposer une version entière (sous-partie 
\ref{4:sec:crisdev}).

\subsection{En caractéristique $p$}
\label{4:sec:crisp}

On se donne ici, un log-schéma fin $X_1$ propre et log-lisse sur la base 
$T_1$ (voir paragraphe \ref{4:subsec:bases}). En particulier, le morphisme 
structural $X_1 \to T_1$ est log-syntomique. On suppose que $X_1$ est du 
type de Cartier et qu'il admet un relèvement $X_2$ fin et log-lisse sur 
$T_2$. Ce relèvement est automatiquement log-syntomique sur $T_2$.

\medskip

Le but, ici, est de prouver que le quadruplet $(H^i ((X_1)_\syn,
\O^\st_1), H^i ((X_1)_\syn, \J^{[r]}_1), \phi_r, N)$ définit un objet de
la catégorie $\Mr$ pour tout $i \leq r$. On montre d'abord l'énoncé 
équivalent avec la base $\tilde E$, \emph{i.e.} que le quadruplet $(H^i 
((X_1)_\syn, \tilde \O^\st), H^i ((X_1)_\syn, \tilde \J^{[r]}), \tilde 
\phi_r, \tilde N)$ définit un objet de $\Mrtilde$. Le point le plus 
difficile est la liberté du $\tilde S$-module $H^i ((X_1)_\syn, \tilde 
\O^\st)$. Les deux paragraphes \ref{4:subsec:isomfais} et 
\ref{4:subsec:isomcoh} y sont consacrés. Dans le paragraphe 
\ref{4:subsec:cristilde}, on explique comment on termine la preuve pour la 
base $\tilde E$ avant d'en déduire dans le paragraphe 
\ref{4:subsec:crise1}, le théorème sur la base $E_1$.

\subsubsection{Des isomorphismes sur les faisceaux} 
\label{4:subsec:isomfais}
On suit pratiquement à la lettre la méthode initiée par Fontaine et 
Messing (\cite{fontaine-messing}) et développée par Breuil 
(\cite{breuil-duke}) dans le cas qui nous intéresse. Pour tout entier $s 
\geq 1$, on note $\J^{\brac s}_1$ le noyau du morphisme composé 
$\xymatrix @C=15pt {\O^\st_{s+1} \ar[r]^-{\phi} & \O^\st_{s+1} \ar[r] 
& \O^\st_s}$. On note $\nu_s : \J^{\brac s}_1 \to \O^\st_1$ la réduction 
modulo $p$ et on définit $\hat f_s : \J^{\brac s}_1 \to \O^\st_1$ par 
$\hat f_s (x) = y$ si $\phi(x) = p^s \hat y$ où $\hat y$ est une section 
(locale) du faisceau $\O^\st_{s+1}$.

On définit ensuite $F^s \O^\st_1 = \im \nu_s$ et $F_s \O^\st_1 = \im 
\hat f_s$. Ce sont deux sous-faisceaux d'anneaux de $\O^\st_1$. On 
montre facilement que la suite $(F^s \O^\st_1)$ est décroissante alors 
que la suite $(F_s \O^\st_1)$ est croissante. En outre, le morphisme 
$\hat f_s$ se factorise en un isomorphisme :
$$\xymatrix {f_s : F^s \O^\st_1 / F^{s+1} \O^\st_1 \ar[r]^-{\sim} &
F_s \O^\st_1 / F_{s-1} \O^\st_1}.$$

On note $F_s^\car \O^\st_1$ le sous-faisceau de $\O^\car_1$-algèbre de 
$\O^\st_1$ engendré par $F^s \O^\st_1$ et $F_s^\car \tilde \O^\st$ 
(resp. $F_s^\car \bar \O^\st$) l'image de $F_s^\car \O^\st_1$ 
dans $\tilde \O^\st$ (resp. dans $\bar \O^\st$). On note également $F^s 
\bar \O^\st$ la réduction de $F^s \O^\st_1$ dans $\bar \O^\st$.

\medskip

Soient $\O_1$ le faisceau structural sur le site $(T_1)_\SYN$, et $\bar 
\O$ sa réduction modulo $\pi$. On a, en reprenant les notations du 
paragraphe \ref{4:sec:crisp}, les descriptions locales suivantes :
\begin{equation}
\O_1 (A^\infty, P^\infty) = A^\infty
\qquad \text{et} \qquad
\bar \O (A^\infty, P^\infty) = A^\infty / \pi.
\end{equation}
On dispose de flèches naturelles $\O_1 \to \O_1^\car$, $\bar \O \to 
\tilde \O^\car$ et $\bar \O \to \bar \O^\car$ qui localement sur les 
descriptions précédentes sont données par l'élévation à la puissance 
$p$. Ces flèches se factorisent toutes les trois par $\calC$.

\bigskip

\begin{prop}
\label{4:prop:fscar}
\begin{enumerate}
\item[i)] Pour tout entier $s$, on a $F^s \bar \O^\st = \bar \J^{[s]}$.

\item[ii)] On a $\bigcup_{s \in \N} F_s^\car \O^\st_1 = \O^\st_1$, 
$\bigcup_{s \in \N} F_s^\car \tilde \O^\st = \tilde \O^\st$ et 
$\bigcup_{s \in \N} F_s^\car \bar \O^\st = \bar \O^\st$.

\item[iii)] Pour tout entier $s$, l'isomorphisme $f_s$ se factorise en 
des isomorphismes de faisceaux sur $(\bar E)_\syn$ :
$$\xymatrix @C=40pt @R=8pt {
\O^\car_1 \otimes_{\O_1} \xybar \J^{[s]} / \xybar \J^{[s+1]} \ar[r]^-{\sim} 
& F_s^\car \O^\st_1 / F_{s-1}^\car \O^\st_1 \\
\xytilde \O^\car \otimes_{\O_1} \xybar \J^{[s]} / \xybar \J^{[s+1]} 
\ar[r]^-{\sim} & F_s^\car \xytilde \O^\st / F_{s-1}^\car \xytilde \O^\st \\
\xybar \O^\car \otimes_{\O_1} \xybar \J^{[s]} / \xybar \J^{[s+1]} 
\ar[r]^-{\sim} & F_s^\car \xybar \O^\st / F_{s-1}^\car \xybar \O^\st }$$
\end{enumerate}
\end{prop}

\begin{preuve}
Elle est entièrement analogue à celle de la proposition 2.2.2.2 de 
\cite{breuil-duke}. Signalons toutefois une subtilité peut-être 
insuffisamment soulignée dans \emph{loc. cit.} : une fois 
prouvés i) et ii), on est amené à évaluer, pour le premier isomorphisme 
de iii), le produit tensoriel :
$$\O^\car_1(A^\infty, P^\infty) \otimes_{A^\infty} \frac{\bar \J^{[s]} 
(A^\infty, P^\infty)} {\bar \J^{[s+1]} (A^\infty, P^\infty)}$$
qui d'après les descriptions précédentes s'identifie à :
$$\O^\car_1(A^\infty, P^\infty) \otimes_{A^\infty} \pa{ 
\bigoplus_{\sum m_i = s} \bar \calC \cdot \gamma_{m_1} (\psi_1) \cdots 
\gamma_{m_t}(\psi_t) \gamma_{m_{t+1}} (X) }.$$
Il faut alors se souvenir que chacun des facteurs du produit tensoriel 
est tué par $u$ et tous les $\psi_i$, de sorte que la structure de 
$A^\infty$-module sur ces deux facteurs se factorise en une structure de 
$\bar \calC$-module. Ayant constaté cela, le produit tensoriel se 
réécrit :
$$\O^\car_1(A^\infty, P^\infty) \otimes_{\bar \calC} \pa{ 
\bigoplus_{\sum m_i = s} \bar \calC \cdot \gamma_{m_1} (\psi_1) \cdots 
\gamma_{m_t}(\psi_t) \gamma_{m_{t+1}} (X) }$$
et donc vaut bien :
$$\bigoplus_{\sum m_i = s} \O^\car_1(A^\infty, P^\infty) \cdot 
\gamma_{m_1} (\psi_1) \cdots \gamma_{m_t}(\psi_t) \gamma_{m_{t+1}} (X)$$
comme annoncé dans \emph{loc. cit.} La fin de la preuve reste inchangée.
\end{preuve}

\bigskip

Notons qu'au passage la preuve de \cite{breuil-duke} donne des 
descriptions locales explicites pour les faisceaux $F_s^\car$, qui sont 
(toujours en gardant les mêmes notations) :
\begin{eqnarray*}
F_s^\car \O^\st_1 & = & \bigoplus_{\substack{m_0 \in \N \\ m_i + \cdots 
+ m_{t+1} \leq s}} \calB \cdot \gamma_{p m_0} (u^e) \gamma_{p m_1} 
(\psi_1) \cdots \gamma_{p m_t} (\psi_t) \gamma_{p m_{t+1}} (X) \\
F_s^\car \tilde \O^\st & = & \bigoplus_{m_1 + \cdots + m_{t+1} \leq s} 
\tilde \calB \cdot \gamma_{p m_1} (\psi_1) \cdots \gamma_{p m_t} 
(\psi_t) \gamma_{p m_{t+1}} (X) \\
F_s^\car \bar \O^\st & = & \bigoplus_{m_1 + \cdots + m_{t+1} \leq s} 
\bar \calB \cdot \gamma_{p m_1} (\psi_1) \cdots \gamma_{p m_t} 
(\psi_t) \gamma_{p m_{t+1}} (X).
\end{eqnarray*}

\bigskip

De ces descriptions explicites, on déduit le théorème suivant :

\begin{theo}
\label{4:th:isomfais}
Soit $s < \frac p e$ un entier (si $p=2$ et $e=1$, on impose $s=0$). On 
a alors des isomorphismes de faisceaux sur le site $(T_2)_\syn$ :
$$\xymatrix @C=40pt @R=8pt {
\xytilde \O^\car \otimes_{\xybar \O} \J_1^{[s]} / \J_1^{[s+1/e]} =
\xytilde \O^\car \otimes_{\xybar \O} \xytilde \J^{[s]} / \xytilde \J^{[s+1/e]} 
\ar[r]^-{\sim}_-{\id \otimes \phi_s} & F_s^\car \xytilde \O^\st \\
\xybar \O^\car \otimes_{\xybar \O} \J_1^{[s]} / \J_1^{[s+1/e]} =
\xybar \O^\car \otimes_{\xybar \O} \xytilde \J^{[s]} / \xytilde \J^{[s+1/e]} 
\ar[r]^-{\sim}_-{\id \otimes \phi_s} & F_s^\car \xybar \O^\st. }$$
\end{theo}

\begin{preuve}
Encore une fois, la preuve est identique à celle du théorème 2.2.2.3 de 
\cite{breuil-duke}. Notons toutefois que la subtilité mentionnée dans la 
preuve de la proposition précédente apparaît à nouveau ici.
\end{preuve}

\bigskip

\noindent
{\it Remarque.} On démontre de façon tout à fait identique que l'on a 
également les isomorphismes suivants, pour $s \leq \frac p e - 1$ :
$$\xymatrix @C=40pt @R=8pt {
\O_1^\car \otimes_{\O_1} \J_1^{[s]} / \J_1^{[s+1]}
\ar[r]^-{\sim}_-{\id \otimes \phi_s} & F_s^\car \O_1^\st }$$
$$\xymatrix @C=40pt @R=10pt {
\xytilde \O^\car \otimes_{\O_1} \J_1^{[s]} / \J_1^{[s+1]} =
\xytilde \O^\car \otimes_{\O_1} \xytilde \J^{[s]} / \xytilde \J^{[s+1]} 
\ar[r]^-{\sim}_-{\id \otimes \phi_s} & F_s^\car \xytilde \O^\st \\
\xybar \O^\car \otimes_{\O_1} \J_1^{[s]} / \J_1^{[s+1]} =
\xybar \O^\car \otimes_{\O_1} \xytilde \J^{[s]} / \xytilde \J^{[s+1]} 
\ar[r]^-{\sim}_-{\id \otimes \phi_s} & F_s^\car \xybar \O^\st. }$$
certainement plus proches de ceux de \cite{breuil-duke} (rappelons que 
$\tilde \O^\car$ et $\bar \O^\car$ sont tués par $\pi$).

\medskip

Par ailleurs, si l'on ne se préoccupe que de la version \og $\J_1$ \fg\ 
et pas de la version \og $\tilde \J$ \fg, les isomorphismes précédents
sont valables pour tout $s \leq p-1$.

\subsubsection{Des isomorphismes sur les groupes de cohomologie}
\label{4:subsec:isomcoh}
Nous aimerions à présent déduire du théorème \ref{4:th:isomfais} des 
isomorphismes sur les groupes de cohomologie et, pour cela, nous allons
projeter ces faisceaux sur le site étale : l'intérêt est que sur ce
site le faisceau $\tilde \O^\car$ (resp. $\bar \O^\car$) se réduit 
simplement à $\tilde S \otimes_{(\phi), k} \bar \O$ (resp. à $\bar S 
\otimes_{(\phi), k} \bar \O$).

\medskip

Soit $X$ un log-schéma fin localement de type fini sur une des trois 
bases $E_1$, $\tilde E$ ou $\bar E$. On dispose dans ces conditions d'un 
morphisme de topoï :
$$\alpha : \widetilde{X_\SYN} \to \widetilde{X_\et}$$
défini de la façon suivante : si $\calF$ est un faisceau sur $X_\SYN$, 
on définit $\alpha_\star \calF$ comme la restriction de $\calF$ au site 
$X_\et$, et on vérifie que l'on obtient ainsi un faisceau pour la 
topologie étale. Réciproquement si $\calF$ est un faisceau sur $X_\ET$, 
on définit $\alpha^\star \calF$ comme le faisceau associé (pour la 
topologie syntomique) au préfaisceau $\calF$.

\medskip

Si $\Gamma_\SYN$ (resp. $\Gamma_\ET$) désigne le foncteur des sections 
globales pour la topologie syntomique (resp. étale), on a évidemment la 
relation $\Gamma_\SYN = \Gamma_\ET \circ \alpha_\star$ d'où 
$R\Gamma_\SYN = R\Gamma_\ET \circ R\alpha_\star$. Ainsi pour calculer la 
cohomologie syntomique d'un faisceau, il suffit de calculer le 
$R\alpha_\star$ de ce faisceau puis de déterminer l'hypercohomologie 
étale du complexe obtenu. C'est ce que nous allons faire.

\paragraph{Calcul des $R\alpha_\star$ de plusieurs faisceaux}

Pour calculer les $R\alpha_\star$ des faisceaux précédemment introduits, 
on aimerait utiliser les résolutions de Berthelot et Kato. Cependant, 
celles-ci
sont valables sur le site cristallin et non sur le site syntomique. Il 
nous faut donc faire un pont entre cohomologie cristalline et 
cohomologie syntomique, pont qui passe par la cohomologie 
cristalline-syntomique.

\medskip

Soit $X$ un log-schéma fin localement de type fini sur une des trois
bases $E_1$, $\tilde E$ ou $\bar E$. Reprenant les notations de 
\cite{breuil-duke} auquel on aura besoin de se référer par la suite, 
on note $\Upsilon$ la base retenue.

On définit le \emph{site cristallin-syntomique} sur $X/\Upsilon$ en 
munissant la catégorie sous-jacente au site cristallin sur $X/\Upsilon$ 
de la topologie syntomique : il s'agit donc d'une catégorie de couples 
$(U \toinj T)$ et on convient qu'une famille de couples $(U_i \toinj 
T_i)$ recouvre $(U \toinj T)$ si les $T_i$ forment un recouvrement 
syntomique de $T$ et si les diagrammes :
$$\xymatrix { U_i \ar[r] \ar[d] & T_i \ar[d] \\ U \ar[r] & T }$$
sont cartésiens.

Bien entendu, selon que l'on considère la catégorie sous-jacente au
petit site cristallin sur $X/\Upsilon$ ou au grand, on obtient
respectivement les petit et grand sites cristallin-syntomiques sur
$X/\Upsilon$. On les note $(X/\Upsilon)_\syncris$ et
$(X/\Upsilon)_\SYNCRIS$. Par les résultats de \cite{breuil-bull} 
(lemme 3.3.1), on a des morphismes de topoï entre les différents 
catégories de faisceaux sur les sites précédents comme le résume le 
carré commutatif suivant :
$$\xymatrix @C=50pt @R=30pt {
(\widetilde {X/\Upsilon})_\SYNCRIS \ar[r]^-{v} \ar[d]_-{w}
& (\widetilde {X/\Upsilon})_\CRIS \ar[d]^-{u} \\
\widetilde{X_\SYN} \ar[r]^-{\alpha} & \widetilde{X_\et} }$$

Par ailleurs, il est possible de définir, comme en \ref{4:subsec:defocar},
sur le site $(X/\Upsilon)_\SYNCRIS$ (resp $(X/\Upsilon)_\syncris$) des
faisceaux $\O_1^\car$, $\tilde \O^\car$ et $\bar \O^\car$ en posant
$\O_1^\car (U \toinj T) = \O_1^\car (U)$ (où le deuxième $\O_1^\car$ est 
celui défini précédemment\footnote{Notez que si $U \toinj T$ est un 
objet de $(X/\Upsilon)_\SYNCRIS$, alors $U \to X$ est étale et donc en 
particulier log-syntomique.}) et des formules analogues pour les autres
faisceaux. (Attention, dans \cite{breuil-duke}, ces faisceaux sont notés 
respectivement $\O_X^{\car, E_1}$, $\O_X^{\car, \tilde E_1}$ et 
$\O_X^{\car, \bar E_1}$.) De même, ces faisceaux vivent également sur 
les sites $(X/\Upsilon)_\CRIS$ (resp. $(X/\Upsilon)_\cris$) et $X_\ET$ 
(resp. $X_\et$).

\bigskip

Soit $\calF$ un faisceau sur un des sites $(X/\Upsilon)_\CRIS$, 
$(X/\Upsilon)_\cris$, $(X/\Upsilon)_\SYNCRIS$ ou 
$(X/\Upsilon)_\syncris$. Pour tout $T$ apparaissant dans un couple
de la forme $(U \toinj T)$, on sait que $\calF$ définit un faisceau
$F_{|U \toinj T}$ sur $T_\et$. On dit que $\calF$ est à 
\emph{composantes quasi-cohérentes} si tous les faisceaux 
$F_{|U \toinj T}$ sont des $\O_T$-modules quasi-cohérents.

On montre (en adaptant la preuve du lemme 3.3.2 de \cite{breuil-bull})
que si $\calF$ est un faisceau à composantes quasi-cohérentes sur 
$(X/\Upsilon)_\SYNCRIS$, alors $R^i v_\star \calF = 0$ pour tout $i \geq 
1$. Autrement dit $R v_\star \calF \simeq v_\star \calF$. De même, par 
un calcul de Cech, on prouve (voir \cite{breuil-duke}, appendice C.1) 
que 
$R w_\star \calF \simeq w_\star \calF$ (toujours en supposant que 
$\calF$ est à composantes quasi-cohérentes).

\bigskip

On dispose d'un résultat de comparaison entre cohomologie 
log-cristalline et cohomologie de de Rham (théorème 6.4 de \cite{kato}) 
qui donne avec les rappels précédents le théorème suivant :

\begin{theo}
\label{4:th:ralphaj}
Soit $X$ un log-schéma fin localement de type fini sur $\Upsilon$. On 
suppose que l'on a une $\Upsilon$-immersion fermée $X \toinj Y$ avec $Y$ 
log-lisse sur $\Upsilon$. Soit $D$ l'enveloppe aux puissances divisées 
de $X$ dans $Y$ (voir paragraphe 6.n de \cite{kato} pour une 
définition). Alors, pour tout entier $s$ :
\begin{itemize}
\item si $\Upsilon = T_1$ :
$$R\alpha_\star \J_1^{[s]} = 
\J_D^{[s]} \to \J_D^{[s-1]} \otimes_{\O_Y} \omega^1_{Y/T_1} \to 
\J_D^{[s-2]} \otimes_{\O_Y} \omega^2_{Y/T_1} \to \cdots$$
\item si $\Upsilon = \tilde T$ :
$$R\alpha_\star \tilde \J^{[s]} = 
\J_D^{[s]} \to \J_D^{[s-1]} \otimes_{\O_Y} \omega^1_{Y/\tilde T} 
\to \J_D^{[s-2]} \otimes_{\O_Y} \omega^2_{Y/\tilde T} \to \cdots$$
\item si $\Upsilon = \bar T$ :
$$R\alpha_\star \bar \J^{[s]} = 
\J_D^{[s]} \to \J_D^{[s-1]} \otimes_{\O_Y} \omega^1_{Y/\bar T} 
\to \J_D^{[s-2]} \otimes_{\O_Y} \omega^2_{Y/\bar T} \to \cdots$$
\end{itemize}
et :
\begin{itemize}
\item si $\Upsilon = T_1$ :
$$R\alpha_\star \pa{\J_1^{[s]} / \J_1^{[s+1]}} = 
\frac{\J_D^{[s]}}{\J_D^{[s+1]}} \to 
\frac{\J_D^{[s-1]}}{\J_D^{[s]}} \otimes_{\O_Y} \omega^1_{Y/T_1} \to 
\frac{\J_D^{[s-2]}}{\J_D^{[s-1]}} \otimes_{\O_Y} \omega^2_{Y/T_1} \to 
\cdots$$
\item si $\Upsilon = \tilde T$ :
$$R\alpha_\star \pa{\tilde \J^{[s]} / \tilde \J^{[s+1]}} = 
\frac{\J_D^{[s]}}{\J_D^{[s+1]}} \to 
\frac{\J_D^{[s-1]}}{\J_D^{[s]}} \otimes_{\O_Y} \omega^1_{Y/\tilde T} 
\to \frac{\J_D^{[s-2]}}{\J_D^{[s-1]}} \otimes_{\O_Y} 
\omega^2_{Y/\tilde T} \to \cdots$$
\item si $\Upsilon = \bar T$ :
$$R\alpha_\star \pa{\bar \J^{[s]} / \bar \J^{[s+1]}} = 
\frac{\J_D^{[s]}}{\J_D^{[s+1]}} \to 
\frac{\J_D^{[s-1]}}{\J_D^{[s]}} \otimes_{\O_Y} \omega^1_{Y/\bar T} 
\to \frac{\J_D^{[s-2]}}{\J_D^{[s-1]}} \otimes_{\O_Y} 
\omega^2_{Y/\bar T} \to \cdots$$
\end{itemize}
\end{theo}

\noindent
{\it Remarque.} En gardant les notations du théorème, si on suppose en 
plus $X$ log-lisse sur $\Upsilon$, on peut choisir $Y = X$ et alors 
$\J_D^{[s]} = \O_X$ si $s \leq 0$ et $\J_D^{[s]} = 0$ sinon. Les 
expressions des $R\alpha_\star$ se simplifient alors considérablement. 
Par exemple, dans le cas où $\Upsilon = \bar T$, on obtient :
$$R\alpha_\star \bar \J^{[s]} =
0 \to \cdots \to 0 \to \omega^s_{X/\bar T} \to \omega^{s+1}_{X/\bar T}
\to \omega^{s+2}_{X/\bar T} \to \cdots $$
$$R\alpha_\star \pa{\bar \J^{[s]} / \bar \J^{[s+1]}} =
0 \to \cdots \to 0 \to \omega^s_{X/\bar T} \to 0 \to \cdots .$$
Retenons ces résolutions que l'on sera amené à réutiliser par la suite.

\bigskip

On peut montrer également (voir appendice B de \cite{breuil-duke}) un 
raffinement du théorème précédent qui permet de tenir compte des 
faisceaux $\O_1^\car$ :

\begin{theo}
\label{4:th:ralphacar}
Soit $X$ un log-schéma fin localement de type fini sur $T_1$. On 
suppose que l'on a une $T_1$-immersion fermée $X \toinj Y$ avec $Y$ 
log-lisse sur $T_1$. Soit $D$ l'enveloppe aux puissances divisées 
de $X$ dans $Y$. Alors, pour tout entier $s$ :
$${\small \begin{array}{l}
R\alpha_\star \pa{\O_1^\car \otimes^\SYN_{\O_1}\J_1^{[s]} / 
\J_1^{[s+1]}} = \\
\hspace{1cm}
\O_1^\car \otimes_{\O_X} \frac{\J_D^{[s]}}{\J_D^{[s+1]}} \to 
\pa{\O_1^\car \otimes_{\O_X} \frac{\J_D^{[s-1]}}{\J_D^{[s]}}}
\otimes_{\O_Y} \omega^1_{Y/T_1} \to \pa{\O_1^\car \otimes_{\O_X} 
\frac{\J_D^{[s-2]}}{\J_D^{[s-1]}}} \otimes_{\O_Y} \omega^2_{Y/T_1} 
\to \cdots
\end{array}}$$
où \og $\otimes^\SYN$ \fg\ signifie que l'on prend le faisceau associé 
au faisceau produit tensoriel pour la topologie syntomique sur $X$.

On a d'autre part des versions analogues avec les bases $\tilde T$ et 
$\bar T$.
\end{theo}

\noindent
{\it Remarque.}
Si $X$ est de plus du type de Cartier, le faisceau $\O_1^\car$ sur 
$X_\et$ s'identifie à $\Sigma \otimes_{(\phi),k} \O_X$ où $\Sigma$ 
est l'anneau de $\Upsilon$. Dans ce dernier cas, on obtient une 
résolution plus simple.

\bigskip

Intéressons-nous maintenant à un analogue du théorème \ref{4:th:ralphaj} 
pour les faisceaux $\J_1^{[q]}$, $\tilde \J^{[q]}$ et $\bar \J^{[q]}$ où 
$q$ n'est pas forcément un entier. Soit $q \in \frac 1 e \N$. Si 
$X/\Upsilon$ est comme précédemment un log-schéma fin, localement de 
type fini, on commence par définir des faisceaux $\J_{X/\Upsilon}^{[q]}$ 
sur les sites $(X/\Upsilon)_\CRIS$ et $(X/\Upsilon)_\SYNCRIS$ en posant :
$$\J_{X/\Upsilon}^{[q]} = u^\delta \J_{X/\Upsilon}^{[s]} + 
\J_{X/\Upsilon}^{[s+1]}$$
si $eq = es + \delta$ est la division euclidienne de $eq$ par $e$. (On 
remarque que lorsque $\Upsilon = \bar E$, on a $u = 0$ et donc le 
premier terme de la somme n'intervient que si $\delta = 0$).

On s'intéresse particulièrement au quotient $\tilde \J^{[s]} / 
\tilde \J^{[q]}$. Par définition, on a immédiatement :
$$\tilde \J^{[s]} / \tilde \J^{[q]} = \frac {\tilde \J^{[s]} / \tilde 
\J^{[s+1]}} {u^\delta (\tilde \J^{[s]} / \tilde \J^{[s+1]})}$$
et une égalité analogue pour les faisceaux $\J_{X/\Upsilon}^{[\cdot]}$.
Du fait que $R v_\star \calF = v_\star \calF$ (resp. $R w_\star 
\calF = w_\star \calF$) si $\calF$ est un faisceau à composantes 
quasi-cohérentes, on déduit :
$$R v_\star \pa{\J_{X/\tilde T}^{[s]} / \J_{X/\tilde T}^{[q]}} = 
\J_{X/\tilde T}^{[s]} / \J_{X/\tilde T}^{[q]}
\quad (\text{resp. }R w_\star \pa{\J_{X/\tilde T}^{[s]} / \J_{X/\tilde 
T}^{[q]}} = \tilde \J^{[s]} / \tilde \J^{[q]}).$$

En vérifiant soigneusement que la démonstration de l'appendice B de 
\cite{breuil-duke} s'applique encore dans ce contexte, on obtient le 
théorème suivant :

\begin{theo}
\label{4:th:ralphajq}
Avec les notations précédentes, on a :
$$\begin{array}{l}
R\alpha_\star \pa{\tilde \J^{[s]} / \tilde \J^{[q]}} = \\
\hspace{0.5cm} \displaystyle
\frac{\J_D^{[s]} / \J_D^{[s+1]}}{u^\delta(\J_D^{[s]} / \J_D^{[s+1]})} 
\to \frac{\J_D^{[s-1]} / \J_D^{[s]}} {u^\delta(\J_D^{[s-1]} / 
\J_D^{[s]})} \otimes_{\O_Y} \omega^1_{Y/\tilde T} \to 
\frac{\J_D^{[s-2]} / \J_D^{[s-1]}} { u^\delta(\J_D^{[s-2]} / 
\J_D^{[s-1]})} \otimes_{\O_Y} \omega^2_{Y/\tilde T} \to \cdots
\end{array}$$
et :
$$\begin{array}{l}
R\alpha_\star \pa{\tilde \O^\car \otimes_{\O_1} \tilde \J^{[s]} / 
\tilde \J^{[q]}} = \\
\hspace{0.5cm} \displaystyle
\tilde \O^\car \otimes_{\O_X} \frac{\J_D^{[s]} / 
\J_D^{[s+1]}}{u^\delta(\J_D^{[s]} / \J_D^{[s+1]})} 
\to \pa{ \tilde \O^\car \otimes_{\O_X} \frac{\J_D^{[s-1]} / \J_D^{[s]}} 
{u^\delta(\J_D^{[s-1]} / \J_D^{[s]})}} \otimes_{\O_Y} \omega^1_{Y/\tilde 
T} 
\to \cdots 
\end{array}.$$
\end{theo}

\noindent
{\it Remarque.} Bien évidemment, on a des versions analogues pour les
bases $T_1$ et $\bar T$. En outre, il est possible d'écrire des 
résolutions de ce type pour beaucoup d'autres faisceaux. Cependant,
pour cet article, nous aurons seulement besoin de celle-ci.

\paragraph{Des isomorphismes à la Deligne-Illusie}

On reprend les notations et les hypothèses du début de 
paragraphe \ref{4:sec:crisp} : $X_1$ désigne un log-schéma propre, 
log-lisse et du type de Cartier sur $T_1$ et on suppose qu'il admet un 
relèvement fin et log-lisse $X_2$ sur $T_2$.

\begin{theo}
\label{4:th:isomcoh}
On garde les notations que l'on vient de rappeler. Alors, pour tout $0 
\leq i \leq s < \frac p e$, on a :
$$\xymatrix @C=40pt @R=8pt {
\xytilde S \otimes_{(\phi),k} H^i((X_1)_\syn, \xytilde \J^{[s]} / 
\xytilde \J^{[s+1/e]}) \ar[r]^-{\sim}_-{\id \otimes \phi_s} & 
H^i((X_1)_\syn, \xytilde \O^\st_1) \\
\hspace{3.7em}
H^i((X_1)_\syn, \xytilde \J^{[s]} / \xytilde \J^{[s+1/e]})
\ar[r]^-{\sim}_-{\phi_s} & H^i((X_1)_\syn, \xybar \O^\st_1) } $$
et des versions équivalentes en remplaçant \og $\tilde \J$ \fg\ par 
\og $\J_1$ \fg.
\end{theo}

\begin{preuve}
Comme les deux faisceaux $\tilde \O^\car$ et $\tilde \J^{[s]} / \tilde 
\J^{[s+1/e]}$ sont tués par $\pi$, le produit tensoriel $\tilde \O^\car 
\otimes_{\O_1} \tilde \J^{[s]} / \tilde \J^{[s+1/e]}$ est isomorphe à
$\tilde \O^\car \otimes_{\bar O} \tilde \J^{[s]} / \tilde \J^{[s+1/e]}$.

Comme $X_1$ est supposé du type de Cartier, le faisceau $\tilde \O^\car$ 
sur le site $(X_1)_\et$ s'identifie à $\tilde S \otimes_{(\phi), k} \bar 
\O$. Dans ces conditions, le théorème \ref{4:th:ralphajq} assure que la
flèche canonique :
\begin{equation}
\label{4:eq:ralpha1}
\tilde S \otimes_{(\phi), k} R\alpha_\star (\tilde \J^{[s]} / \tilde 
\J^{[s+1/e]}) \to R\alpha_\star (\tilde \O^\car \otimes_{\bar O} \tilde 
\J^{[s]} / \tilde \J^{[s+1/e]}) = R\alpha_\star (\tilde \O^\car
\otimes_{\O_1} \tilde \J^{[s]} / \tilde \J^{[s+1/e]})
\end{equation}
est un isomorphisme ($\alpha$ désigne le morphisme de topoï 
$\widetilde{(X_1)_\SYN} \to \widetilde{(X_1)_\et}$). Par le
théorème \ref{4:th:isomfais}, on a un nouvel isomorphisme :
\begin{equation}
\label{4:eq:ralpha2}
R\alpha_\star (\tilde \O^\car \otimes_{\O_1} \tilde \J^{[s]} / \tilde 
\J^{[s+1/e]}) \simeq R\alpha_\star (F_s^\car \tilde \O^\st).
\end{equation}

D'autre part, en combinant la proposition \ref{4:prop:fscar} 
et le théorème \ref{4:th:ralphacar}, on obtient pour tout entier $t$ :
$$R\alpha_\star (F_t^\car \tilde \O^\st / F_{t-1}^\car \tilde \O^\st) 
\simeq \tilde S \otimes_{(\phi), k} R\alpha_\star (\bar \J^{[t]} / \bar 
\J^{[t+1]}) = \tilde S \otimes_{(\phi),k} \omega^t_{\bar X / \bar T} 
[-t]$$
où $\bar X = X \times_{T_1} \bar T$. On en déduit que :
$$\tau_{\leq t-1} R\alpha_\star (F_t^\car \tilde \O^\st) \simeq 
\tau_{\leq t-1} R\alpha_\star (F_{t-1}^\car \tilde \O^\st)$$
et puis par une récurrence immédiate :
$$\tau_{\leq s} R\alpha_\star (F_s^\car \tilde \O^\st) \simeq 
\tau_{\leq s} R\alpha_\star (F_t^\car \tilde \O^\st)$$
pour tout entier $t \geq s$. Le théorème découle alors de la 
propriété $\bigcup_{t \in \N} F_t^\car \tilde \O^\st = \tilde 
\O^\st$ (proposition \ref{4:prop:fscar}) et des isomorphismes 
(\ref{4:eq:ralpha1}) et (\ref{4:eq:ralpha2}).

On traite de manière exactement similaire le cas de \og $\J_1$ \fg.
\end{preuve}

\bigskip

\noindent
{\it Remarque.} En utilisant l'isomorphisme donné par la remarque qui 
suit le théorème \ref{4:th:isomfais}, on obtient des isomorphismes 
analogues qui s'écrivent :
$$\xymatrix @C=40pt @R=8pt {
\xytilde S \otimes_{(\phi),\O_K / p} H^i((X_1)_\syn, \xytilde \J^{[s]} / 
\xytilde \J^{[s+1]}) \ar[r]^-{\sim}_-{\id \otimes \phi_s} & 
H^i((X_1)_\syn, \xytilde \O^\st_1) \\
k \otimes_{(\phi),\O_K / p} H^i((X_1)_\syn, \xytilde \J^{[s]} / \xytilde 
\J^{[s+1]}) \ar[r]^-{\sim}_-{\id \otimes \phi_s} & H^i((X_1)_\syn, \xybar 
\O^\st_1) } $$
où $\Sigma = k \text{ ou } \tilde S$ est vu comme une $\O_K / 
p$-algèbre par la composée $\O_K / p \to k \to \Sigma$, la première
flèche étant la projection canonique et la seconde le Frobenius.
En particulier, mis ensemble ces isomorphismes impliquent que le
morphisme naturel :
$$\frac{H^i((X_1)_\syn, \tilde \J^{[s]} / \tilde \J^{[s+1]})} {u \cdot 
H^i((X_1)_\syn, \tilde \J^{[s]} / \tilde \J^{[s+1]})} \longrightarrow 
H^i((X_1)_\syn, \tilde \J^{[s]} / \tilde \J^{[s+1/e]})$$
est un isomorphisme pour tout $i \leq s < \frac p e$.

\subsubsection{Fin de la preuve}
\label{4:subsec:cristilde}
On garde les notations et les hypothèses introduites au début du 
paragraphe \ref{4:sec:crisp}.

\medskip

Au vu du théorème \ref{4:th:isomcoh}, il reste à prouver, pour s'assurer 
que le quadruplet :
$$(H^i ((X_1)_\syn, \tilde \O^\st), H^i ((X_1)_\syn, \tilde \J^{[r]}), 
\phi_r, N)$$
est un objet de $\Mrtilde$, les deux choses suivantes :
\begin{enumerate}
\item la flèche canonique $H^i ((X_1)_\syn, \tilde \J^{[r]}) \to H^i 
((X_1)_\syn, \tilde \J^{[r]}/\tilde \J^{[r+1/e]})$ est surjective,
\item la flèche canonique $H^i ((X_1)_\syn, \tilde \J^{[r]}) \to H^i 
((X_1)_\syn, \tilde \O^\st)$ est injective,
\end{enumerate}
les autres propriétés de compatibilité étant claires.

Si $\calF$ est un faisceau sur $(X_1)_\syn$, nous notons 
simplement $H^i(\calF)$ pour $H^i ((X_1)_\syn, \calF)$ et si ce dernier 
est un espace vectoriel de dimension finie sur $k$, nous notons 
$h^i(\calF)$ sa dimension.

\paragraph{Propriétés de finitude.}
Pour la suite, nous aurons besoin de raisonner sur les dimensions de 
certains groupes de cohomologie. Il nous faut donc prouver dans un 
premier temps qu'ils sont de dimension finie. Nous commençons par 
donner un résultat agréable sur les faisceaux $u^k \tilde \O^\st$ :

\begin{prop}
\label{4:prop:ukost}
Pour tous entiers $k$ et $i < \frac p e$ (on impose $i=0$ si $p=2$ et 
$e=1$), les $k$-espaces vectoriels $H^i(u^k \tilde \O^\st)$ sont de 
dimension finie et égaux à $u^k H^i(\tilde \O^\st)$.
\end{prop}

\begin{preuve}
En premier lieu, on remarque que si $k \geq p$, tout est nul et donc 
la proposition est trivialement vérifiée.
Montrons que $H^i(\tilde \O^\st)$ est de dimension finie sur 
$k$. Le log-schéma $\bar X = X_1 \otimes_{T_1} \bar T$ est log-lisse sur 
$\bar T = \bar E$ par changement de base et donc, par la remarque qui 
suit le théorème \ref{4:th:ralphaj}, $R \alpha_\star \bar \O^\st = 
\O_{\bar X / \bar E} \to \omega^1_{\bar X / \bar E} \to \cdots$. Comme 
$\bar X$ est propre, les $\omega^j_{\bar X / \bar E}$ sont de dimension 
finie, et il en est donc de même de $H^i(\bar \O^\st)$. Le théorème 
\ref{4:th:isomcoh} permet alors de conclure.

\medskip

Prouvons désormais la proposition par récurrence sur $i$. En écrivant 
la suite exacte longue associée à :
$$\xymatrix @C=25pt {
0 \ar[r] & u^k \xytilde \O^\st \ar[r] & \xytilde \O^\st 
\ar[r]^-{u^{p-k}} & u^{p-k} \xytilde \O^\st \ar[r] & 0 } $$
et en appliquant l'hypothèse de récurrence, on prouve que la flèche 
$H^i(u^k \O^\st) \to H^i(\O^\st)$ est injective. Ainsi $H^i(u^k 
\tilde \O^\st)$ est de dimension finie car inclus dans $H^i(\tilde 
\O^\st)$. Par ailleurs, on vérifie facilement que l'on dispose 
d'une suite exacte de faisceaux, pour tout entier $k \leq p-1$ :
$$\xymatrix @C=25pt {
0 \ar[r] & u^{p-1} \xytilde \O^\st \ar[r] & u^k \xytilde 
\O^\st \ar[r]^u & u^{k+1} \xytilde \O^\st \ar[r] & 0 }.$$
Par l'hypothèse de récurrence, la suite exacte longue associée prend la 
forme :
$$\xymatrix @C=25pt {
0 \ar[r] & H^i(u^{p-1} \xytilde \O^\st) \ar[r] & H^i(u^k \xytilde 
\O^\st) \ar[r] & H^i(u^{k+1} \xytilde \O^\st) }$$
et fournit une inégalité sur les dimensions à savoir $h^i(u^k \tilde 
\O^\st) \leq h^i(u^{k+1} \tilde \O^\st) + h^i(u^{p-1} \tilde \O^\st)$. 
Or $u^{p-1} \tilde \O^\st \simeq \bar \O^\st$, d'où en additionnant les 
inégalités précédentes pour $k$ variant de $0$ à $p$, on obtient 
$h^i(\tilde \O^\st) \leq p h^i(\bar \O^\st)$. Or le théorème 
\ref{4:th:isomcoh} prouve qu'il y a en fait égalité entre les deux nombres 
précédents. Cela implique que toutes les inégalités sommées sont des 
égalités et par suite que l'on a des suites exactes courtes :
$$\xymatrix @C=25pt {
0 \ar[r] & H^i(u^{p-1} \xytilde \O^\st) \ar[r] & H^i(u^k \xytilde
\O^\st) \ar[r] & H^i(u^{k+1} \xytilde \O^\st \ar[r]) & 0 }.$$
Si $k \leq p$, le morphisme de multiplication par $u^k$ se factorise par 
$H^i(\O^\st) \to H^i(u^k \O^\st) \to H^i(\O^\st)$. D'après ce qui 
précède, la première flèche est surjective et la seconde est injective. 
On en déduit le résultat annoncé.
\end{preuve}

\bigskip

Intéressons-nous à présent aux faisceaux $u^k \tilde \J^{[q]}$ et 
commençons par un lemme qui les relie entre eux :

\begin{lemme}
\label{4:lem:exactj}
Soient $k \in \{0, \ldots, p-1\}$ et $q \in \frac 1 e \N$. Notons
$s$ la partie entière de $q$. Alors on a une suite exacte :
$$\xymatrix @C=30pt {
0 \ar[r] & u^{k+1} \xytilde \J^{[q]} \ar[r] & u^k \xytilde \J^{[q+1/e]} 
\ar[r] & \xybar \J^{[s+1]} \ar[r] & 0 }.$$
\end{lemme}

\begin{preuve}
Commençons par préciser les flèches qui apparaissent dans la suite 
exacte. La première $u^{k+1} \tilde \J^{[q]} \to u^k \tilde 
\J^{[q+1/e]}$ est simplement l'inclusion naturelle entre deux 
sous-faisceaux de $\tilde \O^\st$. La seconde flèche est légèrement plus 
subtile. Remarquons en premier lieu, que si l'on note $K$ le noyau de la 
projection (multiplication par $u^k$) $\tilde \J^{[s+1]} \to u^k \tilde 
\J^{[s+1]}$, on dispose d'un diagramme comme suit :
$$\xymatrix @C=25pt {
0 \ar[r] & K \ar[rd] \ar[r] & \xytilde \J^{[s]} \ar[d] \ar[r]^{u^k} & 
u^k \xytilde \J^{[s]} \ar[r] & 0 \\
& & \xybar \J^{[s]} }$$
La flèche diagonale composée est nulle : en effet, d'après les 
descriptions locales, tout élément de $K$ est un multiple de $u$ (au 
moins dans $\tilde \O^\st$) et donc s'annule lorsqu'on le projette dans 
$\bar \J^{[s]}$. On en 
déduit un morphisme de faisceaux $u^k \tilde \J^{[s]} \to \bar 
\J^{[s]}$ (qui correspond moralement à la division par $u^k$).
Par ailleurs, on dispose d'une inclusion $\tilde \J^{[q+1/e]} \to \tilde
\J^{[s]}$ et la seconde flèche de la suite exacte est la composée 
$u^k \tilde \J^{[q+1/e]} \to u^k \tilde \J^{[s]} \to \bar \J^{[s]}$.

Il reste à vérifier que la flèche précédente tombe en fait dans
$\bar \J^{[s+1]}$ et que la suite obtenue ainsi est bien exacte.
Montrons dans un premier temps un résultat analogue sur les faisceaux 
\og $\J_1$ \fg, à savoir que la suite :
$$\xymatrix {
0 \ar[r] & u \J_1^{[q]} \ar[r] & \J_1^{[q+1/e]} \ar[r] & \xybar \J^{[s]}
}$$
est exacte et que l'image de la dernière flèche est $\bar \J^{[s+1]}$.
Notons $f : u \J_1^{[q]} \to \J_1^{[q+1/e]}$ et $g : \J_1^{[q+1/e]} 
\to \bar \J^{[s]}$ les applications qui interviennent. Il est clair que 
$f$ est injective et que $g \circ f = 0$. Notons $q = s + \frac \delta 
e$ avec $0 \leq \delta < e$. En reprenant les notations du paragraphe 
\ref{4:sec:osttilde}, un élément de
$\tilde \J^{[q+1/e]} (A^\infty, P^\infty)$ s'écrit comme une somme de
multiples de termes d'une des deux formes suivantes :
\begin{itemize}
\item[i)] $u^{e m_0 + \delta + 1} \psi_1^{m_1} \cdots \psi_t^{m_t}
X^{m_{t+1}}$ avec $m_0 + \cdots + m_{t+1} \geq s$ ;
\item[ii)] $u^{e m_0} \psi_1^{m_1} \cdots \psi_t^{m_t} X^{m_{t+1}}$ avec
$m_0 + \cdots + m_{t+1} \geq s+1$.
\end{itemize}
On vérifie directement que les éléments du premier type s'envoient
sur $0$ par $g$ et que les éléments du second type s'envoient dans
$\bar \J^{[s+1]}$. Finalement l'application $g$ tombe bien dans $\bar
\J_{[s+1]}$ comme annoncé. En outre, un élément de $\tilde \J^{[s+1]}$
s'écrit comme une somme de multiples de $u_{e m_0} \psi_1^{m_1} \cdots
\psi_t^{m_t} X^{m_{t+1}}$ avec $m_0 + \cdots + m_{t+1} \geq s+1$, et
admet donc un antécédent par $g$.

Il ne reste qu'à vérifier l'exactitude au milieu. Soit $x \in  
\J_1^{[q+1/e]} (A^\infty, P^\infty)$ tel que $g(x) = 0$.
La flèche $u \J_1^{[q]} / \tilde \J^{[s+1]} \to \J_1^{[q+1/e]} / \tilde 
\J^{[s+1]}$ résultant de l'inclusion canonique est un isomorphisme et 
donc, quitte à modifier $x$, on peut supposer qu'il est élément de 
$\J_1^{[s+1]} (A^\infty, P^\infty)$. Notons $\bar x$ l'image de $x$ 
dans :
$$\frac{\tilde \J^{[s+1]} (A^\infty, P^\infty)}{\tilde \J^{[s+2]}
(A^\infty, P^\infty)}
= \bigoplus_{\sum m_i = s+1} (\calB/(u^e, \psi_1, \ldots,
\psi_t, X)) \cdot u^{e m_0} \gamma_{m_1} (\psi_1) \cdots 
\gamma_{m_t} (\psi_t) \gamma_{m_{t+1}} (X) .$$
Le fait que $g(x) = 0$ implique que dans la somme directe précédente,
$\bar x$ n'a des composantes non nulles que sur les $(\calB/(u^e,
\psi_1, \ldots, \psi_t, X)) \cdot u^{e m_0} \psi_1^{m_1} \cdots
\psi_t^{m_t} X^{m_{t+1}}$ avec $m_0 > 0$. Il s'ensuit $x \in u 
\J_1^{[s+1]} + \J_1^{[s+2]} \subset u \J_1^{[q]} + \J_1^{[s+2]}$. Quitte 
à faire une nouvelle modification, on peut donc supposer $x \in 
\J_1^{[s+2]} (A^\infty, P^\infty)$, et par une récurrence immédiate $x 
\in \J_1^{[Np]} (A^\infty, P^\infty)$ pour un certain entier $N$ qu'il 
reste à choisir.

Par ailleurs, on dispose d'une décomposition de $\O_1^\st (A^\infty, 
P^\infty)$ (formule \ref{4:eq:otilde})) :
$$\O_1^\st (A^\infty, P^\infty) = \bigoplus_{m_1, \ldots, m_{t+1}
\in \N} \tilde \calB \cdot \gamma_{p m_1} (\psi_1) \cdots \gamma_{p m_t} 
(\psi_t) \gamma_{p m_{t+1}} (X)$$
et d'une décomposition analogue de $\bar \O^\st (A^\infty, P^\infty)$ :
$$\bar \O^\st (A^\infty, P^\infty) = \bigoplus_{m_1, \ldots, m_{t+1}  
\in \N} \bar \calB \cdot \gamma_{p m_1} (\psi_1) \cdots \gamma_{p m_t}
(\psi_t) \gamma_{p m_{t+1}} (X) .$$
On vérifie directement que $g$ respecte ces décompositions. On
peut donc supposer que $x$ est élément de l'un des termes correspondant
à un uplet $(m_0, \ldots, m_{t+1})$ de la première somme directe. De 
plus comme $x \in \tilde \J^{[Np]} (A^\infty, P^\infty)$, on vérifie 
que si $x \neq 0$, alors $m = m_0 + \cdots + m_{t+1} \geq N-t-2$. Mais 
alors l'hypothèse $g(x) = 0$ implique $x \in u \J_1^{[mp]} (A^\infty, 
P^\infty) \subset u \J_1^{[q]} (A^\infty, P^\infty)$ si $N$ est choisi
suffisamment grand.

\medskip

Pour en déduire le lemme, on raisonne à partir du diagramme commutatif 
suivant :
$$\xymatrix @C=40pt @R=10pt {
u \J_1^{[q]} \ar[r]^-{f} \ar[dd]_-{u^k} &
 \J_1^{[q+1/e]} \ar[rd]^-{g} \ar[dd]_-{u^k} \\
& & \xybar \J^{[s]} \\
u^{k+1} \xytilde \J^{[q]} \ar[r]^-{f_k} &
u^k \xytilde \J^{[q+1/e]} \ar[ur]^-{g_k} }$$
où les multiplications par $u^k$ sont surjectives. On en déduit dans
un premier temps $\im g_k = \im g = \bar \J^{[s+1]}$. Par ailleurs
l'injectivité de la flèche $f_k$ et le fait que $g_k \circ f_k = 0$
sont immédiats. Une chasse au diagramme facile permet alors de 
conclure.
\end{preuve}

\bigskip

En corollaire, on en déduit enfin la proposition suivante :

\begin{prop}
Pour tout entier $k$, tout $q \in \frac 1 e \N$ et tout entier $i$, 
l'espace vectoriel $H^i (u^k \tilde \J^{[q]})$ est de dimension finie 
sur $k$.
\end{prop}

\begin{preuve}
Pour $k \geq p$, l'assertion est évidente puisque $u^k \tilde 
\J^{[q]} = 0$.

Par ailleurs, comme $\bar X = X_1 \times_{T_1} \bar T$ est log-lisse sur 
$\bar T = \bar E$, on a par la remarque qui suit le théorème  
\ref{4:th:ralphaj}, $R \alpha_\star \bar \J^{[s]} = 0 \to \cdots \to 0 \to 
\omega^s_{\bar X / \bar E} \to \omega^{s+1}_{\bar X / \bar E} \to 
\cdots$, et donc puisque $\bar X$ est propre sur $\bar E$, les groupes 
$H^i (\bar \J^{[s+1]})$ sont de dimension finie pour tous entiers $i$ et 
$s$.
Si $k \leq p-1$, le lemme \ref{4:lem:exactj} assure que $H^i(u^k \tilde 
\J^{[q+1/e]})$ est de dimension finie si et seulement si $H^i(u^{k+1} 
\tilde \J^{[q]})$ l'est. On se ramène ainsi à $k = p$ ou $q = 0$.
Le premier cas est traité précédemment et le second par la proposition 
\ref{4:prop:ukost}.
\end{preuve}

\paragraph{Surjectivité de $\mathbf{H^i ((X_1)_\syn, \tilde \J^{[r]}) 
\to H^i ((X_1)_\syn, \tilde \J^{[r]}/\tilde \J^{[r+1/e]})}$.}
On commence par prouver un lemme :

\begin{lemme}
\label{4:lem:noyuj}
Soient $i \in \N$ et $q \in \frac 1 e \N$. Supposons $0 \leq i \leq q < 
\frac p e$ (et $i=q=0$ si $p=2$ et $e=1$). On a des suites exactes 
courtes :
$$\xymatrix @C=25pt {
0 \ar[r] & H^i(\xybar \O^\st) \ar[r] & H^i(\xytilde \J^{[q]}) \ar[r] & 
H^i(u \xytilde \J^{[q]}) \ar[r] & 0 }.$$
\end{lemme}

\begin{preuve}
Montrons tout d'abord que la suite :
$$\xymatrix @C=25pt {
0 \ar[r] & \xybar \O^\st \ar[r] & \xytilde \J^{[q]} \ar[r]^u &
u \xytilde \J^{[q]} \ar[r] & 0 }$$
est exacte. La première flèche résulte de l'inclusion $\bar \O^\st = 
u^{p-1} \tilde \O^\st \subset \tilde \J^{[q]}$ et donc est injective.
La surjectivité est également claire. Par ailleurs, le noyau de la 
multiplication par $u$ sur $\tilde \O^\st$ est $u^{p-1} \tilde \O^\st$
et donc le noyau de $\tilde \J^{[q]} \to u \tilde \J^{[q]}$ s'identifie 
à $(u^{p-1} \tilde \O^\st \cap \tilde \J^{[q]}) = u^{p-1} \tilde 
\O^\st$.

En écrivant la suite exacte longue associée à la suite exacte courte 
précédente, on obtient déjà l'exactitude au milieu dans la suite de 
l'énoncé du lemme. Pour l'injectivité, on remarque que d'après la 
proposition \ref{4:prop:ukost}, on a $H^i(\bar \O^\st) = u^{p-1} 
H^i(\tilde \O^\st)$ et donc la flèche composée $H^i(\bar \O^\st) \to 
H^i(\tilde \J^{[q]}) \to H^i(\tilde \O^\st)$ est injective. Il en est 
donc de même de la flèche $H^i(\bar \O^\st) \to H^i(\tilde \J^{[q]})$. 
La surjectivité découle de l'injectivité car les flèches de bord sont 
nulles.
\end{preuve}

\bigskip

La surjectivité de $H^i(\tilde \J^{[r]}) \to H^i(\tilde \J^{[r]} / 
\tilde \J^{[r+1/e]})$ résulte directement de la proposition 
plus générale suivante :

\begin{prop} 
\label{4:prop:surj} 
Soient $i$ un entier et $q \in \frac 1 e \N$ vérifiant $0 \leq i \leq q 
< \frac p e$ (et $i=q=0$ si $p=2$ et $e=1$). On a des suites exactes 
courtes : 
$$\xymatrix @C=25pt { 0 \ar[r] & H^i(\xytilde \J^{[q+1/e]}) \ar[r] & 
H^i(\xytilde \J^{[q]}) \ar[r] & H^i(\xytilde \J^{[q]} / \xytilde \J^{[q+1/e]}) 
\ar[r] & 0 }.$$
En outre, si $s$ désigne la partie entière de $q$, on a également un 
début de suite exacte : 
$$\xymatrix @C=25pt { 0 \ar[r] & H^{s+1}(\xytilde \J^{[q+1/e]}) \ar[r] & 
H^{s+1}(\xytilde \J^{[q]}) \ar[r] & H^{s+1}(\xytilde \J^{[q]} / \xytilde 
\J^{[q+1/e]})}.$$
\end{prop}

\begin{preuve}
La preuve résulte d'un calcul de dimension. Précisément, on va prouver 
que pour tout $0 \leq i \leq q < \frac p e$, on a $h^i(\tilde \J^{[q]}) 
= h^i(\tilde \J^{[q+1/e]}) + h^i(\tilde \J^{[q]} / \tilde 
\J^{[q+1/e]})$. On déduira alors la proposition par une récurrence 
immédiate sur $i$.

D'après le lemme précédent, on a déjà $h^i(\tilde \J^{[q]}) = h^i(\bar 
\O^\st) + h^i(u \tilde \J^{[q]})$. La formule (\ref{4:eq:jtilde}) montre 
que 
la multiplication par $u^{e(q-s)}$ induit un isomorphisme entre les 
faisceaux $\tilde \J^{[s]} / \tilde \J^{[s+1/e]}$ et $\tilde \J^{[q]} / 
\tilde \J^{[q+1/e]}$. Le théorème \ref{4:th:isomcoh} 
implique alors $h^i(\bar \O^\st) = h^i(\tilde \J^{[q]} / \tilde 
\J^{[q+1/e]})$. D'autre part, on a la suite exacte (lemme 
\ref{4:lem:exactj}) :
$$\xymatrix @C=30pt {
0 \ar[r] & u \xytilde \J^{[q]} \ar[r] & \xytilde \J^{[q+1/e]} \ar[r] & 
\xybar \J^{[s+1]} \ar[r] & 0 }.$$
Or, puisque $R \alpha_\star \bar \J^{[s+1]} = 0 \to \cdots \to 0 \to 
\omega^{s+1}_{\bar X / \bar E} \to \omega^{s+2}_{\bar X / \bar E} \to 
\cdots$ (par la remarque qui suit le théorème \ref{4:th:ralphaj}), on a 
$H^j (\bar \J^{[s+1]}) = 0$ pour tout $j \leq s$. On en déduit 
$H^j (u \tilde \J^{[s]}) = H^j (\tilde \J^{[s+1/e]})$ pour tout $j \leq 
s$. En particulier $h^i (u \tilde \J^{[s]}) = h^i (\tilde \J^{[s+1/e]})$ 
ce qui conclut la démonstration.
\end{preuve}

\paragraph{Injectivité de $\mathbf{H^i ((X_1)_\syn, \tilde \J^{[r]}) \to 
H^i ((X_1)_\syn, \tilde \O^\st)}$.}

Commençons par énoncer le lemme suivant :

\begin{lemme}
\label{4:lem:inj}
Pour tout $i < \frac p e$, l'application $H^i(u^{p-ei} \tilde \J^{[i]}) 
\to H^i(\tilde \J^{[p/e]})$ est injective.
\end{lemme}

\begin{preuve}
Par un raisonnement analogue à celui utilisé pour la preuve du lemme 
\ref{4:lem:exactj}, on montre que l'on a une suite exacte :
$$\xymatrix @C=25pt {
0 \ar[r] & u^{p-ei} \xytilde \J^{[i]} \ar[r] & \xytilde \J^{[p/e]} 
\ar[r]^-{u^{ei}} & u^{ei} \xytilde \J^{[p/e]} \ar[r] & 0 }.$$
Il suffit donc de prouver que $H^{i-1}(u^{ei} \tilde \J^{[p/e]}) = 0$. 
Par le même argument que celui utilisé dans la preuve de la proposition 
\ref{4:prop:surj}, on montre que $H^j(u^{ei} \tilde \J^{[p/e]}) = 
H^j(\tilde \J^{[i+p/e]})$ pour tout $j \leq \frac p e$ et donc \emph{a 
fortiori} pour tout $j \leq i$. 
Posons $q = i + \frac p e$ et notons $s$ la partie entière de $q$. On a 
une suite exacte :
$$\xymatrix @C=25pt {
0 \ar[r] & \xytilde \J^{[q]} \ar[r] & \xytilde \J^{[s]} \ar[r] & 
\xytilde \J^{[s]} / \xytilde \J^{[q]} \ar[r] & 0 }.$$
Étale-localement, on peut relever $X_1$ et un log-schéma $Y_1$ log-lisse 
(auquel on étend les puissances divisées) sur $\tilde E$. Par les 
théorèmes \ref{4:th:ralphaj} et \ref{4:th:ralphajq}, la flèche $R 
\alpha_\star \tilde \J^{[s]} \to R \alpha_\star (\tilde \J^{[s]} / 
\tilde \J^{[q]})$ s'écrit explicitement :
$$\xymatrix {
\J_{Y_1}^{[s]} \ar[r] \ar[d] & 
\J_{Y_1}^{[s-1]} \otimes_{\O_{X_1}} \omega^1_{Y_1 / \xytilde E} \ar[r] 
\ar[d] & 
\J_{Y_1}^{[s-2]} \otimes_{\O_{X_1}} \omega^2_{Y_1 / \xytilde E} \ar[r] 
\ar[d] & \cdots \\
\frac{\J_{Y_1}^{[s]} / \J_{Y_1}^{[s+1]}} {u^\delta 
(\J_{Y_1}^{[s]} / \J_{Y_1}^{[s+1]})} \ar[r] & 
\frac{\J_{Y_1}^{[s-1]} / \J_{Y_1}^{[s]}} {u^\delta 
(\J_{Y_1}^{[s-1]} / \J_{Y_1}^{[s]})} \otimes_{\O_{X_1}} 
\omega^1_{Y_1 / \xytilde E} \ar[r] & 
\frac{\J_{Y_1}^{[s-2]} / \J_{Y_1}^{[s-1]}} {u^\delta 
(\J_{Y_1}^{[s-2]} / \J_{Y_1}^{[s-1]})} \otimes_{\O_{X_1}} 
\omega^2_{Y_1 / \xytilde E} \ar[r] & \cdots } $$
où $\delta$ est le reste de la division euclidienne de $p$ par $e$. 
Or, si $j < i$, on a $s - j \geq \frac p e$ et donc $\gamma_{s-j}(u) = 
0$. Il s'ensuit $\J_{Y_1}^{[s-j]} = 0$. On montre de même que 
$\J_{Y_1}^{[s-i]}$ est tué par $u^\delta$. Ainsi $H^j(\tilde 
\J^{[s]}) = H^j(\tilde \J^{[s]} / \tilde \J^{[q]}) = 0$ pour $j < i$ et 
la flèche $H^i (\tilde \J^{[s]}) \to H^i (\tilde \J^{[s]} / \tilde 
\J^{[q]})$ est un isomorphisme. La nullité de $H^{i-1}(u^{ei} \tilde 
\J^{[p/e]})$ (et donc le lemme) résulte alors d'une écriture de 
la suite exacte longue associée à la suite exacte courte :
$$\xymatrix @C=25pt {
0 \ar[r] & \xytilde \J^{[q]} \ar[r] & \xytilde \J^{[s]} \ar[r] & 
\xytilde \J^{[s]} / \xytilde \J^{[q]} \ar[r] & 0 }.$$
\end{preuve}

\begin{prop}
\label{4:prop:inj}
Pour tout $i$ et $q \in \frac 1 e \N$ tels que $q \leq i \leq r$, 
l'application $H^i (\tilde \J^{[i]}) \to H^i (\tilde \J^{[q]})$ est 
injective.
\end{prop}

\begin{preuve}
Si $p=2$, on a nécessairement $q = i = r = 0$ et le résultat est 
évident. Supposons donc $p \geq 3$.

On raisonne par récurrence descendante sur $q$. Le résultat est trivial 
pour $q = i$. Supposons-le vrai pour un certain $q$ et démontrons-le 
pour $q - \frac 1 e$.
On vérifie facilement que l'on a des suites exactes courtes de faisceaux :
$$\xymatrix @C=30pt @R=15pt {
0 \ar[r] & u^{ei} \xytilde \O^\st \ar@{=}[d] \ar[r] & \xytilde \J^{[i]} 
\ar[d] \ar[r]^-{u^{p-ei}} & u^{p-ei} \xytilde \J^{[i]} \ar[d] \ar[r] & 
0 \\
0 \ar[r] & u^{ei} \xytilde \O^\st \ar[r] & \xytilde \J^{[q-1/e]} 
\ar[r]^-{u^{p-ei}} & u^{p-ei} \xytilde \J^{[q-1/e]} \ar[r] & 0 } $$
qui donnent lieu à de nouvelles suites exactes :
$$\xymatrix @C=30pt @R=15pt {
& H^i(u^{ei} \xytilde \O^\st) \ar@{=}[d] \ar[r] & H^i(\xytilde 
\J^{[i]}) \ar[d] \ar[r] & H^i(u^{p-ei} \xytilde \J^{[i]}) \ar[d] \\
0 \ar[r] & H^i(u^{ei} \xytilde \O^\st) \ar[r] & H^i(\xytilde \J^{[q-1/e]}) 
\ar[r] & H^i(u^{p-ei} \xytilde \J^{[q-1/e]}) } $$
La deuxième suite est exacte à gauche car la flèche composée 
$H^i(u^{ei}  \tilde \O^\st) \to  H^i(\tilde \J^{[q-1/e]}) \to  
H^i(\tilde \O^\st)$ est injective (voir preuve de la proposition 
\ref{4:prop:ukost}). On veut montrer que la flèche verticale du milieu est 
injective, et une chasse au diagramme laissée au lecteur assure que pour 
cela, il suffit de prouver que la flèche verticale de droite l'est.

Or, on peut former le carré commutatif suivant :
$$\xymatrix @C=30pt @R=10pt {
H^i(u^{p-ei} \xytilde \J^{[i]}) \ar[dd] \ar[r] & H^i(\xytilde \J^{[p/e]}) 
\ar[d] \\
& H^i(\xytilde \J^{[i]}) \ar[d] \\
H^i(u^{p-ei} \xytilde \J^{[q-1/e]}) \ar[r] & H^i(\xytilde \J^{[q]}) }$$
La flèche du haut est injective d'après le lemme \ref{4:lem:inj} et celles 
de droite le sont également respectivement d'après la proposition 
\ref{4:prop:surj} et l'hypothèse de récurrence. On en déduit que celle de 
gauche l'est aussi comme on le souhaitait.
\end{preuve}

\bigskip

On a finalement le théorème :

\begin{theo}
\label{4:th:obtilde}
Pour tout $i \leq r$, le quadruplet $(H^i ((X_1)_\syn, \tilde \O^\st), 
H^i ((X_1)_\syn, \tilde \J^{[r]}), \phi_r, N)$ définit un 
objet de la catégorie $\Mrtilde$.
\end{theo}

\begin{preuve}
Il restait à prouver la surjectivité de  $H^i (\tilde \J^{[r]}) \to H^i 
(\tilde \J^{[r]}/\J^{[r+1/e]})$ et l'injectivité de $H^i (\tilde 
\J^{[r]}) \to H^i (\tilde \O^\st)$. Le premier point est une conséquence 
immédiate de la proposition \ref{4:prop:surj}.
Pour le second point, on remarque que le morphisme $H^i (\tilde
\J^{[r]}) \to H^i (\tilde \O^\st)$ se factorise par $H^i (\tilde
\J^{[r]}) \to H^i (\tilde \J^{[i]}) \to H^i (\tilde \O^\st)$. La 
première des deux flèches précédentes est injective par la proposition 
\ref{4:prop:surj} et la seconde est aussi injective par la proposition 
\ref{4:prop:inj}. Ceci clôt la démonstration.
\end{preuve}

\subsubsection{Reformulation sur la base $E_1$}
\label{4:subsec:crise1}
On montre dans ce paragraphe que le quadruplet :
$$(H^i ((X_1)_\syn, \O^\st_1), H^i ((X_1)_\syn, \J^{[r]}_1), \phi_r, 
N)$$
est un objet de la catégorie $\Mr$ pour tout $i \leq r$.
Pour cela, on commence par recopier à la lettre les arguments des 
paragraphes \ref{4:subsec:isomcoh} et \ref{4:subsec:cristilde} pour obtenir 
un équivalent du théorème \ref{4:th:isomcoh} qui s'énonce comme suit :

\begin{theo}
\label{4:th:isomcohs1}
Pour $0 \leq i \leq s \leq p-1$ (et seulement pour $i=s=0$ si $p=2$), on 
a un isomorphisme :
$$\xymatrix @C=30pt {
S_1 \otimes_{(\phi), k[u]/u^e} H^i((X_1)_\syn, 
\J_1^{[s]}/\J_1^{[s+1]}) \ar[r]^-{\sim}_-{\id \otimes \phi_s} &
H^i((X_1)_\syn, \O^\st_1) }.$$
\end{theo}

Dans un premier temps, il nous faut montrer que $H^i((X_1)_\syn, 
\O^\st_1)$ est un $S_1$-module libre et pour cela il suffit de prouver  
que $H^i((X_1)_\syn, \J_1^{[r]}/\J_1^{[r+1]})$ est un 
$k[u]/u^e$-module libre par le théorème \ref{4:th:isomcohs1}.
C'est évident si $e = 1$. \emph{À partir de maintenant et jusqu'à la 
fin de cette partie, on suppose $e \geq 2$}.

Dans ce cas, on est tenté de comparer les deux faisceaux 
$\J_1^{[r]}/\J_1^{[r+1]}$ et $\tilde \J^{[r]} / \tilde \J^{[r+1]}$ 
puisque la version \og ~$\tilde{}$~ \fg\ a déjà été étudiée. Cependant, 
en regardant les descriptions explicites, on se rend compte qu'il n'est 
pas vrai en général que ces deux faisceaux sont isomorphes ; c'est le 
cas simplement lorsque $er \leq p - e$.

\medskip

La solution consiste à introduire une nouvelle catégorie d'objets 
modulo $u^{2p}$ et à procéder en deux étapes : on passe des objets
modulo $u^p$ aux objets modulo $u^{2p}$ puis de ces derniers aux objets 
de $\Mr$.

\paragraph{Les objets modulo $u^{2p}$}

On rappelle que tout au long de ce paragraphe, on suppose $e \geq 2$.

\medskip

On introduit une nouvelle catégorie, notée $\Mrtilde_{(2)}$ dont la 
définition est très proche des autres catégories déjà introduites. 
On pose $\tilde S_{(2)} = k\cro u / u^{2p}$. C'est un anneau muni d'un 
Frobenius $\phi$ semi-linéaire envoyant $u^i$ sur $u^{pi}$ et d'un 
opérateur de monodromie $N$ $k$-linéaire envoyant $u^i$ sur $-i u^i$. 
Un objet de $\Mrtilde_{(2)}$ est la donnée de :
\begin{enumerate}
\item un $\tilde S_{(2)}$-module $\tilde \calM$ libre de rang fini ;
\item un sous-module $\Fil^r \tilde \calM$ de $\tilde \calM$ contenant 
$u^{er} \tilde \calM$ ;
\item une flèche $\phi$-semi-linéaire $\phi_r : \Fil^r \tilde \calM \to 
\tilde \calM$ telle que l'image de $\phi_r$ engendre $\tilde \calM$ en 
tant que $\tilde S_{(2)}$-module ;
\item une application $k$-linéaire $N : \tilde \calM \to \tilde \calM$ 
telle que :
\begin{itemize}
\item pour tout $\lambda \in \tilde S_{(2)}$ et tout $x \in \tilde 
\calM$, $N\pa{\lambda x} = N\pa \lambda x + \lambda N\pa x$
\item $u^e N(\Fil^r \tilde \calM) \subset \Fil^r \tilde \calM$
\item le diagramme suivant commute :
\[\xymatrix @C=50pt {
\Fil^r \xytilde \calM \ar[r]^{\phi_r} \ar[d]_{u^e N} & \xytilde \calM 
\ar[d]^{c_{(2),\pi} N} \\
\Fil^r \xytilde \calM \ar[r]^{\phi_r} & \xytilde \calM}\]
où $c_{(2),\pi}$ est la réduction de $c$ dans $\tilde S_{(2)}$.
\end{itemize}
\end{enumerate}
Les morphismes de $\Mrtilde_{(2)}$ sont les applications $\tilde 
S_{(2)}$-linéaires qui commutent à toutes les structures.

\medskip

On prouve comme dans le paragraphe \ref{4:sec:mrtilde} que les catégories 
$\Mr$, $\Mrtilde_{(2)}$ et $\Mrtilde$ sont toutes les trois équivalentes 
\emph{via} les foncteurs de réduction modulo $u^{2p}$ et modulo $u^p$.

\bigskip

Pour tout rationnel $q \in \frac 1 e \N$, on définit sur le site 
$(T_1)_\SYN$ un faisceau $\tilde \J_{(2)}^{[q]}$ par la formule :
$$\tilde \J_{(2)}^{[q]} = \frac {\tilde \J^{[q]}} {\tilde \J^{[q]} 
\cap u^{2p} \tilde \O^\st}.$$
On définit également $\tilde \O^\st_{(2)} = \tilde \J_{(2)}^{[0]}$. 
Finalement, on vérifie que les applications $\phi_s$ et $N$ passent au 
quotient et définissent des opérateurs encore notés $\phi_s$ et $N$ sur 
ces nouveaux faisceaux.

\bigskip

\noindent
{\it Remarque.} On peut vérifier que si l'on note $\tilde E_{(2)}$ la 
réduction de la base $E_1$ modulo $u^{2p}$, on a :
$$\tilde \J_{(2)}^{[s]}(U) = H^0 ((U/\tilde E_{(2)})_\cris, 
\J^{[s]}_{U/\tilde E_{(2)}}) = H^0 ((U/\tilde E_{(2)})_\CRIS, 
\J^{[s]}_{U/\tilde E_{(2)}})$$
pour tout entier $s$ et tout log-schéma $U$ fin et localement de type 
fini sur $T_1$. En outre, on a également :
$$H^i(X_\syn, \J_{(2)}^{[s]}) = H^i(X_\SYN, \J_{(2)}^{[s]}) = H^i 
((X/\tilde E_{(2)})_\cris, \J^{[s]}_{X/\tilde E_{(2)}}) = 
H^i((X/\tilde E_{(2)})_\CRIS, \J^{[s]}_{X/\tilde E_{(2)}})$$
pour tout entier $i$ et tout log-schéma $X$ fin localement de type fini 
sur $E_n$. On a également la relation $\tilde \J_{(2)}^{[q]} = 
u^\delta \tilde \J_{(2)}^{[s]} + \tilde \J_{(2)}^{[s+1]}$ si $eq = es + 
\delta$ est la division euclidienne de $eq$ par $e$.

\bigskip

Il est alors possible, de manière analogue à ce que nous avons fait dans 
le paragraphe \ref{4:sec:osttilde}, de donner des descriptions locales 
très explicites des faisceaux précédents. On retiendra simplement un 
équivalent de la proposition \ref{4:prop:quotj} :

\begin{prop}
\label{4:prop:isomj2}
Supposons $r > 0$. Sur le site $(T_1)_\syn$ les projections canoniques 
induisent des isomorphismes de faisceaux :
$$\xymatrix @C=25pt {
\xytilde \J_{(2)}^{[r]} / \xytilde \J_{(2)}^{[r+2/e]} \ar[r]^-{\sim} & 
\xytilde \J^{[r]} / \xytilde \J^{[r+2/e]} }
\quad \text{et} \quad \xymatrix @C=25pt {
\J_1^{[r]} / \J_1^{[r+1]} \ar[r]^-{\sim} & \xytilde \J_{(2)}^{[r]} /
\xytilde \J_{(2)}^{[r+1]} }.$$
\end{prop}

\begin{preuve}
Elle est tout à fait analogue à celle de la proposition 
\ref{4:prop:quotj}. Remarquons cependant que l'hypothèse $er \leq p-2$ est 
cruciale pour le premier isomorphisme. Le second, quant à lui, utilise 
la majoration moins fine (car on a supposé $r > 0$) $er + e \leq 2p$.
\end{preuve}

\paragraph{Un objet de la catégorie $\Mrtilde_{(2)}$}

Nous nous proposons de prouver, ici, que le quadruplet :
$$(H^i ((X_1)_\syn, \tilde \O^\st_{(2)}), H^i ((X_1)_\syn, \tilde 
\J^{[r]}_{(2)}), \phi_r, N)$$
est un objet de la catégorie $\Mrtilde_{(2)}$. \emph{On suppose à partir
de maintenant que $r > 0$.}

On recopie encore une fois les arguments des paragraphes 
\ref{4:subsec:isomcoh} et \ref{4:subsec:cristilde} afin obtenir le théorème 
suivant :

\begin{theo}
\label{4:th:isomcoh2}
Pour $0 \leq i \leq s \leq \frac{2p-2} e$ (pour $i=s=0$ si $p=2$), on a 
un isomorphisme :
$$\xymatrix @C=30pt {
\xytilde S_{(2)} \otimes_{(\phi), k[u]/u^2} H^i((X_1)_\syn, \xytilde 
\J_{(2)}^{[s]} / \xytilde \J_{(2)}^{[s+2/e]}) \ar[r]^-{\sim}_-{\id \otimes 
\phi_s} & H^i((X_1)_\syn, \xytilde \O^\st_{(2)}) }.$$
\end{theo}

Comme précédemment, pour simplifier, si $\calF$ est un faisceau sur 
$(X_1)_\syn$, on note $H^i(\calF)$ pour $H^i ((X_1)_\syn, \calF)$ et si
cet espace est de dimension finie sur $k$, on note $h^i(\calF)$ sa 
dimension. On a alors :

\begin{theo}
\label{4:th:libre2}
Pour tout $i \leq r$, le groupe de cohomologie $H^i (\tilde 
\O^\st_{(2)})$ est libre de rang fini sur $\tilde S_{(2)}$.
\end{theo}

\begin{preuve}
D'après le théorème \ref{4:th:isomcoh2}, il suffit de montrer que $H^i 
(\tilde \J_{(2)}^{[r]}/\tilde \J_{(2)}^{[r+2/e]})$ est libre de rang 
fini sur $k[u]/u^2$. Or on a un isomorphisme $\tilde \J_{(2)}^{[r]} / 
\tilde \J_{(2)}^{[r+2/e]} \simeq \tilde \J^{[r]} / \tilde \J^{[r+2/e]}$ 
(proposition \ref{4:prop:isomj2}). Il suffit donc de prouver que 
$H^i (\tilde \J^{[r]} / \tilde \J^{[r+2/e]})$ est libre de rang fini 
sur $k[u]/u^2$.

\medskip

Notons $\calM = H^i(\tilde \O^\st)$, c'est un $k[u]/u^p$-module libre 
de rang fini, disons $d$, d'après le théorème \ref{4:th:obtilde}. Pour 
tout $q \in \frac 1 e \N$, $q \geq r$, le morphisme $H^i(\tilde 
\J^{[q]}) \to \calM$ est injectif d'après les propositions 
\ref{4:prop:surj} et \ref{4:prop:inj}. Notons $\Fil^q \calM$ son image. On 
obtient ainsi une suite décroissante de sous-$\tilde S$-modules de 
$\calM$.

Puisque la multiplication par $u$ se factorise par $\Fil^r \calM \to 
\Fil^{r+1/e}\calM \to \Fil^r \calM$, on a $u \Fil^r \calM \subset 
\Fil^{r+1/e}\calM$.
Par ailleurs, à nouveau la proposition \ref{4:prop:surj} nous dit que le 
quotient $\Fil^r \calM / \Fil^{r+1/e}\calM$ est un $k$-espace vectoriel 
de dimension $d$. Or il en est de même de $\Fil^r \calM / u \Fil^r 
\calM$. 
Il en résulte que les $k$-espaces vectoriels $\Fil^{r+1/e}\calM$ et $u 
\Fil^r \calM$ ont même dimension. L'inclusion trouvée précédemment 
prouve alors qu'ils sont égaux.

\medskip

De même en remplaçant $r$ par $r+1/e$ (et en vérifiant que la 
proposition \ref{4:prop:surj} s'applique encore), on obtient 
$\Fil^{r+2/e}\calM = u \Fil^{r+1/e} \calM = u^2 \Fil^r \calM$. 
L'inégalité $er \leq p-2$ et l'inclusion $u^{er} \calM \subset 
\Fil^r\calM$ assurent que le quotient $\Fil^r \calM / \Fil^{r+2/e}\calM 
= \Fil^r \calM / u^2 \Fil^r \calM$ est libre de rang $d$ sur $k[u]/u^2$.

Par ailleurs, on a une suite exacte longue :
$$\xymatrix @R=25pt {
H^i(\xytilde \J^{[r+2/e]}) \ar[r] & H^i(\xytilde \J^{[r]}) \ar[r] & 
H^i(\xytilde \J^{[r]} / \xytilde \J^{[r+2/e]}) \ar[r] & H^{i+1}(\xytilde 
\J^{[r+2/e]}) \ar[r] & H^{i+1}(\xytilde \J^{[r]}) }$$
et d'après la proposition \ref{4:prop:surj} les première et dernière 
flèches sont injectives. On en déduit que $H^i(\tilde \J^{[r]} / 
\tilde \J^{[r+2/e]})$ s'identifie au quotient $\Fil^r \calM / 
\Fil^{r+2/e}\calM$ et donc qu'il est libre de rang fini sur $k[u]/u^2$.
\end{preuve}

\bigskip

\noindent
{\it Remarque.} La preuve précédente implique $h^i (\tilde \O^\st_{(2)}) 
= p h^i(\tilde \J^{[r+2/e]} / \tilde \J^{[r]}) = 2p h^i(\tilde 
\J^{[r+1/e]} / \tilde \J^{[r]}) = 2p h^i(\bar \O^\st)$.

\bigskip

Avec cette dernière égalité, on peut refaire la démonstration de la 
proposition \ref{4:prop:ukost} et obtenir ainsi :

\begin{prop}
\label{4:prop:ukost2}
Pour tous entiers $k$ et $i < \frac {2p} e$ (on impose $i=0$ si $p=e=2$),
on a $H^i(u^k \tilde \O_{(2)}^\st) = u^k H^i(\tilde \O_{(2)}^\st)$.
\end{prop}

\bigskip

Sans surprise, on dispose d'un analogue de la proposition 
\ref{4:prop:surj} dans cette nouvelle situation :

\begin{prop} 
\label{4:prop:surj2}
Soient $i$ un entier et $q \in \frac 1 e \N$ vérifiant $0 \leq i \leq q 
< E(\frac p e) + 1$ (et $i=q=0$ si $p=2$ et $e=1$) où $E(\frac p e)$ 
désigne la partie entière de $\frac p e$. On a des suites exactes 
courtes : 
$$\xymatrix @C=25pt { 0 \ar[r] & H^i(\xytilde \J_{(2)}^{[q+1/e]}) \ar[r] & 
H^i(\xytilde \J_{(2)}^{[q]}) \ar[r] & H^i(\xytilde \J_{(2)}^{[q]} / \xytilde 
\J_{(2)}^{[q+1/e]}) \ar[r] & 0 }.$$
En outre, si $s$ désigne la partie entière de $q$, on a également un 
début de suite exacte : 
$$\xymatrix @C=25pt { 0 \ar[r] & H^{s+1}(\xytilde \J_{(2)}^{[q+1/e]}) 
\ar[r] & H^{s+1}(\xytilde \J_{(2)}^{[q]}) \ar[r] & H^{s+1}(\xytilde 
\J_{(2)}^{[q]} / \xytilde \J_{(2)}^{[q+1/e]})}.$$
\end{prop}

\begin{preuve}
La démonstration est très proche de celle de la proposition 
\ref{4:prop:surj}.

\medskip

On commence par remarquer que puisque $r > 0$, on a $e 
\leq p-2$ et $E(\frac p e) + 1 \leq \frac{2p-2} e$. Cela implique que
si $s$ désigne la partie entière de $q$, on a les identifications 
$\tilde \J_{(2)}^{[q]} / \tilde \J_{(2)}^{[q+1/e]} \simeq \tilde 
\J_{(2)}^{[s]} / \tilde \J_{(2)}^{[s+1/e]} \simeq \tilde \J^{[s]} / 
\tilde \J^{[s+1/e]}$.
Ainsi, on obtient $h^i(\tilde \J_{(2)}^{[q]} / \tilde 
\J_{(2)}^{[q+1/e]}) = h^i(\tilde \J^{[s]} / \tilde \J^{[s+1/e]}) = 
h^i(\bar \O^\st)$ pour tout entier $i$. (La dernière égalité résulte du 
théorème \ref{4:th:isomcoh}.)

Par ailleurs, une adaptation simple du lemme \ref{4:lem:noyuj} fournit la 
suite exacte :
$$\xymatrix @C=25pt {
0 \ar[r] & H^i(\xybar \O^\st) \ar[r] & H^i(\xytilde \J_{(2)}^{[q]}) \ar[r]^u 
& H^i(u \xytilde \J_{(2)}^{[q]}) \ar[r] & 0 }$$
et donc l'égalité $h^i(\tilde \J_{(2)}^{[q]}) = h^i(u \tilde 
\J_{(2)}^{[q]}) + h^i(\bar \O^\st)$.
De même en adaptant le lemme \ref{4:lem:exactj}, on obtient la suite 
exacte de faisceaux :
$$\xymatrix @C=30pt {
0 \ar[r] & u \xytilde \J_{(2)}^{[q]} \ar[r] & \xytilde \J_{(2)}^{[q+1/e]} 
\ar[r] & \xybar \J^{[s+1]} \ar[r] & 0 }$$
et puisque $R \alpha_\star \bar \J^{[s+1]} = 0 \to \cdots \to 0 \to 
\omega^{s+1}_{\bar X / \bar E} \to \omega^{s+2}_{\bar X / \bar E} \to 
\cdots$, on a $H^j (\bar \J^{[s+1]}) = 0$ pour tout $j \leq s$ puis $H^i 
(u \tilde \J^{[s]}) = H^i (\tilde \J^{[s+1/e]})$. En particulier $h^i (u 
\tilde \J^{[s]}) = h^i (\tilde \J^{[s+1/e]})$ d'où il vient $h^i(\tilde 
\J_{(2)}^{[q]}) = h^i(\tilde \J_{(2)}^{[q+1/e]}) + h^i(\tilde 
\J_{(2)}^{[q]} / \tilde \J_{(2)}^{[q+1/e]})$. On termine alors la 
démonstration en raisonnant par récurrence sur $i$.
\end{preuve}

\bigskip

On a finalement la proposition :

\begin{prop}
\label{4:prop:inj2}
Pour tout $i \leq r$, l'application $H^i(\tilde \J_{(2)}^{[r]}) \to 
H^i(\tilde \O^\st_{(2)})$ est injective.
\end{prop}

\begin{preuve}
Il est possible d'adapter la démonstration de la proposition 
\ref{4:prop:inj}, mais nous pouvons également déduire l'énoncé de la 
proposition \ref{4:prop:inj}.
En effet, on a le diagramme suivant :
$$\xymatrix @C=25pt @R=15pt {
& H^i(u^p \xytilde \O^\st_{(2)}) \ar[r] \ar@{=}[d] & H^i(\xytilde 
\J^{[r]}_{(2)}) \ar[r] \ar[d] & H^i(\xytilde \J^{[r]}) \ar[d] \\
0 \ar[r] & H^i(u^p \xytilde \O^\st_{(2)}) \ar[r] & H^i(\xytilde 
\O^\st_{(2)}) \ar[r] & H^i(\xytilde \O^\st) \ar[r] & 0 } $$
et la suite du bas est exacte d'après la proposition \ref{4:prop:ukost2} 
(notez que la multiplication par $u^p$ sur $\tilde \O^\st_{(2)}$ 
identifie $\tilde \O^\st$ et $u^p \tilde \O^\st_{(2)}$).
Par la proposition \ref{4:prop:inj}, la flèche verticale de droite est 
injective. On vérifie facilement qu'il en est alors forcément de même
de la flèche verticale centrale. D'où la proposition.
\end{preuve}

\bigskip

On en déduit finalement le théorème :

\begin{theo}
\label{4:th:obtilde2}
Pour tout $i \leq r$, le quadruplet $(H^i ((X_1)_\syn, \tilde 
\O_{(2)}^\st), H^i ((X_1)_\syn, \tilde \J_{(2)}^{[r]}), \phi_r, 
N)$ définit un objet de la catégorie $\Mrtilde_{(2)}$.
\end{theo}

\paragraph{Un objet de la catégorie $\Mr$}

On veut ici enfin prouver que $(H^i (\O^\st_1), H^i (\J^{[r]}_1), 
\phi_r, N)$ est un objet de $\Mr$ pour tout $i \leq r$. Pour cela, on
reprend à nouveau les arguments précédents.

\bigskip

On commence par prouver que $H^i (\O^\st_1)$ est un module libre (de 
rang fini) sur $S_1$. D'après le théorème \ref{4:th:isomcohs1}, 
il suffit de prouver que $H^i (\J_1^{[r]}/\J_1^{[r+1]})$ est libre sur
$k[u]/u^e$. Or le faisceau $\J_1^{[r]}/\J_1^{[r+1]}$ s'identifie 
à $\tilde \J_{(2)}^{[r]} / \tilde \J_{(2)}^{[r+1]}$ (proposition 
\ref{4:prop:isomj2}) et il suffit donc de prouver que $H^i(\tilde 
\J_{(2)}^{[r]} / \tilde \J_{(2)}^{[r+1]})$ est libre de rang fini sur
$k[u]/u^e$. Pour cela, on adapte facilement les arguments de la 
preuve du théorème \ref{4:th:libre2} en remplaçant les références aux 
propositions \ref{4:prop:surj} et \ref{4:prop:inj} respectivement par des
références aux proposition \ref{4:prop:surj2} et \ref{4:prop:inj2}.

\medskip

\begin{lemme}
\label{4:lem:proj2}
Pour tous entiers $i$ et $s$ tels que $0 \leq i \leq s < \frac p e$, la 
flèche $H^i(\J^{[s]}_1) \to H^i(\tilde \J^{[s]}_{(2)})$ est surjective.
\end{lemme}

\begin{preuve}
On raisonne par récurrence sur $s$. Pour le cas $s=0$, on a les 
isomorphismes suivants :
\begin{eqnarray*}
H^i(\O^\st_1) & \simeq & S_1 \otimes_{(\phi), k[u]/u^e} H^i(\J^{[i]}_1 
/ \J^{[i+1]}_1) \\
H^i(\tilde \O^\st_{(2)}) & \simeq & \tilde S_{(2)} \otimes_{(\phi), 
k[u]/u^2} H^i(\tilde \J^{[i]}_{(2)} / \tilde \J^{[i+2/e]}_{(2)}).
\end{eqnarray*}
On vérifie que la flèche $H^i(\O^\st_1) \to H^i(\tilde \O^\st_{(2)})$ 
s'obtient à partir de la projection $S_1 \to \tilde S_{(2)}$ et le 
morphisme naturel $H^i(\J^{[i]}_1 / \J^{[i+1]}_1) \to H^i(\tilde 
\J^{[i]}_{(2)} / \tilde \J^{[i+2/e]}_{(2)})$. Il suffit donc de prouver 
que ce dernier morphisme est surjectif.
Or, d'une part, le faisceau $\J^{[i]}_1 / \J^{[i+1]}_1$ s'identifie à
$\tilde \J^{[i]}_{(2)} / \tilde \J^{[i+1]}_{(2)}$ (proposition 
\ref{4:prop:isomj2}) et, d'autre part, on a le diagramme commutatif 
suivant :
$$\xymatrix @C=25pt @R=15pt {
0 \ar[r] & H^i(\xytilde \J^{[i+1]}_{(2)}) \ar[r] \ar[d] & H^i(\xytilde 
\J^{[i]}_{(2)}) \ar[r] \ar@{=}[d] & H^i(\xytilde \J^{[i]}_{(2)} / \xytilde 
\J^{[i+1]}_{(2)} \ar[r]) \ar[d] & 0 \\
0 \ar[r] & H^i(\xytilde \J^{[i+2/e]}_{(2)}) \ar[r] & H^i(\xytilde 
\J^{[i]}_{(2)}) \ar[r] & H^i(\xytilde \J^{[i]}_{(2)} / \xytilde 
\J^{[i+2/e]}_{(2)}) \ar[r] & 0 }$$
où les lignes sont exactes : l'injectivité provient de la proposition 
\ref{4:prop:surj2} et la surjectivité provient de l'injectivité analogue
sur les $H^{i+1}$ (toujours conséquence de la même proposition). On 
en déduit directement la surjectivité de la flèche verticale de droite, 
ce qui conclut le cas $s=0$.

\medskip

On procède ensuite par récurrence sur $s$. On considère le diagramme 
commutatif :
$$\xymatrix @C=25pt @R=15pt {
& H^i(\xytilde \J^{[s+1]}_1) \ar[r] \ar[d] & H^i(\xytilde \J^{[s]}_1) \ar[r] 
\ar[d] & H^i(\xytilde \J^{[s]}_1 / \xytilde \J^{[s+1]}_1) \ar@{=}[d] \\
0 \ar[r] & H^i(\xytilde \J^{[s+1]}_{(2)}) \ar[r] & H^i(\xytilde 
\J^{[s]}_{(2)}) \ar[r] & H^i(\xytilde \J^{[s]}_{(2)} / \xytilde 
\J^{[s+1]}_{(2)}) }$$
La flèche verticale de droite est un isomorphisme par la proposition 
\ref{4:prop:isomj2}. La flèche verticale du milieu est surjective par
hypothèse de récurrence. Une chasse au diagramme prouve facilement que 
la flèche verticale de gauche est aussi surjective, ce qui conclut.
\end{preuve}

\bigskip

\noindent
{\it Remarque.} Le lemme est également vrai pour $s \in \frac 1 e \N$, 
$0 \leq i \leq s < E(\frac p e) + 1$. Ce raffinement se démontre de 
manière analogue en choisissant un pas de $\frac 1 e$ (au lieu de $1$) 
dans la récurrence.

\bigskip

On parvient finalement au but de tout ce paragraphe :

\begin{theo}
\label{4:th:ob}
Pour tout $i \leq r$, le quadruplet $(H^i ((X_1)_\syn, \O_1^\st), H^i 
((X_1)_\syn, \J_1^{[r]}), \phi_r, N)$ définit un objet de la catégorie 
$\Mr$.
\end{theo}

\begin{preuve}
Comme précédemment, il ne reste plus qu'à prouver que la flèche $H^i( 
\J_1^{[r]}) \to H^i(\J_1^{[r]} / \J_1^{[r+1]})$ est surjective et que la 
flèche $H^i(\J_1^{[r]}) \to H^i(\O^\st_1)$ est injective.

\medskip

Pour le premier point, on considère le carré commutatif suivant :
$$\xymatrix {
H^i(\J_1^{[r]}) \ar[r] \ar[d] & H^i(\J_1^{[r]} / \J_1^{[r+1]}) \ar[d] \\
H^i(\xytilde \J_{(2)}^{[r]}) \ar[r] & H^i(\xytilde \J_{(2)}^{[r]} / \xytilde 
\J_{(2)}^{[r+1]}) }$$
La flèche de droite est un isomorphisme (proposition \ref{4:prop:isomj2}).
La flèche de gauche est surjective (lemme \ref{4:lem:proj2}), et celle 
du bas l'est également (proposition \ref{4:prop:surj2}). On en déduit que 
celle du haut l'est aussi.

\medskip

Pour le second point, on reprend les arguments de la preuve de la 
proposition \ref{4:prop:inj2}. On raisonne à partir du diagramme 
suivant :
$$\xymatrix @C=25pt @R=15pt {
0 \ar[r] & K \ar[r] \ar@{=}[d] & \J_1^{[r]} \ar[r] \ar[d] & \xytilde 
\J_{(2)}^{[r]} \ar[r] \ar[d] & 0 \\
0 \ar[r] & K \ar[r] & \O^\st_1 \ar[r] & \xytilde \O^\st_{(2)} \ar[r] & 0 }$$
où $K = u^p \O_1^\st$ et où les deux suites horizontales sont exactes. 
Il donne lieu à un nouveau diagramme :
$$\xymatrix @C=25pt @R=15pt {
& H^i(K) \ar[r] \ar@{=}[d] & H^i(\J_1^{[r]}) \ar[r] \ar[d] & H^i(\xytilde 
\J_{(2)}^{[r]}) \ar[d] \\
0 \ar[r] & H^i(K) \ar[r] & H^i(\O^\st_1) \ar[r] & H^i(\xytilde 
\O^\st_{(2)})}$$
La suite du bas est exacte à gauche car, par le lemme \ref{4:lem:proj2},
la flèche $H^{i-1}(\O^\st_1) \to H^{i-1}(\tilde \O^\st_{(2)})$ est 
surjective. Mais la flèche verticale de droite est injective 
(proposition \ref{4:prop:inj2}). On en déduit que la flèche verticale du 
milieu l'est également. Ceci termine la preuve.
\end{preuve}

\subsubsection{Le cas $r=0$}
\label{4:subsec:crisr0}
Lorsque $r=0$, la condition $er \leq p-2$ est automatique et n'impose
aucune borne sur $e$. Pour obtenir le théorème \ref{4:th:ob}, dans ce cas, 
on a besoin d'introduire plus de catégories-étapes.

\paragraph{Plan de la preuve}

On reprend les constructions du début du paragraphe \ref{4:subsec:crise1} 
en remplaçant \og $(2)$ \fg\ par \og $(t)$ \fg.

\medskip

Pour tout entier $1 \leq t \leq e$, on commence par définir une 
catégorie $\Mrtilde_{(t)}$ d'objets modulo $u^{tp}$. On considère pour
cela l'anneau $\tilde S_{(t)} = k\cro u / u^{tp}$ que l'on munit d'un 
Frobenius $\phi$ semi-linéaire envoyant $u^i$ sur $u^{pi}$ et d'un 
opérateur de monodromie $k$-linéaire envoyant $u^i$ sur $-i u^i$. 
Un objet de $\Mrtilde_{(t)}$ est la donnée de :
\begin{enumerate}
\item un $\tilde S_{(t)}$-module $\tilde \calM$ libre de rang fini ;
\item un sous-module $\Fil^r \tilde \calM$ de $\tilde \calM$ contenant 
$u^{er} \tilde \calM$ ;
\item une flèche $\phi$-semi-linéaire $\phi_r : \Fil^r \tilde \calM \to 
\tilde \calM$ telle que l'image de $\phi_r$ engendre $\tilde \calM$ en 
tant que $\tilde S_{(t)}$-module ;
\item une application $k$-linéaire $N : \tilde \calM \to \tilde \calM$ 
telle que :
\begin{itemize}
\item pour tout $\lambda \in \tilde S_{(t)}$ et tout $x \in \tilde 
\calM$, $N\pa{\lambda x} = N\pa \lambda x + \lambda N\pa x$
\item $u^e N(\Fil^r \tilde \calM) \subset \Fil^r \tilde \calM$
\item le diagramme suivant commute :
\[\xymatrix @C=50pt {
\Fil^r \xytilde \calM \ar[r]^{\phi_r} \ar[d]_{u^e N} & \xytilde \calM 
\ar[d]^{c_{(t),\pi} N} \\
\Fil^r \xytilde \calM \ar[r]^{\phi_r} & \xytilde \calM}\]
où $c_{(t),\pi}$ est la réduction de $c$ dans $\tilde S_{(t)}$.
\end{itemize}
\end{enumerate}
Les morphismes de $\Mrtilde_{(t)}$ sont les applications $\tilde 
S_{(t)}$-linéaires qui commutent à toutes les structures. On prouve 
comme dans le paragraphe \ref{4:sec:mrtilde} que les catégories  
$\Mrtilde_{(t)}$ sont toutes équivalentes à $\Mr$ \emph{via} les 
foncteurs de réduction modulo $u^{tp}$.

\medskip

D'autre part, pour tout rationnel $q \in \frac 1 e \N$, on définit sur 
le site $(T_1)_\SYN$ un faisceau $\tilde \J_{(t)}^{[q]}$ par :
$$\tilde \J_{(t)}^{[q]} = \frac {\tilde \J^{[q]}} {\tilde \J^{[q]} 
\cap u^{tp} \tilde \O^\st}.$$
On pose $\tilde \O^\st_{(t)} = \tilde \J_{(t)}^{[0]}$ et on vérifie 
que les applications $\phi_s$ et $N$ passent au quotient et définissent 
des opérateurs encore notés $\phi_s$ et $N$ sur ces nouveaux faisceaux.

\medskip

On va prouver par récurrence sur $t$ la proposition suivante (on 
rappelle qu'ici $r = 0$) :

\begin{prop}
\label{4:prop:obtildet}
Pour tout $t \leq e$, le quadruplet $(H^0 (\tilde \O_{(t)}^\st), H^0 
(\tilde \O_{(t)}^\st), \phi_0, N)$ définit un objet de la catégorie 
$\MOtilde_{(t)}$ et la suite :
$$\xymatrix @C=25pt { 0 \ar[r] & H^0(\xytilde \J_{(t)}^{[q+1/e]}) \ar[r] & 
H^0(\xytilde \J_{(t)}^{[q]}) \ar[r] & H^0(\xytilde \J_{(t)}^{[q]} / \xytilde 
\J_{(t)}^{[q+1/e]}) \ar[r] & 0 }$$
est exacte pour $q \in \frac 1 e \N$, $q \leq \frac t e$.
\end{prop}

On déduira ensuite du cas $t = e$ un équivalent du théorème \ref{4:th:ob} 
à savoir :

\begin{theo}
\label{4:th:obr0}
Pour tout $t \leq e$, le quadruplet $(H^0 (\O_1^\st), H^0 (\O_1^\st), 
\phi_0, N)$ définit un objet de la catégorie $\MO$.
\end{theo}

\paragraph{La récurrence}

Le but de cette sous-partie est de prouver la proposition 
\ref{4:prop:obtildet}. Pour cela, comme nous l'avons déjà dit, on 
raisonne par récurrence sur $t$ et on suit pas à pas la 
démonstration de la partie \emph{Un objet de la catégorie 
$\Mrtilde_{(2)}$} du paragraphe \ref{4:subsec:crise1}. On redonne 
rapidement les grandes étapes.

\medskip

On considère un entier $t$ compris entre $2$ et $e$.
Les compatibilités entre les opérateurs ne posant pas de problème, il 
suffit de montrer que $H^0 (\tilde \O_{(t)}^\st)$ est libre de rang fini 
sur $\tilde S_{(t)}$ et que l'image de $\phi_0$ engendre tout $H^0 
(\tilde \O_{(t)}^\st)$. Or comme précédemment, on dispose du théorème
suivant :

\begin{theo}
\label{4:th:isomcoht}
Avec les notations précédentes, on a un isomorphisme :
$$\xymatrix @C=30pt {
\xytilde S_{(t)} \otimes_{(\phi), k[u]/u^t} H^0(\xytilde 
\O^\st_{(t)} / \xytilde \J_{(t)}^{[t/e]}) \ar[r]^-{\sim}_-{\id \otimes 
\phi_0} & H^0(\xytilde \O^\st_{(t)}) }.$$
\end{theo}

\noindent
Il suffit donc de prouver que $H^0(\tilde \O^\st_{(t)} / \tilde 
\J_{(t)}^{[t/e]})$ est libre de rang fini sur $k[u]/u^t$ et que 
la suite :
$$\xymatrix @C=25pt { 0 \ar[r] & H^0(\xytilde \J_{(t)}^{[q+1/e]}) \ar[r] & 
H^0(\xytilde \J_{(t)}^{[q]}) \ar[r] & H^0(\xytilde \J_{(t)}^{[q]} / \xytilde 
\J_{(t)}^{[q+1/e]}) \ar[r] & 0 }$$
est exacte pour $q \in \frac 1 e \N$, $q \leq t$. (On remarque que 
seule la surjectivité n'est \emph{a priori} pas claire).

\medskip

On dispose d'un équivalent de la proposition \ref{4:prop:isomj2} qui se 
démontre de façon tout à fait analogue :

\begin{prop}
Sur le site $(T_1)_\syn$ la projection canonique induit un isomorphisme
de faisceaux :
$$\xymatrix @C=25pt {
\xytilde \O^\st_{(t)} / \xytilde \J_{(t)}^{[t/e]} \ar[r]^-{\sim} & 
\xytilde \O^\st_{(t-1)} / \xytilde \J_{(t-1)}^{[t/e]} }. $$
\end{prop}

\noindent
{\it Remarque.} Pour la preuve de la proposition précédente, on utilise 
l'inégalité $t \leq p(t-1)$ qui est bien vérifiée lorsque $t \geq 
2$.

\bigskip

\begin{lemme}
Le module $H^0 (\tilde \O^\st_{(t)})$ est libre de rang fini sur $\tilde 
S_{(t)}$.
\end{lemme}

\begin{preuve}
Elle est assez semblable à la preuve du théorème \ref{4:th:libre2}.

Il suffit de démontrer que $H^0(\tilde \O^\st_{(t)} / \tilde 
\J_{(t)}^{[t/e]})$ est libre sur $k[u]/u^t$. Par ailleurs d'après
la proposition précédente, le faisceau $\tilde \O^\st_{(t)} / \tilde
\J_{(t)}^{[t/e]}$ s'identifie à $\tilde \O^\st_{(t-1)} / \tilde
\J_{(t-1)}^{[t/e]}$.

On utilise à ce niveau l'hypothèse de récurrence qui assure que 
$\calM = H^0 (\tilde \O^\st_{(t-1)})$ muni de ses structures 
supplémentaires est un objet de $\MOtilde_{(t-1)}$. C'est en particulier 
un $\tilde S_{(t-1)}$-module libre de rang fini, disons $d$. Notons, 
pour $q \in \frac 1 e \N$, $\Fil^q \calM =  H^0 (\tilde 
\J^{[q]}_{(t-1)})$. Ils forment une suite décroissante de sous-$\tilde 
S_{(t-1)}$-modules de $\calM$.

Par ailleurs, on a $u \calM \subset \Fil^{1/e}\calM$, et encore 
l'hypothèse de récurrence fournit la suite exacte :
$$\xymatrix @C=25pt { 0 \ar[r] & \Fil^{1/e} \calM \ar[r] &
\calM \ar[r] & H^0(\xytilde \O_{(t-1)}^\st / \xytilde \J_{(t-1)}^{[1/e]}) 
\ar[r] & 0 } .$$
Par ailleurs on vérifie que pour tout $q \in \frac 1 e \N$, $q <
\frac t e$, la multiplication par $u^q$ induit un isomorphisme entre
les faisceaux $\tilde \O_{(t-1)}^\st / \tilde \J_{(t-1)}^{[1/e]}$ et
$\tilde \J_{(t-1)}^{[q]} / \tilde \J_{(t-1)}^{[q+1/e]}$. On en déduit :
$$\dim_k H^0(\tilde \O_{(t-1)}^\st / \tilde \J_{(t-1)}^{[1/e]})
= \frac 1 t \dim_k H^0(\tilde \O_{(t-1)}^\st / \tilde 
\J_{(t-1)}^{[t/e]}) = d$$
la dernière égalité étant une conséquence du théorème \ref{4:th:isomcoht}
appliqué pour $t-1$. D'autre part, $\dim_k (\calM / u\calM) = d$ d'où il 
vient $\Fil^{1/e}\calM = u \calM$. On montre de même par récurrence que 
$\Fil^q\calM = u^{eq} \calM$ pour tout $q \in \frac 1 e \N$, $q \leq 
\frac t e$.

\medskip

Finalement, on a une suite exacte longue :
$$\small \xymatrix @C=20pt {
0 \ar[r] & H^0(\xytilde \J_{(t-1)}^{[t/e]}) \ar[r] & H^0(\xytilde 
\O_{(t-1)}^\st) \ar[r] & H^0(\xytilde \O_{(t-1)}^\st / \xytilde 
\J_{(t-1)}^{[t/e]}) \ar[r] & H^1(\xytilde \J_{(t-1)}^{[t/e]}) \ar[r] & 
H^1(\xytilde \O_{(t-1)}^\st) }.$$
D'après l'hypothèse de récurrence, toutes les flèches $H^1(\tilde 
\J_{(t-1)}^{[q+1/e]}) \to H^1(\tilde \J_{(t-1)}^{[q]} )$ sont 
injectives pour $q \leq \frac {t-1} e$ et donc la flèche $H^1(\tilde 
\J_{(t-1)}^{[t/e]}) \to H^1(\tilde \O_{(t-1)}^\st)$ l'est aussi. On en 
déduit que $H^0(\tilde \O_{(t-1)}^\st / \tilde \J_{(t-1)}^{[t/e]})$ 
s'identifie à $\calM / \Fil^{t/e} \calM = \calM / u^t \calM$ qui est un 
bien un $k[u]/u^t$-module libre de rang fini (en l'occurrence $d$).
\end{preuve}

\bigskip

Les propositions \ref{4:prop:ukost2} et \ref{4:prop:surj2} ont des 
équivalents transparents dans ce nouveau contexte qui se démontrent de 
façon analogue. On conclut comme cela la récurrence.

\paragraph{La fin de la preuve}

Pour déduire le théorème \ref{4:th:obr0} du cas $t = e$ de la 
proposition \ref{4:prop:obtildet}, on dégage tout d'abord le
lemme suivant (preuve analogue à celle de la proposition 
\ref{4:prop:quotj}) :

\begin{lemme}
Sur le site $(T_1)_\syn$ la projection canonique induit un isomorphisme
de faisceaux :
$$\xymatrix @C=25pt {
\O^\st_1 / \J_1^{[1]} \ar[r]^-{\sim} & \xytilde \O^\st_{(e)} / \xytilde 
\J_{(e)}^{[1]} }. $$
\end{lemme}

\medskip

En combinant ce précédent lemme au théorème \ref{4:th:isomcohs1}, on 
prouve la liberté sur $S_1$ de $H^0(\O^\st_1)$. Par suite, on démontre 
que la flèche $H^0(\O_1^\st) \to H^0(\tilde \O_{(e)}^\st)$ est 
surjective en utilisant les mêmes arguments que ceux présentés dans 
l'étape d'initialisation de la récurrence du lemme \ref{4:lem:proj2}. 
Finalement,
en recopiant la première partie de la preuve du théorème \ref{4:th:ob},
on parvient à prouver la surjectivité du morphisme $H^0(\O_1^\st) 
\to H^0(\O_1^\st / \J_1^{[1]})$. Cela conclut la preuve du théorème
\ref{4:th:obr0}.

\subsection{Dévissages}
\label{4:sec:crisdev}
Le but de cette partie est de déduire à partir du cas $n=1$ traité 
précédemment le cas $n$ quelconque. Exactement, si $X_n$ désigne un 
log-schéma propre et log-lisse sur $T_n$, et que $X_1 = X_n \times_{T_n} 
T_1$ est du type de Cartier (si $n=1$, on suppose en outre que $X_1$ 
admet un relèvement propre et log-lisse sur $T_2$), nous allons 
démontrer le théorème suivant :

\begin{theo}
\label{4:th:obn}
Pour tout $i < r$ et tout entier $n$, le quadruplet :
$$(H^i ((X_n)_\syn, \O_n^\st), H^i ((X_n)_\syn, \J_n^{[r]}), \phi_r, 
N)$$ définit un objet de la catégorie $\Mr$.
\end{theo}

\noindent
{\it Remarque.} On rappelle que le morphisme $\phi_r$ n'est pas défini 
sur $(X_n)_\syn$ mais seulement sur $(X_{n+r})_\syn$ mais que cela n'est 
pas grave du fait des identifications canoniques $H^i ((X_n)_\syn,  
\O_n^\st) \simeq H^i ((X_{n+r})_\syn, \O_n^\st)$ et $H^i 
((X_n)_\syn, \J_n^{[r]}) \simeq H^i ((X_{n+r})_\syn, 
\J_n^{[r]})$.

\bigskip

\begin{preuve}
On raisonne par récurrence sur $n$. L'initialisation est donnée par
le théorème \ref{4:th:ob}. Pour l'hérédité, on considère les suites
exactes courtes suivantes :
$$\xymatrix @C=30pt @R=15pt {
0 \ar[r] & \J_1^{[r]} \ar[r]^{p^n} & \J_{n+1}^{[r]} \ar[r] & \J_n^{[r]} 
\ar[r] & 0 \\
0 \ar[r] & \O_1^\st \ar[r]^{p^n} & \O_{n+1}^\st \ar[r] & \O_n^\st \ar[r] 
& 0. }$$
Elles fournissent deux suites exactes longues de cohomologie qui 
s'insèrent dans le diagramme commutatif suivant :
$$\xymatrix @C=30pt @R=25pt {
H^{i-1}(\J_n^{[r]}) \ar[r] \ar[d]_{\phi_r} & H^i(\J_1^{[r]}) \ar[r] 
\ar[d]_{\phi_r} & H^i(\J_{n+1}^{[r]}) \ar[r] \ar[d]_{\phi_r} & 
H^i(\J_n^{[r]}) \ar[r] \ar[d]_{\phi_r} & H^{i+1}(\J_1^{[r]}) 
\ar[d]_{\phi_r} \\
H^{i-1}(\O_n^\st) \ar[r] & H^i(\O_1^\st) \ar[r] & H^i(\O_{n+1}^\st)
\ar[r] & H^i(\O_n^\st) \ar[r] & H^{i+1}(\O_1^\st)
}$$
où tous les groupes de cohomologie sont calculés sur le site 
$(X_n)_\syn$. Par hypothèse de récurrence, les deux quadruplets 
$(H^{i-1}(\O_n^\st), H^{i-1}(\J_n^{[r]}), \phi_r, N)$ et 
$(H^i(\O_1^\st), H^i(\J_1^{[r]}), \phi_r, N)$ sont des objets de la
catégorie $\Mr$. Comme celle-ci est abélienne, il en est de même de
leur image que l'on note $(\calM', \Fil^r \calM', \phi_r, N)$. De même,
les deux quadruplets $(H^i(\O_n^\st), H^i(\J_n^{[r]}), \phi_r,
N)$ et $(H^{i+1}(\O_1^\st), H^{i+1}(\J_1^{[r]}), \phi_r, N)$ sont des
objets de $\Mr$, le premier en vertu de l'hypothèse de récurrence et
le second par le théorème \ref{4:th:ob}. Leur noyau $(\calM'', \Fil^r
\calM'', \phi_r, N)$ est donc aussi objet de $\Mr$.

En outre, on dispose d'un diagramme :
$$\xymatrix @C=30pt @R=25pt {
0 \ar[r] & \Fil^r \calM' \ar[r] \ar[d]_{\phi_r} & H^i(\J_{n+1}^{[r]})
\ar[r] \ar[d]_{\phi_r} & \Fil^r \calM'' \ar[r] \ar[d]_{\phi_r} & 0 \\
0 \ar[r] & \calM' \ar[r] & H^i(\O_{n+1}^\st) \ar[r] & \calM''
\ar[r] & 0 }$$
où les deux lignes horizontales sont exactes. Une adaptation directe du
lemme 2.3.1.2 de \cite{breuil-duke} entraîne alors que
$(H^i(\O_{n+1}^\st), H^i(\J_{n+1}^{[r]}), \phi_r, N)$ est un objet de
$\Mr$, ce qui achève la récurrence.
\end{preuve}

\bigskip

\noindent
{\it Remarque.} Il est fort probable que le théorème précédent reste
vrai lorsque $i=r$, mais ce cas particulier échappe à la preuve que
l'on vient de donner.

\subsubsection*{Le cas entier}

Après avoir obtenu un théorème modulo $p^n$ pour tout entier $n$, il
est tentant de passer à la limite projective. Précisément, posons, au
moins pour $i < r$ :
\begin{eqnarray*}
\calM & = & \varprojlim_{n \geq 1} H^i ((X_n)_\syn, \O_n^\st) \\
\Fil^r \calM & = & \varprojlim_{n \geq 1} H^i ((X_n)_\syn, 
\J_n^{[r]}).
\end{eqnarray*}
Les applications $\phi_r$ et $N$ passent à la limite pour finir
respectivement des applications $\Fil^r \calM \to \calM$ et
$\calM \to \calM$ que l'on note encore $\phi_r$ et $N$.

Soit $\calM_\tors$ l'ensemble des éléments de $\calM$ tué par une
puissance de $p$, et $\calM_\free = \calM/\calM_\tors$. On munit sans
problème ces modules d'un $\Fil^r$, d'un $\phi_r$ et d'un $N$, et en
copiant les arguments du paragraphe 4.1 de \cite{breuil-duke}, on
obtient le théorème suivant :

\begin{theo}
\begin{enumerate}
\item[i)] Le module $\calM_\tors$ muni des structures supplémentaires 
est un objet de $\Mr$.
\item[ii)] Le module $\calM_\free$ muni des structures supplémentaires 
est un module fortement divisible\footnote{Pour une définition, on 
pourra se reporter à \cite{breuil-invent}.}
\end{enumerate}
\end{theo}

\section{Calcul de la cohomologie étale}
\label{4:sec:etale}

On fixe toujours un entier $r$ vérifiant $er < p-1$. On se donne de 
plus $X_K$ un schéma (au sens classique) propre et lisse sur $K$ 
et on suppose que $X_K$ admet un modèle propre et semi-stable $X$ sur 
l'anneau des entiers $\O_K$. Le diviseur donné par la fibre spéciale 
fait de $X$ un log-schéma défini sur la base $T = (\spec \O_K, 
\O_K\setminus \{0\})$ propre, log-lisse, et dont la fibre spéciale est 
du type de 
Cartier. Nous sommes donc en situation d'utiliser les résultats de la 
partie précédente. En particulier, si $X_n = X \times_T T_n$, le 
théorème \ref{4:th:obn} s'applique et assure que pour tout $n$ et pour 
tout $i \leq r$ le quadruplet :
$$(H^i ((X_n)_\syn, \O_n^\st), H^i ((X_n)_\syn, \J_n^{[r]}), \phi_r,
N)$$
est un objet de $\Mr$. D'autre part, par le théorème \ref{4:th:ob}, le 
résultat demeure pour $i = r$ lorsque $n = 1$.

Le but de cette partie est de démontrer le théorème \ref{4:th:comp}, 
dont nous précisons l'énoncé :

\begin{theo}
\label{4:th:comprappel}
Pour tout entier $r$ tel que $er < p-1$, pour tout entier $n$ et pour 
tout $0 \leq i < r$ (et aussi $i=r$ si $n=1$), on a un isomorphisme 
canonique de modules galoisiens :
$$\xymatrix { H^i((X_{\xybar K})_\et, \Z/p^n\Z)(r) \ar[r]^-{\sim} &
{\Tst}_\star (H^i ((X_n)_\syn, \O_n^\st), H^i ((X_n)_\syn,
\J_n^{[r]}), \phi_r, N) }.$$
\end{theo}
 
Fixons avant tout quelques notations. Si $L$ est une extension 
algébrique de $K$, définissons $T_L = (\spec \O_L, \O_L\setminus 
\{0\})$ et si $n$ est un entier et $Y$ est un log-schéma sur $T$, 
posons $Y_n = Y \times_T T_n$, $Y_L = Y \times_T T_L$ et $Y_{n,L} = Y_n 
\times_T T_L = Y_L \times_T T_n$.

\medskip

Le premier (et principal) ingrédient de la preuve est un résultat de 
Kato et Tsuji qui s'énonce comme suit :

\begin{theo}
Pour $0 \leq i \leq s \leq p-2$, on a des isomorphismes canoniques 
compatibles à l'action de Galois :
$$\xymatrix {
H^i((X_{n, \xybar K})_\et, s_{n,X_{\xybar K}}^\log(s)) \ar[r]^-{\sim} &
H^i((X_{\xybar K})_\et, \Z/p^n\Z)(s) } .$$
\end{theo}

\noindent
Dans le théorème précédent, $s_{n,X_{\bar K}}^\log(s)$ désigne un 
certain complexe de faisceaux étales sur $X_{n, \bar K}$ construit par 
Kato (voir \cite{kato-vanishing}). Par ailleurs, Breuil démontre (lemme 
3.2.4.3 de \cite{breuil-duke} --- la démonstration est écrite dans le 
cas non ramifié, mais elle fonctionne de la même façon dans le cas 
général) le théorème suivant :

\begin{theo}
Pour tout entier $i$, et tout $s \in \{0, \cdots, p-1\}$, on a des 
isomorphismes canoniques compatibles à l'action de Galois :
$$\xymatrix {
\displaystyle \xyvarinjlim L H^i((X_{n+s, L})_\syn, \S_n^s) 
\ar[r]^-{\sim} 
& H^i((X_{n, \xybar K})_\et, s_{n,X_{\xybar K}}^\log(s))} $$
où la limite inductive est prise sur toutes les extensions finies $L$ de 
$K$.
\end{theo}

\noindent
Ici, $\S_n^s$ désigne un certain faisceau sur le site log-syntomique 
dont le rappel de la définition est l'objet du paragraphe 
\ref{4:subsec:sns}. Forts de cela, il ne reste plus pour conclure qu'à 
prouver :

\begin{prop}
\label{4:prop:isom}
Pour $0 \leq i < r$ (et aussi $i=r$ lorsque $n=1$), on a des 
isomorphismes canoniques compatibles à l'action de Galois :
\begin{equation}
\label{4:eq:isom}
\xymatrix {
\displaystyle \xyvarinjlim L H^i((X_{n+r, L})_\syn, \S_n^r) 
\ar[r]^-{\sim} 
& {\Tst}_\star (H^i ((X_n)_\syn, \O_n^\st), H^i ((X_n)_\syn, 
\J_n^{[r]}), \phi_r, N) } .
\end{equation}
\end{prop}

\noindent
La démonstration de cette proposition est l'objet du paragraphe 
\ref{4:subsec:preuve}. Finalement, le paragraphe \ref{4:subsec:serre} qui 
termine cette partie explique comment on déduit des résultats précédents 
le théorème \ref{4:th:serre}.

\subsection{Les faisceaux $\S_n^s$}
\label{4:subsec:sns}
On note $T^\triv_n$ le log-schéma $\spec \O_K/p^n$ muni de la
log-structure triviale. Comme on avait défini les faisceaux
$\J_n^{[s]}$ sur le site $(T_n)_\syn$ (voir paragraphe 
\ref{4:subsec:ostlocal}), on définit sur le site $(T^\triv_n)_\syn$ des 
faisceaux $\J_n^{\cris, [s]}$ en posant :
$$\J_n^{\cris, [s]}(U) = H^0((U/T^\triv_n)_\cris, \J^{[s]}_{U / 
T^\triv_n}) = H^0((U/T^\triv_n)_\CRIS, \J^{[s]}_{U / T^\triv_n})$$
pour tout $U$ log-syntomique sur $T^\triv_n$. On pose également
$\O_n^\cris = \J_n^{\cris, [0]}$.

Le morphisme naturel $T_n \to T^\triv_n$ est log-syntomique, et donc les
faisceaux précédents définissent par restriction des faisceaux sur le
site $(T_n)_\syn$ encore notés $\J_n^{\cris, [s]}$ et $\O_n^\cris$.
On dispose de descriptions explicites locales des faisceaux
précédents :

\begin{prop}
\label{4:prop:descocris}
En reprenant les notations ($A^\infty$, $P^\infty$ et $W_n^{\cris,\DP}
(A^\infty, P^\infty)$) du paragraphe \ref{4:subsec:ostlocal}, on a un
isomorphisme canonique :
$$\xymatrix{\O^\cris_n (A^\infty, P^\infty) \ar[r]^-{\sim} & W_n^{\cris, 
\DP} (A^\infty, P^\infty)}.$$
Par ailleurs, cet isomorphisme respecte la filtration donnée à gauche
par les $\J^{\cris, [s]}_n (A^\infty, P^\infty)$ et à droite par la
filtration canonique par les puissances divisées.
\end{prop}

\medskip

Comme dans le cas \og $\st$ \fg, on définit pour $s \leq p-1$, des
applications $\phi_s : \J_n^{\cris, [s]} \to \O_n^\cris$. Finalement, on
appelle $\S_n^s$ le noyau de l'application $\phi_s - \id$.

\bigskip

La proposition suivante réunit deux suites exactes importantes à propos
des faisceaux introduits précédemment :

\begin{prop}
\label{4:prop:sns}
Pour tous entiers $n$ et $s$, on a une suite exacte courte de faisceaux 
sur le site $(T_n)_\syn$ :
$$\xymatrix @C=35pt {
0 \ar[r] & \J_n^{\cris, [s]} \ar[r] & \J_n^{[s]} \ar[r]^-{N} &
\J_n^{[s-1]} \ar[r] & 0 } .$$
Pour tout entier $n$ et tout $s \in \{0, \ldots, p-1\}$, on a une suite
exacte courte de faisceaux sur le site $(T_{n+s})_\syn$ :
$$\xymatrix @C=35pt {
0 \ar[r] & \S_n^s \ar[r] & \J_n^{\cris, [s]} \ar[r]^-{\phi_s - \id} &
\O_n^\cris \ar[r] & 0 } .$$
\end{prop}

\begin{preuve}
La première suite exacte résulte des descriptions précédentes si
l'on se rappelle que :
$$N\pa{\frac{X^i}{i!}} = (1+X) \frac{X^{i-1}}{(i-1)!}.$$
On pourra consulter la preuve de la proposition 3.1.3.1 de
\cite{breuil-duke} pour plus de détails.

Pour la seconde suite exacte, il suffit de prouver la surjectivité de
$\phi_s - \id$ et avec les descriptions locales précédentes, on
construit explicitement un antécédent (local pour la topologie 
log-syntomique) à tout élément de $\O_n^\cris (A^\infty, P^\infty)$.
Exactement, la preuve est identique à celle de la proposition 3.1.4.1 de 
\cite{breuil-duke}, sauf le dernier argument qui est remplacé par
celui du lemme \ref{4:lem:filconoyau}.
\end{preuve}

\subsection{La preuve}
\label{4:subsec:preuve}
Le but de ce chapitre est de donner une preuve de la proposition
\ref{4:prop:isom}, ce qui est suffisant comme nous l'avons expliqué, 
pour démontrer le théorème \ref{4:th:comp}. On suit encore une fois de 
très près la démonstration de \cite{breuil-duke} valable pour le cas 
$e=1$.

\medskip

On note $\calM = H^i ((X_n)_\syn, \O_n^\st)$. Pour tout entier $s$, on a
un morphisme (pas nécessairement injectif) $H^i ((X_n)_\syn, \J_n^{[t]})
\to \calM$ et on note $\Fil^t \calM$ son image. On vérifie que l'on
obtient ainsi une filtration admissible (voir définition
\ref{4:def:filadm}) sur $\calM$. Par définition (voir paragraphe
\ref{4:subsec:tstcov}), le membre de droite de l'isomorphisme 
\ref{4:eq:isom} s'identifie à :
$$\Fil^r (\Ast \otimes_S \calM)^{\phi_r=1}_{N=0}$$
avec :
$$\Fil^r (\Ast \otimes_S \calM) = \sum_{t=0}^r \Fil_X^t \Ast \otimes_S
\Fil^{s-t}\calM$$
où on rappelle que les $\Fil^t_X$ sont définis par la formule
(\ref{4:eq:filx}) et qu'ils sont plats sur $S_n$. On rappelle également
que l'on dispose du lemme \ref{4:lem:filconoyau} qui permet de voir
le module $\Fil^r (\Ast \otimes_S \calM)$ comme le conoyau d'un
morphisme.

\bigskip

La preuve de la proposition \ref{4:prop:isom} passe par les calculs
successifs des modules $\Fil^r (\Ast \otimes_S \calM)$, $\Fil^r (\Ast
\otimes_S \calM)_{N=0}$ et finalement $\Fil^r (\Ast \otimes_S
\calM)_{N=0}^{\phi_r=1} = {\Tst}_\star (\calM)$. Ceux-ci sont traités
respectivement dans les paragraphes \ref{4:subsec:preuve1},
\ref{4:subsec:preuve2} et \ref{4:subsec:preuve3}. Le paragraphe
\ref{4:subsec:jstar}, quant à lui, rappelle quelques préliminaires
nécessaires pour la gestion des limites inductives.

\subsubsection{Le foncteur $j_\star$}
\label{4:subsec:jstar}

Dans ce paragraphe, on rappelle comment construire des faisceaux sur le
site $(X_n)_\syn$ dont la cohomologie s'identifie à $\varinjlim_L H^i
((X_{n,L})_\syn, \J_n^{[s]})$.

\medskip

Si $L$ est une extension de $K$, on a un morphisme canonique $j_L :
T_L \to T$. On montre (lemme 3.1.1.1 de \cite{breuil-duke}) qu'il est
log-syntomique et donc qu'il induit un morphisme de topoï $(\widetilde
{T_{n,L}})_\syn \to (\widetilde {T_n})_\syn$. Soit $\calF$ un faisceau
de groupes abéliens sur $(T_n)_\syn$. Pour tout $L$, on considère le
faisceau $j_{L\star} j_L^\star \calF$ et on remarque que si $L'$ est une
extension finie de $L$, on a un morphisme $j_{L\star} j_L^\star \calF
\to j_{L'\star} j_{L'}^\star \calF$. On pose finalement :
$$j_\star \calF = \varinjlim_L j_{L\star} j_L^\star \calF.$$
C'est un faisceau sur $(T_n)_\syn$ et on montre (corollaire 3.1.1.4 de
\cite{breuil-duke}) que l'on a une identification canonique :
$$\varinjlim_L H^i((X_{n,L})_\syn, \J_n^{[s]}) = H^i((X_n)_\syn, j_\star
\J_n^{[s]})$$
pour tout entier $s$.

\subsubsection{Le calcul de $\Fil^r (\Ast \otimes_S \calM)$}
\label{4:subsec:preuve1}

Le but de ce paragraphe est de donner une description en terme de 
conoyau (analogue à celle du lemme \ref{4:lem:filconoyau}) de $\Fil^r
(\Ast \otimes_S \calM)$. Pour cela, on commence par rappeler que
l'anneau $\Ast$ admet une interprétation cohomologique incarnée par
l'isomorphisme canonique suivant : 
$$\Ast/p^n = \varinjlim_L H^0((T_{n,L})_\syn, \O_n^\st)$$
où la limite inductive est prise sur les extensions finies $L$ de $K$.
Il existe aussi un isomorphisme analogue pour décrire la filtration sur
$\Ast$ qui est :
$$\Fil^t (\Ast/p^n) = \Fil^t \Ast/p^n = \varinjlim_L H^0((T_{n,L})_\syn,
\J_n^{[t]}).$$
Ces isomorphismes permettent de construire une application canonique :
$$\Fil^t \Ast/p^n \otimes_{S_n} H^0((X_n)_\syn,
\J_n^{s-t}) \to \varinjlim_L H^0((X_{n,L})_\syn, \J_n^{[s]}) = 
H^0((X_n)_\syn, j_\star \J_n^{[s]})$$
et donc un morphisme de faisceaux :
$$\bigoplus_{t=0}^s \Fil_X^t \Ast/p^n \otimes_{S_n} \J_n^{[s-t]} \to
\bigoplus_{t=0}^s \Fil^t \Ast/p^n \otimes_{S_n} \J_n^{[s-t]} \to
j_\star \J_n^{[s]}.$$

\medskip

On a alors le lemme suivant, à mettre en parallèle avec le lemme
\ref{4:lem:filconoyau} :

\begin{lemme}
\label{4:lem:preuve1}
Pour tout entier $s \leq r$, on a un diagramme commutatif :
$${\tiny \xymatrix @R=10pt @C=20pt {
0 \ar[r] & \displaystyle \xybigoplus {t=1} s \Fil^t_X \xyAst/p^n
\otimes_{S_n} H^i(\J_n^{[s+1-t]}) \ar[d] \ar[r] & \displaystyle 
\xybigoplus {t=0} s \Fil^t_X \xyAst/p^n \otimes_{S_n} H^i(\J_n^{[s-t]})
\ar[r] \ar[d] & H^i (j_\star \J_n^{[s]}) \ar[r] \ar[d] & 0 \\
0 \ar[r] & \displaystyle \xybigoplus {t=1} s \Fil^t_X \xyAst/p^n
\otimes_{S_n} \Fil^{s+1-t} \calM \ar[r] & \displaystyle 
\xybigoplus {t=0} s \Fil^t_X \xyAst/p^n \otimes_{S_n} \Fil^{s-t} \calM 
\ar[r] & \Fil^s(\xyAst \otimes_{S} \calM) \ar[r] & 0 }}$$
où tous les morphismes respectent l'action de Galois, et où les deux 
lignes sont exactes et les flèches verticales surjectives. (Notez 
que tous les groupes de cohomologie sont calculés sur le site 
$(X_n)_\syn$.)
\end{lemme}

\begin{preuve}
Tout d'abord, précisons les flèches. Dans la suite exacte
du haut, la première flèche a déjà été définie dans l'énoncé du
lemme \ref{4:lem:filconoyau}. La flèche correspondante dans la suite
exacte du bas a une définition tout à fait analogue. Les autres
flèches ne posent pas de problème, à part \emph{a priori} la flèche
verticale de droite. Cependant, elle n'en posera plus lorsque l'on aura
prouvé l'exactitude des deux suites (puisque ce sera alors simplement
la flèche induite sur les conoyaux).

L'exactitude de la ligne du haut n'est autre que l'objet du lemme 
\ref{4:lem:filconoyau}. Les surjectivités des deux premières flèches
verticales sont immédiates. Il ne reste donc plus qu'à prouver
l'exactitude de la suite exacte du bas (de laquelle résultera
directement la surjectivité de la flèche verticale de droite).

\medskip

On commence par prouver que la suite de faisceaux sur le site
$(T_n)_\syn$ :
$$\xymatrix @C=20pt {
0 \ar[r] & \displaystyle \xybigoplus {t=1} s \Fil^t_X \xyAst /p^n
\otimes_{S_n} \J_n^{[s+1-t]} \ar[r] & \displaystyle 
\xybigoplus {t=0} s \Fil^t_X \xyAst/p^n \otimes_{S_n} \J_n^{[s-t]}
\ar[r] & j_\star \J_n^{[s]} \ar[r] & 0 }$$
est exacte. Par un dévissage, on se ramène dans un premier temps au
seul cas $n=1$. De plus, en recopiant les arguments de la preuve de la
proposition 3.1.2.3 de \cite{breuil-duke}, on se ramène au cas $s=0$.
Il s'agit donc de montrer que $\Ast/p^n \otimes_{S_n} \O^1_\st \simeq
j_\star \O_1^\st$. C'est à nouveau un calcul local pour la topologie
log-syntomique, en tout point analogue à celui mené dans la
démonstration du lemme 3.1.2.2 de \cite{breuil-duke}. 
\end{preuve}

\subsubsection{Le calcul de $\Fil^r (\Ast \otimes_S \calM)_{N=0}$}
\label{4:subsec:preuve2}

Le but de cette partie est de démontrer le lemme suivant qui
constitue la deuxième étape de la preuve.

\begin{lemme}
\label{4:lem:preuve2}
On a des isomorphismes de modules galoisiens :
$$\xymatrix @R=8pt @C=35pt {
\displaystyle \xyvarinjlim L H^i((X_{n,L})_\syn, \O_n^\cris) 
\ar[r]^-{\sim} & (\xyAst \otimes_S \calM)_{N=0} \\
\displaystyle \xyvarinjlim L H^i((X_{n,L})_\syn, \J_n^{\cris,[r]}) 
\ar[r]^-{\sim} & \Fil^r (\xyAst \otimes_S \calM)_{N=0} }$$
\end{lemme}

\begin{preuve}
La démonstration est identique à celle du corollaire 3.2.3.5 de
\cite{breuil-duke}. Nous redonnons simplement les grandes lignes.
Tout d'abord, on montre comme dans le lemme 3.1.1.2 de
\cite{breuil-duke} la nullité de $\varinjlim_L R^i j_{L\star} j_L^\star
\J_n^{[s]}$ d'où on déduit, à partir de la première suite exacte de
la proposition \ref{4:prop:sns}, une suite exacte de faisceaux :
$$\xymatrix {
0 \ar[r] & j_\star \J_n^{\cris,[s]} \ar[r] & j_\star \J_n^{[s]}
\ar[r]^-{N} & j_\star \J_n^{[s-1]} \ar[r] & 0 }$$
de laquelle on déduit une suite exacte courte sur les groupes de
cohomologie :
$$\xymatrix {
0 \ar[r] & H^i(j_\star \J_n^{\cris,[s]}) \ar[r] &
H^i(j_\star \J_n^{[s]}) \ar[r]^-{N} & H^i(j_\star \J_n^{[s-1]}) \ar[r] &
0 }$$
où tous les groupes de cohomologie sont calculés sur le site
$(X_n)_\syn$ (l'argument est le même que celui utilisé pour la
proposition 3.2.3.1 de \cite{breuil-duke}).

\medskip

La suite de la preuve consiste à reprendre le diagramme du lemme
\ref{4:lem:preuve1} et à procéder à une étude relativement fine des noyaux
des flèches verticales. Précisément, si on note :
\begin{eqnarray*}
K^t & = & \ker (H^i(\J_n^{[t]}) \to H^i(\O_n^\st)) = \ker
(H^i(\J_n^{[t]}) \to \Fil^t \calM) \\
\bar K^t & = & \ker (H^i(j_\star \J_n^{[t]}) \to \Fil^t (\Ast
\otimes_S \calM))
\end{eqnarray*}
on peut compléter le diagramme de la façon suivante :
$${\tiny \xymatrix @R=10pt @C=20pt {
& 0 \ar[d] & 0 \ar[d] & 0 \ar[d] \\
0 \ar[r] & \displaystyle \xybigoplus {t=1} s \Fil^t_X \xyAst/p^n
\otimes_{S_n} K^{s+1-t} \ar[d] \ar[r] & \displaystyle 
\xybigoplus {t=0} s \Fil^t_X \xyAst/p^n \otimes_{S_n} K^{s-t}
\ar[r] \ar[d] & \xybar K^s \ar[r] \ar[d] & 0 \\
0 \ar[r] & \displaystyle \xybigoplus {t=1} s \Fil^t_X \xyAst/p^n
\otimes_{S_n} H^i(\J_n^{[s+1-t]}) \ar[d] \ar[r] & \displaystyle 
\xybigoplus {t=0} s \Fil^t_X \xyAst/p^n \otimes_{S_n} H^i(\J_n^{[s-t]})
\ar[r] \ar[d] & H^i (j_\star \J_n^{[s]}) \ar[r] \ar[d] & 0 \\
0 \ar[r] & \displaystyle \xybigoplus {t=1} s \Fil^t_X \xyAst/p^n
\otimes_{S_n} \Fil^{s+1-t} \calM \ar[r] \ar[d] & \displaystyle 
\xybigoplus {t=0} s \Fil^t_X \xyAst/p^n \otimes_{S_n} \Fil^{s-t} \calM 
\ar[r] \ar[d] & \Fil^s(\xyAst \otimes_{S} \calM) \ar[r] \ar[d] & 0 \\
& 0 & 0 & 0 }}$$
L'opérateur $N$ induit un morphisme entre le diagramme précédent et
son équivalent lorsque l'on remplace $s$ par $s-1$ (en particulier, il
induit une application $N : \bar K_s \to \bar K_{s-1}$ pour tout $s \in
\{0, \ldots, r\}$ en convenant que $K_{-1} = K_0$). La suite exacte du
haut implique en prenant $s=0$ que $\bar K_0 = 0$. Par ailleurs, une
étude un peu minutieuse de cette même suite exacte (voir lemmes 3.2.3.3
et 3.2.3.4 de \cite{breuil-duke}) montre que l'application $N : \bar K_r
\to \bar K_{r-1}$ est un isomorphisme.

En considérant la suite exacte de droite, ceci implique que pour $s=0$
et $s=r$, on a des isomorphismes :
$$\xymatrix @C=35pt {
\displaystyle \xyvarinjlim L H^i((X_{n,L})_\syn, \J_n^{\cris,[s]}) =
H^i((X_n)_\syn, j_\star \J_n^{\cris,[s]}) \ar[r]^-{\sim} &
\Fil^s (\xyAst \otimes_S \calM)_{N=0} }$$
ce qui termine la preuve du lemme.
\end{preuve}

\subsubsection{Le calcul de $\Fil^r (\Ast \otimes_S
\calM)_{N=0}^{\phi_r=1}$}
\label{4:subsec:preuve3}

Il n'est plus difficile à présent de terminer la preuve de la
proposition \ref{4:prop:isom}. En effet, la deuxième suite exacte de la
proposition \ref{4:prop:sns} nous fournit une suite exacte longue :
$${\small
\xymatrix { \cdots \ar[r] & H^i((X_{n+r})_\syn, \S_n^r) \ar[r] &
H^i((X_n)_\syn, \J_n^{\cris,[r]}) \ar[r]^-{\phi_r - \id} \ar[r] &
H^i((X_n)_\syn, \O_n^\cris) \ar[r] & \cdots }}$$
et puis, comme le foncteur $\varinjlim_L$ est exact (la limite est
filtrante) on obtient une nouvelle suite exacte longue :
$$\xymatrix { \cdots \ar[r] & \displaystyle \xyvarinjlim L
H^i((X_{n+r})_\syn, \S_n^r) \ar[r] & \displaystyle \xyvarinjlim L
H^i((X_n)_\syn, \J_n^{\cris,[r]}) \ar[r]^-{\phi_r - \id} \ar[r] &
\displaystyle \xyvarinjlim L H^i((X_n)_\syn, \O_n^\cris) \ar[r] & \cdots
} .$$
Par ailleurs, la flèche $\phi_r - \id : \varinjlim_L H^i((X_n)_\syn,
\J_n^{\cris,[r]}) \to \varinjlim_L H^i((X_n)_\syn, \O_n^\cris)$
s'identifie \emph{via} les isomorphismes du lemme \ref{4:lem:preuve2} à la
flèche :
$$\xymatrix @C=50pt {
\Fil^r (\xyAst \otimes_S \calM)_{N=0} \ar[r]^-{\phi_r-\id} & 
(\xyAst \otimes_S \calM)_{N=0} }$$
et on sait par le lemme \ref{4:lem:phirsurj} que celle-ci est surjective.
On en déduit que la suite exacte longue se coupe en suites exactes
courtes :
$$\xymatrix { 0 \ar[r] & \displaystyle \xyvarinjlim L H^i((X_{n+r})_\syn,
\S_n^r) \ar[r] & \displaystyle \xyvarinjlim L H^i((X_n)_\syn,
\J_n^{\cris,[r]}) \ar[r]^-{\phi_r - \id} \ar[r] & \displaystyle
\xyvarinjlim L H^i((X_n)_\syn, \O_n^\cris) \ar[r] & 0 }$$
ce qui termine la preuve.

\subsection{Une conjecture de Serre}
\label{4:subsec:serre}
On montre dans ce paragraphe comment la théorie développée au long de
cet article permet de résoudre complètement la conjecture de l'inertie
modérée de Serre formulée dans le paragraphe 1.13 de \cite{serre}.

\bigskip

Avant de rappeler l'énoncé de la conjecture, faisons quelques 
préliminaires et profitons-en pour fixer les notations (pour plus de 
précisions, voir le paragraphe 1 de \cite{serre}). Soit $V$ une
$\F_p$-représentation de dimension finie irréductible du sous-groupe
d'inertie $I$ de groupe de Galois absolu de $K$. Par un résultat 
classique de théorie des groupes, du fait que $V$ a pour cardinal 
un multiple de $p$, le sous-groupe d'inertie sauvage (qui est un 
pro-$p$-groupe distingué) agit trivialement. Ainsi l'action de $I$ se 
factorise à travers une action du groupe d'inertie modérée $I_t$.

Par ailleurs, puisque $V$ est supposée irréductible, l'anneau $E$ des
endomorphismes équivariants de $V$ est un corps fini et $V$ hérite
d'une structure d'espace vectoriel de dimension $1$ sur ce corps. La
représentation de départ fournit un caractère $\rho : \It \to 
E^\star$. Notons $q = p^h$ le cardinal de $E$ et $\F_q$ le sous-corps de
$\bar k$ formé des solutions de l'équation $x^q = x$. On dispose 
de l'application suivante appelée \emph{caractère fondamental de
niveau $h$} :
$$\begin{array}{rcl}
\theta_h : \quad \It & \to & \mu_{q-1}(\bar K) \simeq \F_q^\star \\
g & \mapsto & \frac{g(\eta)}{\eta}
\end{array}$$
où $\eta$ désigne une racine $(q-1)$-ième de l'uniformisante $\pi$. 

Les corps $E$ et $\F_q$ sont finis de même cardinal et donc isomorphes 
(non canoniquement). Si l'on compose $\theta_h$ par
un tel isomorphisme $f$, on obtient une application $\theta_{h,f} : \It
\to E^\star$ et on montre facilement (voir la proposition 5 du
paragraphe 1 de \cite{serre}), que $\rho = \theta_{h,f}^n$ pour un
certain entier $n$ compris entre $0$ et $q-2$. L'entier $n$ dépend de
l'isomorphisme $f$ choisi mais les chiffres de son écriture en base $p$,
eux, n'en dépendent pas. Ce sont par définition les \emph{exposants de
l'inertie modérée} de la représentation $V$.

\medskip

La conjecture de Serre s'énonce alors comme suit :

\begin{theo}
Soit $X$ un schéma propre et lisse sur $K$ à réduction semi-stable 
sur $\O_K$ et soit $r$ un entier. Les exposants de l'inertie
modérée sur un quotient de Jordan-Hölder de la restriction au groupe
d'inertie de la représentation galoisienne $H^r_\et\pa{X_{\bar K},
\Z/pZ}^\vee$ (où \og $\vee$ \fg\ désigne le dual) sont compris entre
$0$ et $er$.
\end{theo}

\begin{preuve}
On remarque dans un premier temps que le résultat est évident si $er
\geq p-1$ (des chiffres en base $p$ sont nécessairement inférieurs ou
égaux à $p-1$). On peut donc supposer $er < p-1$ et appliquer les
résultats de cet article.

Par le théorème \ref{4:th:comp}, la représentation galoisienne
$V = H^r_\et\pa{X_{\bar K}, \Z/pZ}$ est dans l'image essentielle du
foncteur ${\Tst}_\star$. Puisque cette image essentielle est stable par
sous-objets et quotients (théorème \ref{4:th:tstcont}), tout quotient de
Jordan-Hölder de $V$ est également dans l'image essentielle de
${\Tst}_\star$. Par ailleurs, un tel quotient de Jordan-Hölder est par
définition irréductible et donc ne peut être l'image par
${\Tst}_\star$ que d'un objet simple.

Le théorème résulte à ce niveau du théorème 5.2.2 de \cite{caruso}.
\end{preuve}

\bigskip

\noindent
{\it Remarque.} On a un résultat équivalent avec les groupes de
cohomologie $H^i_\et\pa{X_{\bar K}, \Z/p^nZ}^\vee$ lorsque $i < r$ et
$er < p-1$ (ou autrement dit lorsque $i < E((p-2)/e)$, $E(\cdot)$
désignant la partie entière).

\nocite{fontaine-messing}
\nocite{fontaine-ast1}
\nocite{fontaine-ast2}
\nocite{fontaine-laffaille}
\nocite{breuil-griffiths}
\nocite{breuil-messing}
\nocite{tsuji}
\nocite{faltings}
\nocite{faltings2}
\nocite{berthelot-ogus}
\bibliography{conjserre}
\bibliographystyle{amsalpha}

\end{document}